\documentclass[11pt]{article}
\usepackage{amssymb,amsmath,amsthm,amsfonts,mathrsfs}
\usepackage[usenames]{color}
\usepackage{indentfirst}
\usepackage{mathtools}
\usepackage{cite}
\usepackage{verbatim}
\topmargin by-\headsep  \textheight 246mm
\usepackage[active]{srcltx}
\usepackage{hyperref}
\usepackage[margin=1in]{geometry}
\usepackage[ruled,vlined]{algorithm2e}
\usepackage{algorithmicx,algpseudocode}
\algdef{SE}[DOWHILE]{Do}{doWhile}{\algorithmicdo}[1]{\algorithmicwhile\ #1}%
\usepackage{tikz}
\usepackage{graphicx}
\usepackage{float}
\usepackage[justification=raggedright]{caption}
\usepackage{subcaption}
\usepackage{wrapfig}
\usepackage{diagbox}
\usepackage{bm}
\usepackage{multirow,multicol}
\usepackage{comment}
\usepackage{booktabs}
\usepackage{threeparttable}
\graphicspath{{./pics/}}

\numberwithin{equation}{section}
\newtheorem{thm}{Theorem}[section]

\newtheorem{lem}{Lemma}[section]

\theoremstyle{remark}
\newtheorem{rem}{Remark}[section]

\newtheorem{example}{Example}[section]

\newcommand{\dx}{{\rm d}x}
\newcommand{\ds}{\mathrm{ds}}
 \def\p{\partial} 
\def\to{\rightarrow}
\def\Om{\Omega}  
    \def\R{{\mathbb R}} \def\T{\mathcal{T}}

\def\m{\mbox} \def\t{\times}

\def\bold{\boldsymbol}

\def\bb{\begin{equation}} \def\ee{\end{equation}}

\def\beqn{\begin{eqnarray}}  \def\eqn{\end{eqnarray}}

\def\beqnx{\begin{eqnarray*}} \def\eqnx{\end{eqnarray*}}

\begin{document}
\begin{center}
{\Large\bf Adaptive finite element approximations of the first eigenpair associated with $p$-Laplacian}
\end{center}

\centerline{Guanglian Li\footnote{Department of Mathematics, The University of Hong Kong, Hong Kong Special Administrative Region, China. ({\tt lotusli@maths.hku.hk})}\quad Jing Li\footnote{School of Mathematical Sciences, East China Normal University, Shanghai 200241, China. ({\tt betterljing@163.com})}\quad Julie Merten\footnote{Computational and Numerical Mathematics, Bernoulli Institute, University of Groningen, the Netherlands. ({\tt j.y.merten@rug.nl})}\quad Yifeng Xu\footnote{Corresponding author. Department of Mathematics \& Scientific Computing Key Laboratory of Shanghai Universities, Shanghai Normal University, Shanghai 200234, China. ({\tt
yfxu@shnu.edu.cn})}\quad Shengfeng Zhu\footnote{Key Laboratory of MEA (Ministry of Education) \& Shanghai Key Laboratory of PMMP, School of Mathematical Sciences, East China Normal University, Shanghai 200241, China. ({\tt sfzhu@math.ecnu.edu.cn})}}

\begin{abstract}
In this paper, we propose an adaptive finite element method for computing the first
eigenpair of the $p$-Laplacian problem. We prove that starting from a fine initial mesh our proposed adaptive algorithm produces a sequence of discrete first eigenvalues that converges to the first eigenvalue of the continuous problem and the distance between discrete eigenfunctions and the normalized eigenfunction set corresponding to the first eigenvalue in $W^{1,p}$-norm also tends to zero. Extensive numerical examples are provided to show the effectiveness and efficiency.
\end{abstract}


\vspace{0.2cm}


\section{Introduction}
In this paper, we consider the eigenvalue problem of the $p$-Laplacian operator with homogeneous Dirichlet boundary condition:
\begin{equation}\label{diff_eq}
    \left\{\begin{array}{ll}
        -\boldsymbol{\nabla}\cdot(|\boldsymbol{\nabla}u|^{p-2}\boldsymbol{\nabla}u)=\lambda|u|^{p-2}u\quad&\m{in}~\Om,\\
        u=0\quad&\m{on}~\p\Om,
    \end{array}
    \right.
\end{equation}
where $\Omega$ is an open bounded Lipschitz polygonal/polyhedral connected domain in $\mathbb{R}^d$ ($d=2,3$) and $1<p<\infty$.
From the perspective of applications, the $p$-Laplacian operator arises from non-Newtonian
fluids \cite{GlowinskiRappaz:2003}  and power-law materials \cite{AtkinsonChampion:1984}. As usual, integration by parts yields the weak formulation of \eqref{diff_eq}: Find $(\lambda,u)\in
\mathbb{R}\times V$ such that
\begin{equation}\label{vp_eigen}
    \int_\Om|\boldsymbol{\nabla}u|^{p-2}\boldsymbol{\nabla}u\cdot \boldsymbol{\nabla}v\dx=\lambda\int_\Om|u|^{p-2}uv\dx\quad\forall~v\in V.
\end{equation}
Here we denote $V:=W_{0}^{1,p}(\Om)$.

The existing theory developed in \cite{gp,le,lind}
asserts that Problem \eqref{vp_eigen} has a nondecreasing sequence of positive eigenvalues
$\{\lambda_n\}_{n\geq1}$ diverging to $+\infty$. It should be noted that the first eigenvalue $\lambda_1$
is simple and isolated \cite{le,lind}, and is equivalent to the minimum of the Rayleigh quotient
\begin{equation}\label{min_cont}
\lambda_1=\inf_{v\in V\setminus \{0\}}\mathcal{J}(v):=\frac{\int_\Om|\boldsymbol{\nabla}v|^p\dx}{\int_\Om|v|^p\dx}.
\end{equation}
The existence of a minimizer to Problem \eqref{min_cont} can be established by the
standard minimization approach (cf. \cite{BS,lind}). Moreover, the reciprocal of $\lambda_1$ is the best constant in Poincar\'e inequality, which implies that $\lambda_1>0$. A normalized
eigenfunction set with respect to $\lambda_1$ is defined as
\[E_{\lambda_1}:=\bigg\{u\in
V\bigg\vert\int_\Om|\boldsymbol{\nabla}u|^{p-2}\boldsymbol{\nabla}u\cdot
\boldsymbol{\nabla}v\dx=\lambda_1\int_\Om|u|^{p-2}uv\dx,~\forall~v\in
V,~\|u\|_{L^p(\Omega)}=1\bigg\}.\]

Some attempts have been made in numerically computing eigenpairs of Problem \eqref{diff_eq}  \cite{BBEM,BEM,BoSz,Horak:2011,LeftonWei:1997,YaoZhou:2007}. But due to the degenerate structure of the operator and the existence of possible reentrant
corners in the computational domain $\Om$, the solution to Problem \eqref{diff_eq} features local singularities.
As a remedy, adaptive techniques are preferred in numerical simulation for accuracy and
efficiency. Generally speaking, a standard adaptive finite element method (AFEM) comprises the following four modules in every loop:
\begin{equation}\label{afemloop}
    \textrm{SOLVE}\to\textrm{ESTIMATE}\to\textrm{MARK}\to\textrm{REFINE}.
\end{equation}
The most prominent advantage of AFEM is to make efficient
use of computer resources to attain the given error tolerance with minimum degrees of freedom, so it has become an effective tool in practice of scientific computing and
engineering. Since the seminal work \cite{br} by Babu\v{s}ka and Rheinboldt in 1978, there
has been much and rapid progress in the mathematical theory of this field. In particular,
the understanding of a posteriori estimation, the main ingredient in the module ESTIMATE, is now on a mature level; see e.g. \cite{ao,ver}. Moreover, great efforts have been put
into the study of AFEM itself in terms of convergence and complexity over the past three decades.
For linear elliptic problems one may refer to two survey papers \cite{cfpp,nsv} and
the references therein for an overview. In the case of linear or nonlinear eigenvalue
problems, where nonlinearity consists in low order terms,
we are aware that existing works are only limited to the linear Laplacian/diffusion/bi-Laplacian operator; see  \cite{BoffiGallistlGardiniGastaldi:2017,BonitoDemlow:2016,CarstensenGallistlSchedensack:2015,cg1,cg2,CarPut:2022,ChenDaiGongHeZhou:2014,cghz,chz,DaiHeZhou:2015,dxz,Gallistl:2014,Gallistl:2015a,Gallistl:2015b,gm,gmz1}. For the nonlinear Laplacian
equation, we mention \cite{bdk, CarLiuYan:2006, dk, LiuChen:2020,ly,veeser} for results on a posteriori error estimation and adaptive computations.

The aim of this paper is to develop adaptive finite element approximations of the first eigenvalue $\lambda_1$ to Problem
\eqref{diff_eq}. To be specific, we propose Algorithm \ref{afem} of standard form \eqref{afemloop}, which facilitates implementation in practical applications, for the first eigenpair of Problem \eqref{vp_eigen} in Section \ref{sec:AFEM} and establish the convergence of its resulting first
discrete eigenpairs $\{(\mu_k,u_k)\}_{k\geq0}$ in Section \ref{sec:conv}. It is demonstrated in Theorem \ref{thm_afem_conv} that the whole sequence $\{\mu_k\}_{k\geq0}$ converges to $\lambda_1$ and the $W^{1,p}(\Om)$-norm between $L^p(\Om)$-normalized sequence $\{u_k\}_{k\geq0}$ and $E_{\lambda_1}$ tends to
zero. 

In addition to standard arguments for linear problems \cite{gm,gmz1},
minimization techniques for nonlinear elliptic problems \cite{am,BS} are utilized to deal with the nonlinear structure of Problem \eqref{diff_eq} in the convergence analysis. By introducing an auxiliary minimization problem \eqref{min_lim} over the limiting space given by
the adaptive process \eqref{afemloop}, we first prove in Theorem \ref{thm_min_lim} and Theorem \ref{thm_min_conv} that the sequence of discrete
eigenpairs $\{(\mu_k,u_k)\}_{k\geq0}$ converges (up to a subsequence) to $(\mu_\infty,u_\infty)$, where $\mu_\infty$ denotes the minimum
to Problem \eqref{min_lim}, and $u_\infty$ the $L^p(\Om)$-normalized minimizer. Then we further prove $(\mu_\infty,u_\infty)$ satisfies the variational formulation \eqref{vp_eigen} (see Section \ref{subsect:aux} and Step 1 in the proof of Theorem \ref{thm_afem_conv}), which indicates that $u_\infty$ is a critical point of $\mathcal{J}$ over $V$. Finally, the
convergence of the whole sequence $\{(\mu_k,u_k)\}_{k\geq0}$ to $\lambda_1$ and $E_{\lambda_1}$ is proved under the assumption that the initial mesh is sufficiently fine (see Step 2 and Step 3 in the proof of Theorem \ref{thm_afem_conv}). The unquantifiable fineness requirement on the initial mesh is precisely undesirable in adaptive computations, but it seems inevitable even in the analysis of AFEM for linear eigenvalue problems. It should also be pointed out that some computable quantities adopted in the module
ESTIMATE are derived (see the proof of Lemma \ref{lem_res_conv}), although they do not provide an upper bound of the error, in our convergence analysis, where a practical assumption (see \eqref{marking_max} in Algorithm \ref{afem}) imposed in the module MARK as for linear cases \cite{gmz1, gm} is utilized.

The remaining of this paper is organized as follows. A numerical scheme built on the finite
element method for Problem \eqref{min_cont} is presented in Section 2. In Section 3, we
introduce a standard adaptive finite element method with a general yet reasonable requirement on the marking strategy, the convergence of which is investigated in Section 4. Section 5 deals with the implementation of our proposed algorithm and contains some numerical results illustrating the efficiency. The paper is ended with some concluding remarks in Section 6. Throughout the paper, we use standard
notation for $L^{p}(\Omega)$ or $L^\infty(\Omega)$ space, the  Sobolev space $W_0^{1,p}(\Omega)$ ($W^{2,\infty}(\Omega)$)  and its dual space $W^{-1,q}(\Om)$ with $q=p/(p-1)\in (1,\infty)$ as well as their related (semi-)norms. Moreover, the upper-case letter $C$, with or
without subscript, denotes a generic constant independent of the mesh size and it may take a different value at each occurrence.

\section{Discrete Problem}
In this section, we introduce a discrete problem to approximate the minimization problem \eqref{min_cont}. For this purpose, let $\mathcal{T}$ be a conforming
triangulation of $\overline{\Om}$ into a set of closed triangles or tetrahedra with a discretization parameter
$h_T:=|T|^{1/d}$ for each $T\in\mathcal{T}$.  Let $V_\mathcal{T}$ be the associated conforming space of continuous piecewise linear functions vanishing on boundary $\p\Om$ given by
\[
    V_\mathcal{T}:=\{v\in C(\overline{\Omega})~|~v|_{\p\Om}=0,v|_{T}\in P_1(T),~\forall~T\in V_{\mathcal{T}}\}.
\]
Then the finite element approximation of \eqref{min_cont} is seeking $u_\mathcal{T}\in V_\mathcal{T}\setminus\{0\}$, satisfying
\begin{equation}\label{min_disc}
    \mu_\T:=\mathcal{J}(u_\mathcal{T})=\inf_{v\in V_\mathcal{T}\setminus \{0\}}\mathcal{J}(v).
\end{equation}
Note that the property $V_\T\subset V$ implies that
\begin{equation}\label{eigva_lb}
    0<\lambda_1\leq\mu_\T.
\end{equation}
\begin{thm}\label{thm_min_disc}
Let $\mu_{\mathcal{T}}$ be the solution to Problem \eqref{min_disc}, then $\mu_{\mathcal{T}}$ is positive and attained by some nonnegative function
$u_\mathcal{T}\in V_\mathcal{T}\setminus \{0\}$.
\end{thm}

\begin{proof}
Our proof follows from the argument in \cite{BS}, which is concerned with the continuous problem \eqref{min_cont}. The Poincar\'{e} inequality reads
\[
     \frac{\int_\Om|\boldsymbol{\nabla}u_\T|^p\dx}{\int_\Om|u_\T|^p\dx}\geq C>0\quad\forall~v\in
     V_\mathcal{T}\setminus\{0\},
\]
which provides a positive lower bound. Hence, $\mu_\T$ is positive.

Let $\{v_m\}_{m\geq 0}\subset
V_\mathcal{T}\setminus\{0\}$ be a minimizing sequence to Problem \eqref{min_disc}. Since $\{|v_m|\}_{m\geq0}$ is also a
minimizing sequence, then we can assume that $v_m\geq 0$ a.e. in $\Om$ for all $m$. Moreover, the homogeneity of the objective functional $\mathcal{J}$ allows for the normalization
$\int_\Om|v_m|^p\dx=1$. Consequently, the minimizing sequence $\{|v_m|\}_{m\geq0}$ is bounded in the finite dimensional space
$V_\mathcal{T}$. This guarantees the existence of a subsequence, still denoted by $\{v_m\}_{m\geq0}$, and some $u_{\mathcal{T}}\in V_\mathcal{T}$, satisfying
\begin{equation*}
    v_m\to u_\T~\mbox{{strongly} in}~W_0^{1,p}(\Om),\quad v_m\to u_\T~\mbox{a.e. in}~\Om.
\end{equation*}
By this, we derive that $u_\mathcal{T}\geq0$ a.e. in $\Om$, $\int_\Om|u_\T|^p\dx=1$ and $\mathcal{J}(u_\T)=\lim_ {m\to\infty}\mathcal{J}(v_m)=\mu_{\T}$. The proof is complete.
\end{proof}
Theorem \ref{thm_min_disc} implies that the discrete minimizer $u_\T$ can be normalized as $\|u_\T\|_{L^p(\Omega)}=1$ in our adaptive algorithm below and the subsequent analysis. 
By the differential calculus for $\mathcal{J}$ \cite{BS}, it is easy to see that solution
$(\mu_{\T},u_\T)$ to Problem \eqref{min_disc} satisfies the discrete formulation of Problem \eqref{vp_eigen},
\begin{equation}\label{vp_disc}
    \int_\Om|\boldsymbol{\nabla}u_\T|^{p-2}\boldsymbol{\nabla}u_\T\cdot\boldsymbol{\nabla}v\dx=\mu_\T\int_\Om|u_\T|^{p-2}u_\T v\dx\quad\forall~v\in V_\T.
\end{equation}

\section{Adaptive Finite Element Method}\label{sec:AFEM}

Let $\mathbb{T}$ be the set of all possible conforming triangulations of
$\overline{\Omega}$ obtained from some initial mesh $\mathcal{T}_0$ by successive use of bisection \cite{koss,nsv,stev1}. This refinement process ensures that the set $\mathbb{T}$ is uniformly shape regular, i.e. the
shape regularity of any $\mathcal{T}\in\mathbb{T}$ is uniformly bounded by a constant depending on the initial mesh $\mathcal{T}_0$
\cite{nsv,traxler}. 
$\mathcal{T}'\in \mathbb{T}$ is referred to as a refinement of
$\mathcal{T}\in\mathbb{T}$ if $\mathcal{T}'$ is
produced from $\mathcal{T}$ by a finite number of bisections. The collection of all interior faces in $\mathcal{T}\in\mathbb{T}$ is
denoted by $\mathcal{F}_\mathcal{T}(\Omega)$ and the scalar $h_F:=|F|^{1/(d-1)}$ stands for the diameter-equivalent mesh-size 
of each $F\in\mathcal{F}_{\mathcal{T}}(\Omega)$,
which is associated with a fixed normal unit vector $\boldsymbol{n}_F$.
The union of elements neighbouring some $T\in\mathcal{T}$ is
denoted by $D_{T}$, i.e.
\begin{align*}
D_{T}=\bigcup_{T'\in\mathcal{T}: \partial T \cap \partial T' \ne \emptyset }T'.
\end{align*}
First, let $(\mu_\T,u_\T)\in \R\t V_\T$ be the solution to Problem \eqref{min_disc} with $\|u_\T\|_{L^p(\Omega)}=1$ for $\mathcal{T}\in \mathbb{T}$.
We define an element residual and a jump residual with respect to an element $T\in\mathcal{T}$ and a face $F\in\mathcal{F}_\T(\Om)$ by
\[
    R_T(\mu_\T,u_\T):=\mu_\T |u_\T|^{p-2}u_\T,\quad J_F(u_\T):=[|\boldsymbol{\nabla}u_\T|^{p-2}\boldsymbol{\nabla}u_\T]\cdot\boldsymbol{n}_F
\]
with 
\[
    [|\boldsymbol{\nabla}u_\T|^{p-2}\boldsymbol{\nabla}u_\T]:=(|\boldsymbol{\nabla}u_\T|^{p-2}\boldsymbol{\nabla}u_\T)|_{F\subset \partial T^+}-(|\boldsymbol{\nabla}u_\T|^{p-2}\boldsymbol{\nabla}u_\T)|_{F\subset \partial T^-}
\]
denoting the jump across an interior face $F$ shared by $T^+, T^- \in \T$. Then the local error indicator on each element $T\in\mathcal{T}$ is defined by
\begin{align}\label{eq:err-ind}
    \eta_{\mathcal{T}}^{q}(\mu_\T,u_{\mathcal{T}};T):=h_T^{q}\|R_{T}(\mu_\T,u_{\mathcal{T}})\|_{L^{q}(T)}^{q}+\sum_{F\subset\partial T\cap\Omega}h_F\|J_F(u_{\mathcal{T}})\|^q_{L^q(F)}\quad \text{with}~q = p / (p-1).
\end{align}
Over some element patch $\mathcal{M}\subseteq\mathcal{T}$, the error estimator is defined by
\[
    \eta_{\mathcal{T}}(\mu_\T,u_{\mathcal{T}};\mathcal{M}):=\bigg(\sum_{T\in\mathcal{M}}\eta^q_{\mathcal{T}}(\mu_\T,u_{\mathcal{T}};T)\bigg)^{1/q}.
\]
When $\mathcal{M}=\mathcal{T}$, we abbreviate $\eta_{\mathcal{T}}(\mu_\T,u_{\mathcal{T}};\mathcal{T})$ to $ \eta_{\mathcal{T}}(\mu_\T,u_{\mathcal{T}})$.

Next, we propose an AFEM for Problem \eqref{min_cont}. In what follows, all dependence on a triangulation $\mathcal{T}_k$ is replaced by the mesh refinement level $k$ in the subscript, e.g. $V_k:=V_{\mathcal{T}_k}$.
\begin{algorithm}
\caption{AFEM for the 1st eigenvalue of $p$-Laplacian}\label{afem}
\begin{algorithmic}[1]
    \State {(INITIALIZE)} Specify an initial conforming mesh $\mathcal{T}_0$ and set counter $k:=0$.

    \State {(SOLVE)} Solve Problem \eqref{min_disc} on $\mathcal{T}_k$ for $(\mu_{k},u_{k})\in \mathbb{R}\times V_{k}$ s.t., $\|u_k\|_{L^p(\Omega)}=1$.

    \State {(ESTIMATE)} Compute the error estimator $\eta_{k}(\mu_k,u_k; T)$ by \eqref{eq:err-ind} for each $T\in\mathcal{T}_k$.

    \State {(MARK)} Mark a subset $\mathcal{M}_k\subseteq\mathcal{T}_k$ such that $\mathcal{M}_k$ contains at least one element $T_k$ with the largest error indicator over $\mathcal{T}_k$, i.e.
                \begin{equation}\label{marking_max}
                \eta_{k}(\mu_k,u_k;T_k):=\max_{T\in\mathcal{T}_k}\eta_{k}(\mu_k,u_k;T).
                \end{equation}

    \State {(REFINE)} Refine each $T\in\mathcal{M}_k$ by bisection to get $\mathcal{T}_{k+1}$.

    \State Set $k:=k+1$ and go to Step 1.

    \end{algorithmic}
\end{algorithm}

We note that a general yet reasonable assumption is included in the module MARK of Algorithm \ref{afem}. This requirement in marking elements for a further refinement helps to improve the computing efficiency in practical applications and is fulfilled by several popular marking strategies \cite{nsv}, e.g. the maximum strategy, the equi-distribution strategy, the modified equi-distribution strategy and the practical D\"{o}rfler's strategy. In the numerical implementation of Algorithm \ref{afem}, a tolerance for the estimator or a bound for the number of DOFs is usually prescribed as a stopping criterion. Without it, an infinite sequence $\{\mu_k\}_{k\geq0}$ is generated by Algorithm \ref{afem} and a natural question is whether it converges to the first eigenvalue of \eqref{diff_eq}, which will be examined at extensive length in the next section. For this purpose, we end this section with a stability estimate for the local error indicator.
\begin{lem}\label{lem_stab_err}
    Let $\{(\mu_k,u_k)\}_{k\geq0}$ be the sequence of discrete solutions generated by Algorithm \ref{afem}.
    For any $T\in\mathcal{T}_k$, there holds
    \begin{equation}\label{stab_err}
\begin{aligned}
        \eta_k^q(u_k,\mu_k;T)&\leq C\left( h_T^q\|u_k\|_{L^{p}(T)}^{p}+\|\boldsymbol{\nabla}u_k\|_{L^{p}(D_T)}^{p}\right),\\
        \eta_k(u_k,\mu_k)&\leq \widetilde{C},
\end{aligned}
    \end{equation}
    where the constants $C$ and $\widetilde{C}$ depend on $\mu_0$.
\end{lem}

\begin{proof}
 On the one hand, since $\|u_k\|_{L^p(\Omega)}=1$ and since the sequence $\{V_k\}_{k\geq0}$ is nested, then the discrete eigenvalues $\{\mu_k\}_{k\geq0}$ form a monotonously decreasing positive sequence satisfying
\begin{equation}\label{stab_eig}
    |u_k|^p_{W^{1,p}(\Omega)}=\mu_k\leq\mu_0\quad\forall~k\geq0.
\end{equation}
Note that $q=p/(p-1)>1$, then a straightforward calculation leads to
\begin{align}\label{eq:xxx}
h_T^{q}\mu_k^{q}\int_T|u_k|^p\dx\leq h_T^{q}\mu_0^{q}\|u_k\|_{L^p(T)}^p.
\end{align}
On the other hand, on each $F\in \mathcal{F}_k(\Omega)$ shared by two adjacent elements $T$, $T'\in\mathcal{T}_k$, an application of the scaled trace theorem \cite{ver} and the inverse estimate \cite{ciarlet} reveals
\[
    \begin{aligned}
        h_F\| [|\boldsymbol{\nabla}u_k|^{p-2}\boldsymbol{\nabla}u_k]\cdot\boldsymbol{n}_F \|^q_{L^q(F)} &\leq 2^{q-1} h_F \left( \int_F|(|\boldsymbol{\nabla}u_k|^{p-2}\boldsymbol{\nabla}u_k)_T|^{\frac{p}{p-1}}\ds + \int_{F}|(|\boldsymbol{\nabla}u_k|^{p-2}\boldsymbol{\nabla}u_k)_{T'}|^{\frac{p}{p-1}}\ds \right)\\
        &= 2^{q-1}  h_F \left(\|(\boldsymbol{\nabla}u_k)_T\|_{L^p(F)}^p + \|(\boldsymbol{\nabla}u_k)_{T'}\|_{L^p(F)}^p\right)\\
        &\leq C \|\boldsymbol{\nabla}u_k\|^p_{L^p(T\cup T')}.
    \end{aligned}
 \]
Then we get
 \[
    \sum_{F\subset\partial T\cap\Om}h_F\int_F|[|\boldsymbol{\nabla}u_k|^{p-2}\boldsymbol{\nabla}u_k]\cdot\boldsymbol{n}_F|^{\frac{p}{p-1}}\ds
    \leq C\int_{D_T}|\boldsymbol{\nabla}u_k|^{p}\dx.
\]
This, together with \eqref{eq:xxx}, leads to the first assertion.

Summing up the first inequality in \eqref{stab_err} over $T\in \mathcal{T}_k$, in combination with $\displaystyle h_T \leq \max_{T\in\mathcal{T}_0}{h}_T$, \eqref{stab_eig} and $\|u_k\|_{L^p(\Omega)}=1$, yields the second assertion. This completes the proof.
\end{proof}

\section{Convergence}\label{sec:conv}

In this section, we are concerned with the convergence of Algorithm \ref{afem} in the sense that the
sequence of discrete eigenvalues $\{\mu_k\}_{k\geq0}$ converges to the first eigenvalue $\lambda_1$
of Problem \eqref{vp_eigen} and the distance in $W^{1,p}(\Om)$-norm between $E_{\lambda_1}$
and the sequence of discrete eigenfunctions $\{u_k\}_{k\geq0}$ tends to zero. As in \cite{gmz1, gm} for eigenvalue problems associated with the linear diffusion operator, our analysis starts with an artificial minimization problem in Section \ref{subsect:min_lim}, with its solution being proved to
be the limit of the sequence of discrete eigenfunctions generated by Algorithm \ref{afem}. Then we invoke some auxiliary results on the error estimator $\eta_k(\mu_k,u_k)$ in Section \ref{subsect:aux}, and finally prove the desired convergence in Section \ref{subsect:conv}. It is interesting to note that $\eta_k(u_k,\mu_k)$ given by \eqref{eq:err-ind} is not a reliable estimator and the marking assumption \eqref{marking_max} will play an important role in the subsequent analysis.

\subsection{Limiting Behaviour}\label{subsect:min_lim}
With the sequence of $\{V_k\}_{k\geq0}$ generated by Algorithm \ref{afem}, we define a limiting space $V_\infty:=\overline{\bigcup_{k\geq0}V_k}\quad\mbox{in}~W^{1,p}(\Om)\mbox{-norm}$. It is not
difficult to know that $V_\infty$ is a closed subspace of $V$. We  now consider a limiting minimization problem: find $u_\infty\in V_\infty$ such that
\begin{equation}\label{min_lim}
    \mathcal{J}(u_\infty)=\inf_{v\in V_\infty\setminus\{0\}}\mathcal{J}(v).
\end{equation}

\begin{thm}\label{thm_min_lim}
 Problem \eqref{min_lim} has a nonnegative solution $u_\infty$ in
    $V_\infty\setminus\{0\}$.
\end{thm}

\begin{proof}
    Let $\{(\mu_k,u_k)\}_{k\geq0}$ be the sequence of discrete solutions given by Algorithm \ref{afem}. Since $\{V_k\}_{k\geq0}$ is nested, then by \eqref{eigva_lb} $\{\mu_k\}_{k\geq0}$ is a decreasing sequence bounded from below by $\lambda_1$. Consequently, there is $\mu_\infty>0$ such that $\mu_k\to\mu_\infty$. The identity $\|u_k\|_{L^{p}(\Omega)}=1$, together with $|u_k|_{W^{1,p}(\Omega)}^p=\mu_k$ in \eqref{stab_eig}, leads to the assertion that $\{\|u_k\|_{W^{1,p}(\Omega)}\}_{k\geq0}$ is bounded. 

On one hand, note that $V_\infty$ is a closed subspace of $W_0^{1,p}(\Omega)$, then the reflexivity and Sobolev compact embedding theorem \cite{adams} imply the existence of a subsequence $\{u_{k_j}\}_{j\geq 0}$ and some $u_\infty\in V_\infty$, satisfying
    \begin{equation}\label{pf1_thm_min_lim}
\begin{aligned}
        &u_{k_j}\to u_\infty~\m{weakly in}~W_0^{1,p}(\Om),\\
        & u_{k_j}\to u_\infty~\m{strongly in}~L^{p}(\Om),\\
        & u_{k_j}\to u_\infty~\m{a.e. in}~\Om.
\end{aligned}
    \end{equation}
    By $\|u_{k_j}\|_{L^{p}(\Omega)}=1$ and the second (i.e. the strong) convergence in \eqref{pf1_thm_min_lim}, we have
    \begin{equation}\label{pf2_thm_min_lim}
    \|u_\infty\|_{L^p(\Om)}=1.
    \end{equation}
    On the other hand, thanks to the definition of $V_\infty$, any $v\in V_\infty\setminus\{0\}$ admits a
    sequence $\{v_k\}_{k\geq0}\subset\bigcup_{k\geq 0}V_k$ with each $v_k\in V_k$ such that
    \begin{equation}\label{pf3_thm_min_lim}
        v_k \to v\quad\m{strongly in}~W^{1,p}_0(\Om).
    \end{equation}
    As each $u_k$ is a minimizer of $\mathcal{J}$ over $V_k\setminus\{0\}$, then 
    \begin{equation}\label{pf4_thm_min_lim}
        \mathcal{J}(u_k)\leq\mathcal{J}(v_k).
    \end{equation}
    Now using the first (i.e. the weak) convergence in \eqref{pf1_thm_min_lim} and collecting \eqref{pf2_thm_min_lim}-\eqref{pf4_thm_min_lim}, we arrive at
\begin{equation}\label{pf5_thm_min_lim}
    \begin{aligned}
        \mathcal{J}(u_\infty)=|u_\infty|^p_{W^{1,p}(\Omega)}&\leq\liminf_{j\to\infty}|u_{k_j}|^p_{W^{1,p}(\Omega)}\\
        &=\liminf_{j\to\infty}\mathcal{J}(u_{k_j})\\
        &\leq\limsup_{j\to\infty}\mathcal{J}(u_{k_j})\\
        &\leq\limsup_{k\to\infty}\mathcal{J}(u_k)\leq\limsup_{k\to\infty}\mathcal{J}(v_k)=\mathcal{J}(v)\quad\forall~v\in V_\infty\setminus\{0\}.
    \end{aligned}
\end{equation}
    This implies that $u_\infty$ is a minimizer over $V_\infty\setminus\{0\}$. As each $u_k$ is nonnegative, the third pointwise convergence in \eqref{pf1_thm_min_lim} implies that $u_\infty\geq 0$ a.e. in $\Om$. This completes the proof.
\end{proof}
\begin{thm}\label{thm_min_conv}
Let  $\{(\mu_k,u_k)\}_{k\geq0}$ be the sequence of discrete solutions generated by Algorithm \ref{afem}, then there holds
     \begin{equation}\label{min_conv1}
        \mu_k\to\mu_\infty=\mathcal{J}(u_\infty).
     \end{equation}
    Moreover, there exists a subsequence $\{u_{k_j}\}_{j\geq0}$ such that
     \begin{equation}\label{min_conv2}
        \|u_{k_j}-u_\infty\|_{W^{1,p}(\Omega)}\to 0.
     \end{equation}
\end{thm}

\begin{proof}
    At the beginning of the proof of Theorem \ref{thm_min_lim}, we showed that there is $\mu_\infty > 0$ such that $\mu_k\to\mu_\infty$. Taking $v=u_\infty$ in \eqref{pf5_thm_min_lim} implies
    \begin{equation}\label{pf1_thm_lim_conv}
    \mu_{k_j} = \mathcal{J}(u_{k_j}) \to \mathcal{J}(u_\infty). 
    \end{equation}
    By the uniqueness of the limit, we conclude that $\mu_\infty=\mathcal{J}(u_\infty)$. For the second assertion, thanks to the uniform convexity of $W_0^{1,p}(\Om)$ \cite{adams} and the weak convergence in \eqref{pf1_thm_min_lim}, it suffices to prove the norm convergence $\|u_{k_j}\|_{W^{1,p}(\Om)} \to \|u_{\infty}\|_{W^{1,p}(\Om)}$, which is an immediate consequence of the $L^p(\Om)$ strong convergence in \eqref{pf1_thm_min_lim}, easy facts that $\mu_{k_j}=|u_{k_j}|^p_{W^{1,p}(\Omega)}$ and $\mathcal{J}(u_\infty) = |u_{\infty}|^p_{W^{1,p}(\Omega)}$ as well as \eqref{pf1_thm_lim_conv}.
\end{proof}

\begin{rem}\label{rem_min_conv}
    In fact, using the arguments in the proof of Theorems \ref{thm_min_lim} and \ref{thm_min_conv}, from any subsequence $\{u_{k_j}\}_{j\geq0}$ of $\{u_k\}_{k\geq0}$ we can extract another subsequence $\{u_{k_{j_m}}\}_{m\geq0}$ converging in $W_0^{1,p}(\Om)$ to some $\widetilde{u}_\infty\in V_\infty$ satisfying $\mu_\infty=\mathcal{J}(\widetilde{u}_\infty)$.
\end{rem}

\subsection{Auxiliary Result}\label{subsect:aux}
Let $\{\mathcal{T}_{k}\}_{k\geq0}$ be the triangulation sequence generated by Algorithm \ref{afem}. First, we introduce the following notation,
\[
    \mathcal{T}_{k}^{+}:=\bigcap_{\ell\geq k}\mathcal{T}_{\ell},\quad
    \mathcal{T}_{k}^{0}:=\mathcal{T}_{k}\setminus\mathcal{T}_{k}^{+},\quad
    \Omega_{k}^{+}:=\bigcup_{T\in\mathcal{T}^{+}_{k}}D_{T},\quad
    \Omega_{k}^{0}:=\bigcup_{T\in\mathcal{T}^{0}_{k}}D_{T}.
\]
By definition, $\mathcal{T}_{k}^{+}$ consists of all elements not refined after the $k$-th iteration and
any element in $\mathcal{T}_{k}^{0}$ are refined at least once after the $k$-th iteration.
Note that $\mathcal{T}_{\ell}^{+}\subset\mathcal{T}_{k}^{+}\subset \mathcal{T}_{k}$ for $\ell<k$.

Next, we define a mesh-size function
$h_{k}:\Omega\rightarrow\mathbb{R}^{+}$ almost everywhere by $h_{k}(x)=h_{T}$ for $x$ in the interior of an element $T\in\mathcal{T}_{k}$ and
$h_{k}(x)=h_{F}$ for $x$ in the relative interior of face $F\in\mathcal{F}_{k}(\Om)$.
This mesh-size function has the following property \cite{nsv},
\begin{equation}\label{con_0_mz}
    \lim_{k\rightarrow\infty}\|h_{k}\chi^{0}_{k}\|_{L^\infty(\Omega)}=0,
\end{equation}
where $\chi^{0}_{k}$ is the characteristic function of $\Omega_{k}^{0}$. Note that the uniform refinement strategy corresponds to $\|h_k\|_{L^\infty(\Omega)}\to0$ as $k\to\infty$, and the resulting sequences of nested spaces satisfies $W_0^{1,p}(\Omega)=\overline{\bigcup_{k\geq0}V_k}$.

\begin{lem}\label{lem_maxerr}
 Let $\{(\mu_{k_j},u_{k_j})\}_{j\geq0}$ be the convergent subsequence defined in Theorem \ref{thm_min_conv} and let $\{\mathcal{M}_{k_j}\}_{j\geq0}$ be the associate sequence of marked element patches, then there holds
    \begin{equation}\label{con_0_me}
        \lim_{j\rightarrow\infty}\max_{T\in\mathcal{M}_{k_j}}\eta_{k_j}(\mu_{k_j},u_{k_j};T)=0.
    \end{equation}
\end{lem}

\begin{proof}
  Let $T_j\in\mathcal{M}_{k_j}$ be the element with the largest error indicator over $\mathcal{M}_{k_j}$ for each $k_j$. As
    $\mathcal{M}_{k_j}\subset\mathcal{T}_{k_j}^0$, it is not difficult to know from
    \eqref{con_0_mz} that
    \begin{equation}\label{pf1_lem_maxerr}
\begin{aligned}
    h_{T_j}&\leq\|h_{k_j}\chi^{0}_{k_j}\|_{L^\infty(\Omega)}\to 0 &&\quad\m{as}~j\to\infty,\\
 |T_j|\leq |D_{T_j}|&\leq C\|h_{k_j}\chi^{0}_{k_j}\|_{L^\infty(\Omega)}^{d}\rightarrow 0&&\quad\m{as}~j\to\infty.
\end{aligned}
    \end{equation}
    By virtue of \eqref{stab_err} in Lemma \ref{lem_stab_err}, we have
    \begin{align*}
    \eta^q_{k_j}(\mu_{k_j},u_{k_j};T_j)&\leq C(h^q_{T_j}\|u_{k_j}\|_{L^p(T_j)}^{p}+\|\boldsymbol{\nabla}u_{k_j}\|_{L^p(D_{T_j})}^{p})\\
    &\leq C\Big(\|u_{k_j}-u_\infty\|_{L^{p}(\Omega)}^{p}+\|u_\infty\|_{L^p(T_j)}^{p}+\|\boldsymbol{\nabla}(u_{k_j}-u_\infty)\|_{L^p(\Omega)}^{p}+\|\boldsymbol{\nabla}u_{\infty}\|_{L^p(D_{T_j})}^{p}\Big).
    \end{align*}
    Therefore, the desired vanishing limit comes from \eqref{pf1_lem_maxerr}, \eqref{min_conv2} in Theorem \ref{thm_min_conv} and the absolute continuity of $\|\cdot\|_{L^p(\Omega)}$ with respect to the Lebesgue measure.
\end{proof}

Next, we introduce the residual with respect to the eigenpair $(\mu_k,u_k)$,
\[
    \langle\mathcal{R}(\mu_k,u_k),v\rangle:=\int_\Om|\bold{\nabla}u_k|^{p-2}\bold{\nabla}u_k\cdot\bold{\nabla}v\dx
    -\mu_k\int_\Om|u_k|^{p-2}u_kv\dx\quad\forall~v\in W_0^{1,p}(\Om).
\]
To establish the convergence of this residual, we need to first recall the nodal interpolation operator $I_j: W^{2,\infty}(\Omega)\cap W_0^{1,p}(\Omega)\to V_j$ and the Scott-Zhang quasi-interpolation $\widetilde{I}_j: W^{1,p}_0(\Omega)\to V_j$, which have the following approximation properties \cite{ciarlet,sz},
\begin{align}
    \|\nabla(v-I_j v)\|_{L^p(T)}&\leq C h_T^{1+d/p}
|v|_{W^{2,\infty}(T)}&& \forall~T\in \mathcal{T}_j,\label{eq:nodal-interp}\\
\|v-\widetilde{I}_j v\|_{L^p(T)} + \sum_{F\textcolor{red}{\subset}\partial T\cap \Omega}h_F^{1/p} \|v-\widetilde{I}_j v\|_{L^p(F)} & \leq C h_T\|\nabla v\|_{L^p(D_T)} && \forall~T\in \mathcal{T}_j.\label{eq:sz-quasi-interp}
\end{align}

\begin{lem}\label{lem_res_conv}
Let $\{(\mu_{k_j},u_{k_j})\}_{j\geq0}$ be the convergent subsequence defined in Theorem \ref{thm_min_conv}, there holds
    \begin{equation}\label{res_conv}
        \lim_{j\to\infty}\langle\mathcal{R}(\mu_{k_j},u_{k_j}),v\rangle=0\quad\forall~v\in C_0^{\infty}(\Om).
    \end{equation}
\end{lem}

\begin{proof}
 For the sake of brevity, $k_j$ is abbreviated to $j$. For any $v\in C_0^\infty(\Omega)$, invoking the nodal interpolation operator $I_j$ and the Scott-Zhang quasi-interpolation $\widetilde{I}_j$ associated with $V_j$, and noting the eigenpair $(\mu_{j},u_{j})$ satisfies \eqref{vp_disc} over $\mathcal{T}_j$, we derive
    \begin{align*}
        \left|\langle\mathcal{R}(\mu_{j},u_{j}),v\rangle\right|
        &=\left|\langle\mathcal{R}(\mu_{j},u_{j}),v-I_jv\rangle\right|
        =\left|\langle\mathcal{R}(\mu_{j},u_{j}),w-\widetilde{I}_jw\rangle\right|.
    \end{align*}
Here, we denote $w:=v-I_jv$. Combined with the elementwise integration by parts and H\"{o}lder inequality, this leads to
\begin{align*}
        \left|\langle\mathcal{R}(\mu_{j},u_{j}),v\rangle\right|
&=\left|-\sum_{T\in\mathcal{T}_j}\int_T R_{T}(\mu_j,u_j)(w-\widetilde{I}_jw)\dx-\sum_{F\in\mathcal{F}_j(\Om)}\int_F J_F(u_j) (w-\widetilde{I}_jw)\ds\right|\\
        &\leq \sum_{T\in\mathcal{T}_j}\|R_T(\mu_j,u_j)\|_{L^q(T)}\|w-\widetilde{I}_jw\|_{L^p(T)}+\sum_{F\in\mathcal{F}_j(\Om)}\|J_F(u_j)\|_{L^q(F)}\|w-\widetilde{I}_jw\|_{L^{p}(F)}.
    \end{align*}
We further proceed by the error estimate \eqref{eq:sz-quasi-interp}, $\frac{1}{p}+\frac{1}{q}=1$ and a split of $\mathcal{T}_j$ into $\mathcal{T}_{\ell}^+$ and $\mathcal{T}_j\setminus\mathcal{T}_{\ell}^+$ for some $\ell<j$ to find
    \begin{align*}
        \left|\langle\mathcal{R}(\mu_{j},u_{j}),v\rangle\right|&\leq
        C\sum_{T\in\mathcal{T}_j}\left(h_T\|R_T(\mu_j,u_j)\|_{L^q(T)}\|\bold{\nabla}w\|_{L^p(D_T)}+\sum_{F\subset\partial T\cap\Om}h_F^{1/q}\|J_F(u_j)\|_{L^q(F)}\|\boldsymbol{\nabla}w\|_{L^p(D_{T})}\right)\\
        &\leq C\sum_{T\in\mathcal{T}_j}\eta_{j}(\mu_j,u_j;T)\|\boldsymbol{\nabla}(v-I_jv)\|_{L^p(D_T)}\\
        &\leq C\left(\eta_j(\mu_j,u_j;\mathcal{T}_j\setminus\mathcal{T}_{\ell}^+)\|\nabla(v-I_jv)\|_{L^{p}(\Omega_{\ell}^0)}+\eta_j(\mu_j,u_j;\mathcal{T}_{\ell}^+)\|\nabla(v-I_jv)\|_{L^p(\Omega_{\ell}^+)}\right).
    \end{align*}
   With $\T_j(\Omega_\ell^0)$ denoting the restriction of $\T_j$ over $\Omega_{\ell}^0$, the error estimate for $I_j$ \eqref{eq:nodal-interp} and the fact $|T|=h_T^d$ for any $T\in \T_j$ imply that for any $v \in C_0^\infty(\Omega)$,
   \begin{align*}
        \| \boldsymbol{\nabla} ( v - I_j v ) \|_{L^p(\Omega_\ell^0)}^p &\leq C \sum_{T\in\T_j(\Omega_\ell^0)} h_T^{p+d} |v|_{W^{2,\infty}(T)}^p
\leq C |\Omega| \|h_j\|^p_{L^\infty(\Om_{\ell}^0)} \|v\|^p_{W^{2,\infty}(\Omega)},\\
    \| \boldsymbol{\nabla} ( v - I_j v ) \|_{L^p(\Omega_\ell^+)}^p
&\leq C \sum_{T\in\T_j} h_T^{p+d} |v|_{W^{2,\infty}(T)}^p
\leq  C \|v\|^p_{W^{2,\infty}(\Omega)}.
   \end{align*}
   Therefore, by the stability estimate \eqref{stab_err} we arrive at
    \begin{equation}\label{pf1_lem_res_conv}
        \left|\langle\mathcal{R}(\mu_{j},u_{j}),v\rangle\right|\leq
        C_1\|h_j\|_{L^\infty(\Om_{\ell}^0)}\|v\|_{W^{2,\infty}(\Omega)}+C_2\eta_{j}(\mu_j,u_j;\mathcal{T}_{\ell}^+)\|v\|_{W^{2,\infty}(\Omega)}\quad\forall~v\in C_0^{\infty}(\Om).
    \end{equation}
 Let $\ell\to\infty$, then the first term on the right hand side of \eqref{pf1_lem_res_conv} goes to zero due to the monotonicity $h_j\leq h_{\ell}$ and \eqref{con_0_mz}.

To estimate the second term, the marking assumption \eqref{marking_max} in Algorithm \ref{afem}, combined with the  fact $\mathcal{T}_{\ell}^+\subset\mathcal{T}_j^+\subset\mathcal{T}_j$ for $j>\ell$, leads to
    \begin{equation*}
\begin{aligned}
        \eta_{j}(\mu_j,u_j;\mathcal{T}_{\ell}^+)&\leq(\#\mathcal{T}_{\ell}^+)^{1/q}\max_{T\in\mathcal{T}_{\ell}^+}\eta_j(\mu_j,u_j;T)\\
 &\leq (\#\mathcal{T}_{\ell}^+)^{1/q}\max_{T\in\mathcal{T}_j^+}\eta_j(\mu_j,u_j;T)\\
        &\leq (\#\mathcal{T}_{\ell}^+)^{1/q}\max_{T\in\mathcal{M}_j}\eta_j(\mu_j,u_j;T),
    \end{aligned}
\end{equation*}
where $\# \mathcal{T}_l^+ $ denotes the number of all elements in $\mathcal{T}_l^+$. Then Lemma \ref{lem_maxerr} implies the second term tends to zero when $j\to\infty$. This proves the assertion.
\end{proof}

\begin{lem}\label{lem_lim-con}
 Let $\{(\mu_{k_j},u_{k_j})\}_{j\geq0}$ be the convergent subsequence given by Theorem \ref{thm_min_conv}, then the following identity holds
    \begin{equation}\label{lim-con}
        \lim_{j\to\infty}\langle\mathcal{R}(\mu_{k_j},u_{k_j}),v\rangle=
        \int_\Om|\bold{\nabla}u_\infty|^{p-2}\bold{\nabla}u_\infty\cdot\bold{\nabla}v\dx
        -\mu_\infty\int_\Om|u_\infty|^{p-2}u_\infty v\dx\quad\forall~v\in W_0^{1,p}(\Om).
    \end{equation}
\end{lem}

\begin{proof}
    We introduce two functionals
    $\mathcal{F}(v):=\int_{\Om}|\bold{\nabla}v|^p\dx/p$ and       $\mathcal{G}(v):=\int_\Om|v|^p\dx/p$ on $W_0^{1,p}(\Omega)$.
   As $\mathcal{F}$ and $\mathcal{G}$ are both $C^1$ functionals (cf. e.g. \cite{BS}) with
\begin{equation}
    \begin{aligned}
        &\mathcal{F}'(w)v=\int_\Om|\bold{\nabla}w|^{p-2}\bold{\nabla}w\cdot\bold{\nabla}v\dx,\\
&\mathcal{G}'(w)v=\int_\Om|w|^{p-2}wv\dx\quad\forall ~w, v\in W_0^{1,p}(\Om).
    \end{aligned}
\end{equation}
     Theorem \ref{thm_min_conv} leads to the desired assertion.
\end{proof}

\subsection{Main Result}\label{subsect:conv}

Now we are in a position to state the main result of this paper.

\begin{thm}\label{thm_afem_conv}
Assume that the initial mesh $\mathcal{T}_0$ is sufficiently fine, i.e. $\|h_0\|_{L^\infty(\Omega)}\ll 1$. Let $\{(\mu_k,u_k)\}_{k\geq0}$ be the sequence of discrete eigenpairs produced by Algorithm \ref{afem}, then the following identities hold
    \begin{equation}\label{afem_conv}
\begin{aligned}
        &\lim_{k\to\infty}\mu_k=\lambda_1,\\
&\lim_{k\to\infty}\inf_{v\in E_{\lambda_1}}\|u_k-v\|_{W^{1,p}(\Omega)}=0.
\end{aligned}
    \end{equation}
\end{thm}

\begin{proof}
The proof is divided into three steps.

    \noindent\textbf{Step 1.} By Theorem \ref{thm_min_conv}, $\mu_k\to\mu_\infty$ and there exists a subsequence
    $\{u_{k_j}\}$ such that $\|u_{k_j}-u_\infty\|_{W^{1,p}(\Omega)}\to0$. First we prove that
    $(\mu_\infty,u_\infty)$ is an eigenpair of Problem \eqref{vp_eigen}.
    By the H\"{o}lder inequality and the stability estimate \eqref{stab_eig},
    \begin{align}
        \left|\langle\mathcal{R}(\mu_k,u_k),v\rangle\right|&=\left|\int_\Om|\bold{\nabla}u_k|^{p-2}\bold{\nabla}u_k\cdot\bold{\nabla}v\dx
        -\mu_k\int_\Om|u_k|^{p-2}u_kv\dx\right|\nonumber\\
        &\leq|u_k|^{p-1}_{W^{1,p}(\Omega)}|v|_{W^{1,p}(\Omega)}+\mu_k\|u_k\|^{p-1}_{L^{p}(\Omega)}\|v\|_{L^{p}(\Omega)}\\
        &\leq C(\mu_0)\|v\|_{W^{1,p}(\Omega)}\quad\forall~v\in W_0^{1,p}(\Om).\nonumber
    \end{align}
 Hence, the sequence $\left\{\|\mathcal{R}(\mu_k,u_k)\|_{W^{-1,q}(\Om)}\right\}_{k\geq0}$ is uniformly bounded. This, together with Lemma \ref{lem_res_conv} and the density of $C_0^{\infty}(\Om)$ in $W_0^{1,p}(\Omega)$, implies that the convergent subsequence $\{(u_{k_j},\mu_{k_j})\}_{j\geq0}$ satisfies
    \[
        \lim_{j\to\infty}\langle\mathcal{R}(\mu_{k_j},u_{k_j}),v\rangle=0\quad\forall~v\in W_0^{1,p}(\Om).
    \]
In combination with Lemma \ref{lem_lim-con}, we derive
    \begin{equation}\label{pf3_thm_afem_conv}
        \int_\Om|\bold{\nabla}u_\infty|^{p-2}\bold{\nabla}u_\infty\cdot\bold{\nabla}v\dx
        =\mu_\infty\int_\Om|u_\infty|^{p-2}u_\infty v\dx\quad\forall~v\in W_0^{1,p}(\Om).
    \end{equation}
This means that $(\mu_\infty,u_\infty)$ is an eigenpair of Problem \eqref{vp_eigen}.

\noindent\textbf{Step 2.} In view of \eqref{pf3_thm_afem_conv} and Theorem \ref{thm_min_conv}, the first result in \eqref{afem_conv} is true once $\mu_\infty=\lambda_1$ is proved. To this end we define $E:=\{u\in W_0^{1,p}(\Om)|~u~\m{satisfies \eqref{vp_eigen} with}~\|u\|_{L^p(\Omega)}=1~\m{for some }\lambda\in\mathbb{R}\}$, i.e. $E$ consists of all $L^p(\Omega)$-normalized eigenfunctions of Problem \eqref{vp_eigen}. Problem \eqref{vp_eigen} has a nondecreasing sequence of positive eigenvalues tending to $+\infty$ \cite{gp, le}, so there holds $E_{\lambda_1}\subsetneq E$. Since any $u\in E_{\lambda_1}$ is a minimizer of $\mathcal{J}$ over $W_0^{1,p}(\Om)$, then we derive $\lambda_1=\mathcal{J}(u)\leq\inf_{v\in E\setminus
    E_{\lambda_1}}\mathcal{J}(v)$. We prove by contradiction that the equality cannot happen. Let
    $\{w_n\}_{n\geq0}\subset E\setminus E_{\lambda_1}$ be a minimizing sequence such that
    \[
        \mathcal{J}(w_n)\to\inf_{v\in E\setminus
        E_{\lambda_1}}\mathcal{J}(v) = \lambda_1.
    \]  
    Then $\mathcal{J}(w_n)$ is a sequence of
    eigenvalues with $\lambda_1$ as its limit, contradicting the fact that $\lambda_1$ is isolated \cite{lind}. Thus, $\lambda_1 = \mathcal{J}(u) < \inf_{v\in E\setminus E_{\lambda_1}}\mathcal{J}(v)$ for any $u\in E_{\lambda_1}$.

Next, we justify the fineness condition on the initial mesh $\mathcal{T}_0$ in Algorithm \ref{afem}. Assuming $\{\mathcal{T}_k\}_{k\geq0}$ is a sequence of uniformly refined meshes, at this point $\|h_k\|_{L^\infty(\Omega)}\to 0$ and $W_0^{1,p}(\Omega)=\overline{\bigcup_{k\geq0}V_k}$ as mentioned in Section \ref{subsect:aux}. Therefore, for any $u\in E_{\lambda_1}$ there exists a sequence $\{v_{\ell}\}_{\ell\geq0}$ with each $v_\ell\in V_{\ell}$ such that $v_{\ell}\to u$ strongly in $W_0^{1,p}(\Om)$. Noting $\mathcal{J}$ is continuous over $W_0^{1,p}(\Om)$ and $\mathcal{J}(u) = \lambda_1 <\inf_{v\in E\setminus
    E_{\lambda_1}}\mathcal{J}(v)$, we have $\mathcal{J}(v_{\ell})<\inf_{v\in E\setminus
    E_{\lambda_1}}\mathcal{J}(v)$ for sufficiently large $\ell$ or sufficiently small mesh-size $\|h_{\ell}\|_{L^\infty(\Omega)}$. This observation and Theorem \ref{thm_min_conv} imply that for the sequence of adaptively generated meshes $\{\mathcal{T}_k\}_{k\geq0}$ by Algorithm \ref{afem}, we may choose a fine enough initial mesh $\mathcal{T}_0$, over which there holds
    \begin{equation}\label{pf4_thm_afem_conv}
        \mathcal{J}(u_\infty)\leq\mathcal{J}(u_0)<\inf_{v\in E\setminus
        E_{\lambda_1}}\mathcal{J}(v).
   \end{equation}
   On the other hand, by \eqref{pf3_thm_afem_conv} in Step 1, $u_\infty\in E$. If $u_\infty \in\hspace{-9pt}\slash \hspace{3pt} E_{\lambda_1}$, then $\mathcal{J}(u_\infty)\geq \inf_{v\in E\setminus E_{\lambda_1}}\mathcal{J}(v)$ and \eqref{pf4_thm_afem_conv} yields an obvious contradiction. Therefore, $u_\infty \in E_{\lambda_1}$ and $\mu_\infty= \mathcal{J}(u_\infty) =\lambda_1$. 

    \noindent\textbf{Step 3.} To prove the second result in \eqref{afem_conv}, we also proceed with mathematical contradiction. If the
    result is false, there exist a number $\varepsilon>0$ and a subsequence $\{u_{k_j}\}_{j\geq0}$ of \textcolor{red}{$\{u_k\}_{k\geq 0}$} such that 
    \[
        \inf_{v\in E_{\lambda_1}}\|u_{k_j}-v\|_{W^{1,p}(\Omega)}\geq\varepsilon\quad \text{for all}~k_j.
    \]
    As discussed in Remark \ref{rem_min_conv}, we may extract another subsequence of $\{u_{k_{j_m}}\}_{m\geq0}$ converging to some $\widetilde{u}\in V_\infty$. Using the argument in Step 1 and the
    first result in \eqref{afem_conv}, we further know \textcolor{red}{that} $\widetilde{u}$ satisfies
    \eqref{pf3_thm_afem_conv} with $\mu_\infty=\lambda_1$, i.e. $\widetilde{u}\in E_{\lambda_1}$. This is a contradiction.
\end{proof}

\section{Numerical Examples}
To demonstrate the performance of Algorithm \ref{afem}, we consider three 2-d numerical tests with a unit disk, a unit square and a L-shaped domain and three 3-d numerical examples with a unit cube, a 3-d L-shaped domain and a torus as the computational domain respectively in this section. All experiments are implemented in MATLAB R2020a on a personal computer with a 13th Gen Intel(R) Core(TM) i7-13700 (24 CPUs) and 32GB RAM, while the Partial Differential Equation Toolbox and the MATLAB package $i$FEM \cite{Chen:2009} are used in numerical simulations of 2-d and 3-d examples, respectively. 

We utilize a normalized inverse iteration of sublinear supersolutions (IISS) \cite[Algorithm 2]{BBEM} to solve Problem \eqref{min_disc} over each mesh level, which is repeated in Algorithm \ref{alg:alg1} for the sake of completeness.
\begin{algorithm}
\begin{algorithmic}[1]
	\State Solve (\textit{torsion funcion})
\begin{equation*}
\left\{
\begin{aligned}
-\Delta_pu_0:=-\boldsymbol{\nabla}\cdot(|\boldsymbol{\nabla}u_0|^{p-2}\boldsymbol{\nabla}u_0)&=1&&\text{ in }\Omega\\
 u_0&=0&&\text{ on }\partial\Omega.
\end{aligned}
\right.
\end{equation*}
\State $m\leftarrow 0$.
\State {$\lambda_{m} = 1/\|u_{m}\|_{L^\infty(\Omega)}^{p-1}$}. 
\Do
    \State $m\leftarrow m+1$.
	\State Solve (\textit{inverse iteration})
    \begin{equation*}
        \left\{
            \begin{aligned}
                -\Delta_pu_{m}&=({u_{m-1}}/{\|u_{m-1}\|_{L^\infty(\Omega)}})^{p-1}\ &&{\rm in}\ \Omega,\\
                u_{m}&=0\ &&{\rm on} \ \partial\Omega.
            \end{aligned}
        \right.
    \end{equation*}
    \State  {$\lambda_{m} = 1/\|u_{m}\|_{L^\infty(\Omega)}^{p-1}$}. 
\doWhile{$|\lambda_{m}-\lambda_{m-1}|/|\lambda_{m-1}|\geqslant\epsilon_{M}$}
	\State Return {$\lambda_m$} and {$u_{m}/{\|u_{m}\|_{L^{\infty}(\Omega)}}$}.\quad(\textit{first eigenvalue and first eigenfunction})
\end{algorithmic}
\caption{Normalized IISS \cite[Algorithm 2]{BBEM}}
\label{alg:alg1}
\end{algorithm}
Note that Algorithm \ref{alg:alg1} involves solving a $p$-Laplacian problem for the torsion function in Step 1 and the inverse iteration sequence in Step 6. To do so, we call a decomposition coordination algorithm \cite{Aragon} as presented in Algorithm \ref{alg:alg2} with $f$ being the right hand side of the $p$-Laplacian problem in Steps 1 and 6 of Algorithm \ref{alg:alg1} and $g=0$ in the current situation.
\begin{algorithm}
    \begin{algorithmic}[1]
        \State Define two vector fields $\boldsymbol{\xi}_{1}$, $\boldsymbol{\nu}_{0}$: $\Omega \rightarrow \mathbb{R}^{d}$; $n\leftarrow 0$.

\Do
    \State $n\leftarrow n+1$.
 \State Compute  $u_{n}$  by solving a linear problem
\begin{equation*}
\left\{
\begin{aligned}
&-\Delta u_{n}=\nabla \cdot\left(\boldsymbol{\xi}_{n} - \boldsymbol{\nu}_{n-1}\right)+f ~&&\text{in}~\Omega,\\
&u_{n}=g~&&\text{on}~\partial \Omega.
\end{aligned}
\right.
\end{equation*}
 \State Compute $\boldsymbol{\nu}_{n}$  by solving the algebraic nonlinear equation $\left|\boldsymbol{\nu}_{n}\right|^{p-2} \boldsymbol{\nu}_{n}+\boldsymbol{\nu}_{n}=\boldsymbol{\xi}_{n}+\nabla u_{n}$.

  \State Compute $\boldsymbol{\xi}_{n+1}$ as $\boldsymbol{\xi}_{n+1} = \boldsymbol{\xi}_{n}+\nabla u_{n} - \boldsymbol{\nu}_{n}$.

     \doWhile{$n=1$ or  $\frac{\|u_{n}-u_{n-1}\|_{L^2(\Omega)}}{\| u_{n-1}\|_{L^2(\Omega)}}>\epsilon_{N}$}
 	
        \State Return ${u_{n}}$.

    \end{algorithmic}
\caption{Decomposition Coordination \cite{Aragon}}
\label{alg:alg2}
\end{algorithm}

\begin{figure}[htbp]
    \centering
        \begin{subfigure}{.3\textwidth}
            \centering
            \includegraphics[scale=0.3]{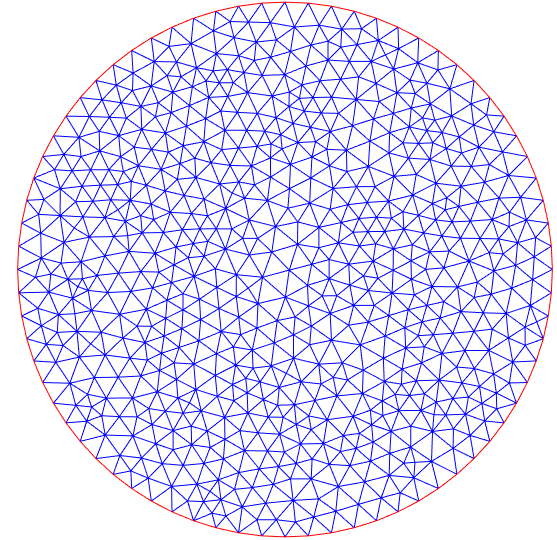}
            \label{fig:initial-mesh_example1}
        \end{subfigure}
        \begin{subfigure}{.3\textwidth}
            \centering
            \includegraphics[scale=0.3]{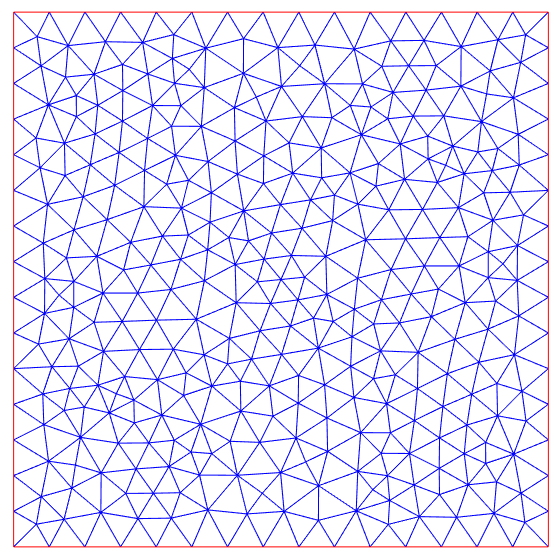}
            \label{fig:initial-mesh_example2}
        \end{subfigure}
        \begin{subfigure}{.3\textwidth}
            \centering
            \includegraphics[scale=0.3]{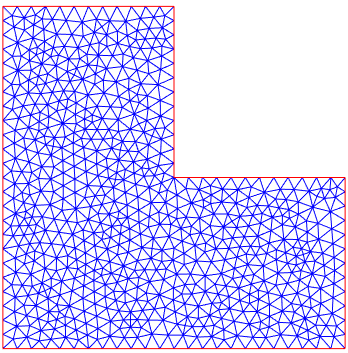}
            \label{fig:initial-mesh_example3}
        \end{subfigure}
    \caption{Initial meshes ($k=0$ for $\mathcal{T}_0$) with the number of vertices being 682, 365, and 741 for Examples \ref{example1}-\ref{example3} from left to right.}\label{fig:initial-mesh}
\end{figure}

\begin{figure}[htbp]
    \centering
        \begin{subfigure}{.04\textwidth}
            \centering
            \includegraphics[scale=0.2]{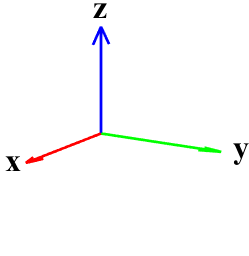}
        \end{subfigure}
        \begin{subfigure}{.27\textwidth}
            \centering
            \includegraphics[scale=0.3]{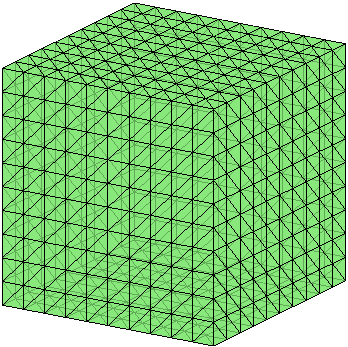}
            \label{fig:initial-mesh3d_example4}
        \end{subfigure}
        \hspace{0cm}
        \begin{subfigure}{.04\textwidth}
            \centering
            \includegraphics[scale=0.2]{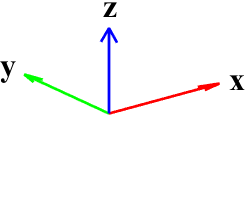}
        \end{subfigure}
        \begin{subfigure}{.26\textwidth}
            \centering
            \includegraphics[scale=0.35]{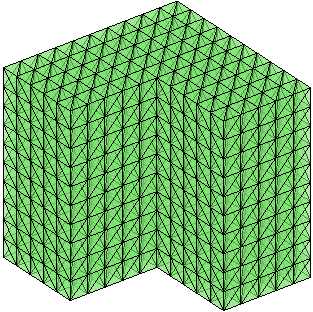}
            \label{fig:initial-mesh3d_example5}
        \end{subfigure}
        \begin{subfigure}{.04\textwidth}
            \centering
            \includegraphics[scale=0.2]{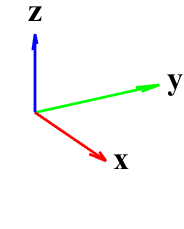}
        \end{subfigure}
        \begin{subfigure}{.26\textwidth}
            \centering
            \includegraphics[scale=0.35]{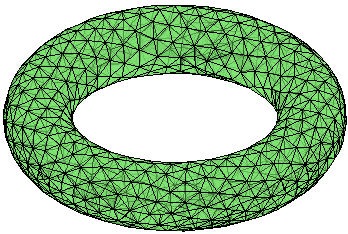}
            \label{fig:initial-mesh3d_example6}
        \end{subfigure}
    \caption{Initial meshes ($k=0$ for $\mathcal{T}_0$) with the number of vertices being 1331, 585, and 959 for Examples \ref{example4}-\ref{example6} from left to right.}\label{fig:initial-mesh3d}
\end{figure}

In all experiments, the module MARK of Algorithm \ref{afem}, utilizing D\"{o}rfler's strategy with $\theta=0.6$ for the 2-d examples and $\theta=0.8$ for the 3-d examples, yields a subset $\mathcal{M}_k$ such that $\eta_{k}(\mu_k,u_k,\mathcal{M}_k)\geq \theta \eta_{k}(\mu_k,u_k)$. Algorithm \ref{afem} proceeds until the relative error for two consecutive approximate eigenvalues is below a prescribed tolerance $\epsilon_{K}$, i.e. $|\mu_{k-1} - \mu_k|/ \mu_{k-1} <  \epsilon_{K}$. For large $p$ in the 2-d examples and all $p$ in the 3-d examples, an upper bound $K$ is specified for the counter $k$ of adaptive refinement steps. In Algorithm \ref{alg:alg2}, we fix tolerance $\epsilon_{N}=10^{-5}$, 
and each component of $\boldsymbol{\xi}_{1}$ and $\boldsymbol{\nu}_{0}$ is an independent sample following the uniform distribution $U(0,0.5)$. Figures \ref{fig:initial-mesh}-\ref{fig:initial-mesh3d} display all initial meshes used in Examples \ref{example1}-\ref{example6}. 

\begin{figure}[htb!]
    \centering
    \subcaptionbox{$\mathcal{T}_1$ (1050) \label{fig:UnitDiskP1_1Mesh2TH}}    {\includegraphics[scale=0.31]{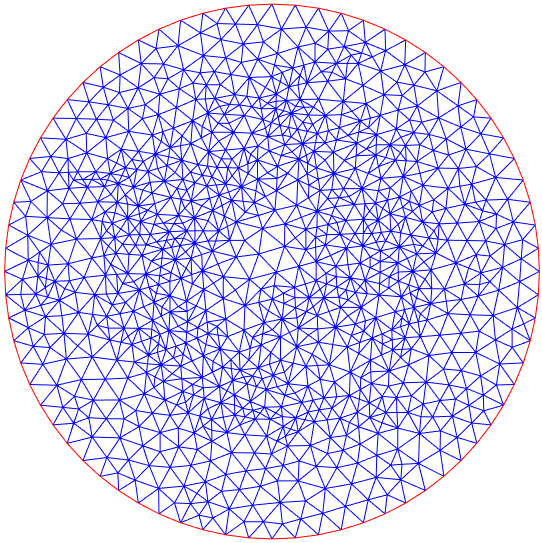}}\hspace{0.5cm}
    \subcaptionbox{$\mathcal{T}_3$ (2498) \label{fig:UnitDiskP1_1Mesh4TH}}
    {\includegraphics[scale=0.31]{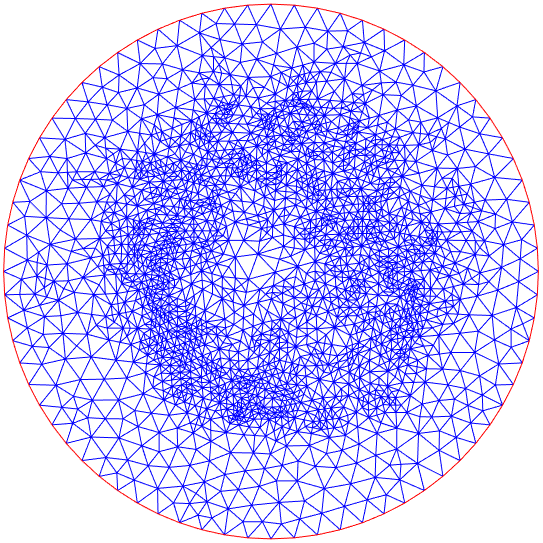}}\hspace{0.5cm}
    \subcaptionbox{$u_3$\label{fig:UnitDiskP1_1U4TH}}
    {\includegraphics[scale=0.31]{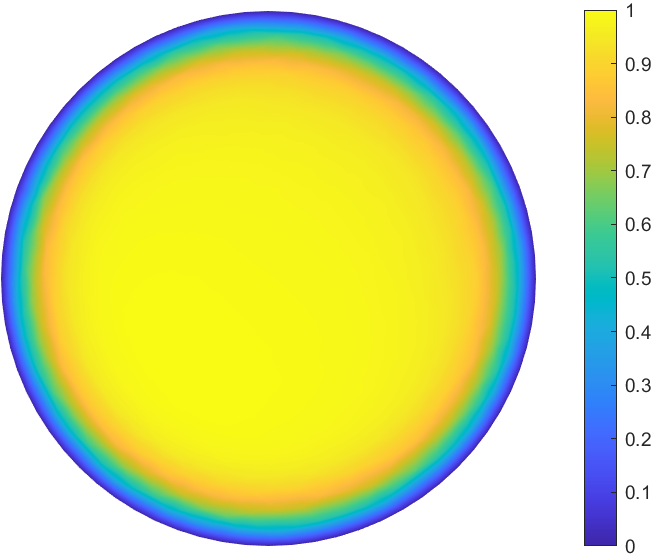}}\hspace{0.5cm}
	\subcaptionbox{reference\label{fig:UnitDiskP1_1U1THUniform}}
    {\includegraphics[scale=0.31]{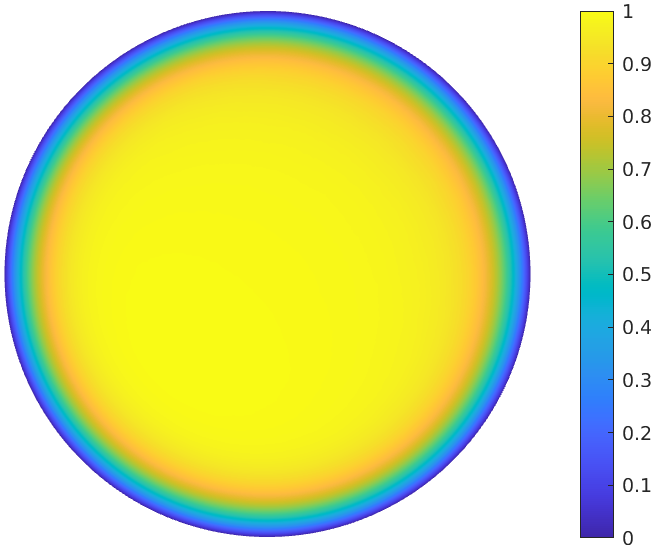}}
\\
    \subcaptionbox{$\mathcal{T}_2$ (2188) \label{fig:UnitDiskP1_5Mesh3TH}}
    {\includegraphics[scale=0.3]{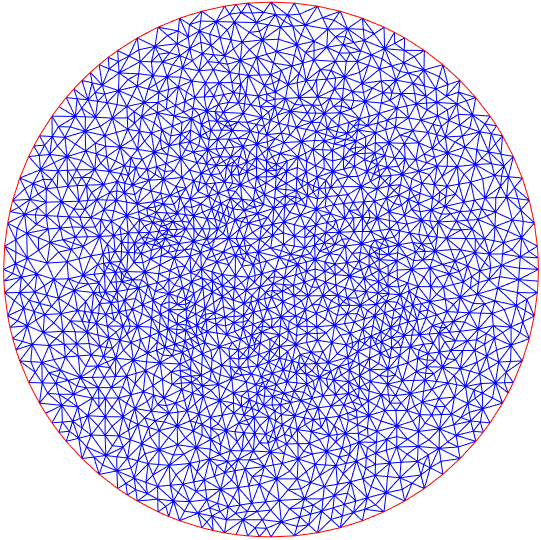}}\hspace{0.5cm}
    \subcaptionbox{$\mathcal{T}_5$ (12469) \label{fig:UnitDiskP1_5Mesh6TH}}
    {\includegraphics[scale=0.3]{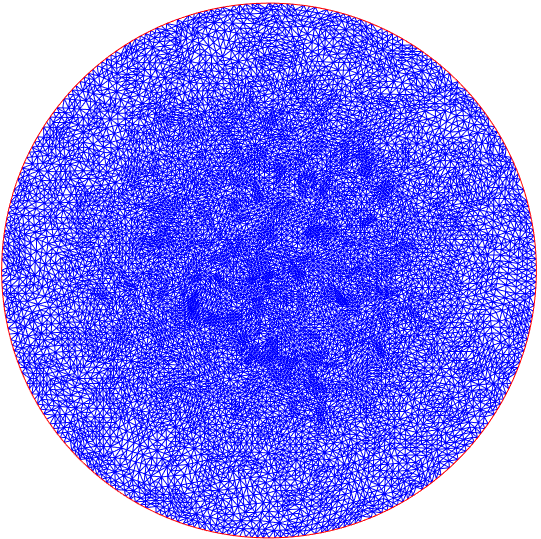}}\hspace{0.5cm}
    \subcaptionbox{$u_7$ \label{fig:UnitDiskP1_5U8TH}} {\includegraphics[scale=0.3]{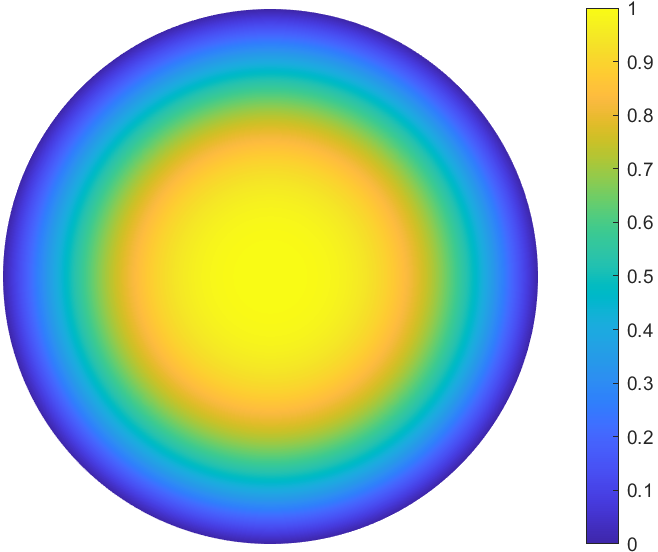}}\hspace{0.5cm}
    \subcaptionbox{reference\label{fig:UnitDiskP1_5U1THUniform}} {\includegraphics[scale=0.3]{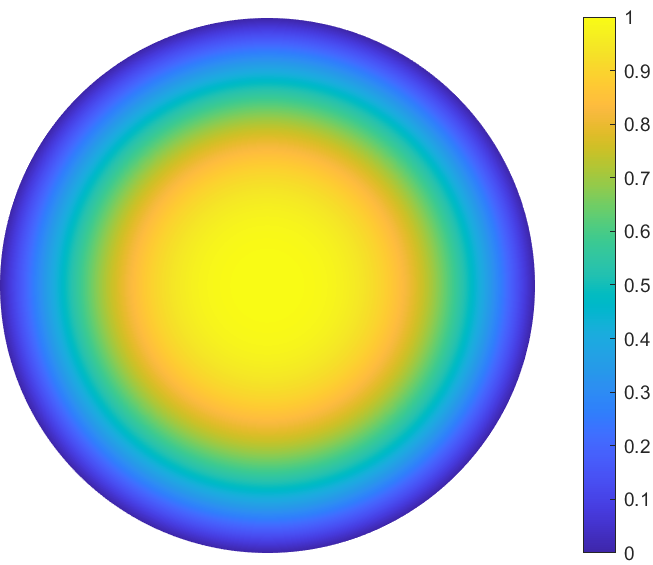}}
\\
    \subcaptionbox{$\mathcal{T}_1$ (1196) \label{fig:UnitDiskP2Mesh2TH}} {\includegraphics[scale=0.3]{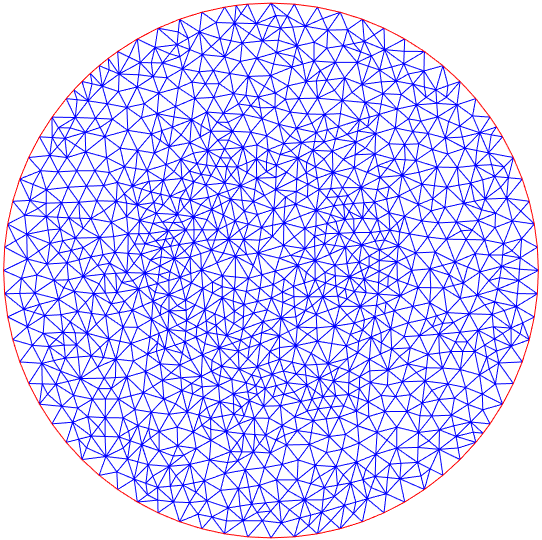}}\hspace{0.5cm}
    \subcaptionbox{$\mathcal{T}_5$ (11639) \label{fig:UnitDiskP2Mesh6TH}} {\includegraphics[scale=0.3]{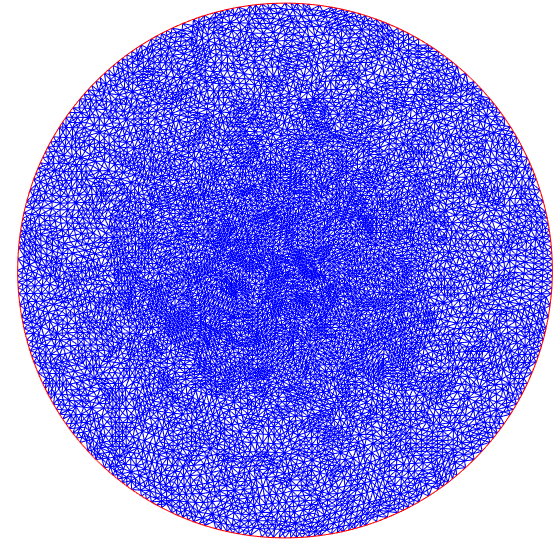}}\hspace{0.5cm}
    \subcaptionbox{$u_7$ \label{fig:UnitDiskP2U8TH}} {\includegraphics[scale=0.3]{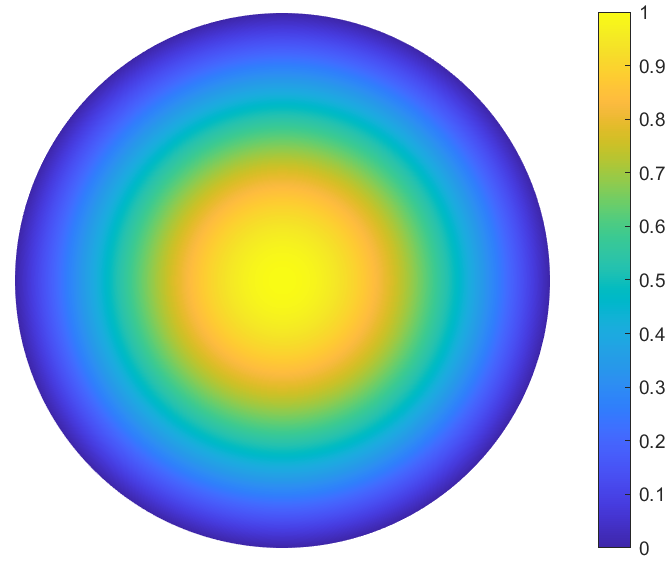}}\hspace{0.5cm}
    \subcaptionbox{reference \label{fig:UnitDiskP2U1THUniform}} {\includegraphics[scale=0.3]{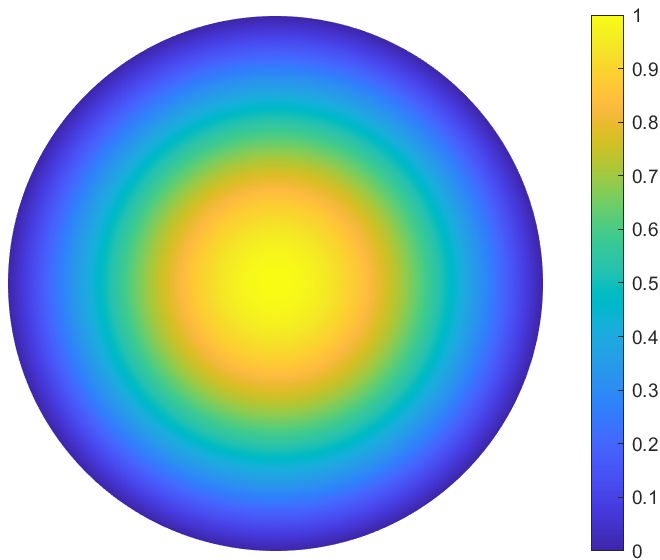}}
    \caption{Adaptive mesh refinement level with the number of vertices over each mesh as well as final computed first eigenfunctions by adaptive refinements and references for $p=1.1, 1.5, 2$ from top to bottom in Example \ref{example1}.}
    \label{UnitDiskP1to2MeshRefinement}
\end{figure}

\begin{figure}[htb!]
    \centering
    \subcaptionbox{$\mathcal{T}_1$ (1180)\label{fig:UnitDiskP2_5Mesh2TH}}
    {\includegraphics[scale=0.29]{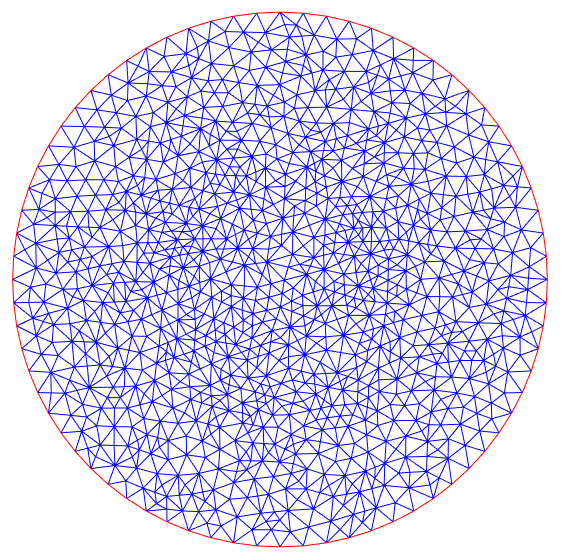}}\hspace{0.5cm}
    \subcaptionbox{$\mathcal{T}_3$ (3576)\label{fig:UnitDiskP2_5Mesh4TH}} {\includegraphics[scale=0.29]{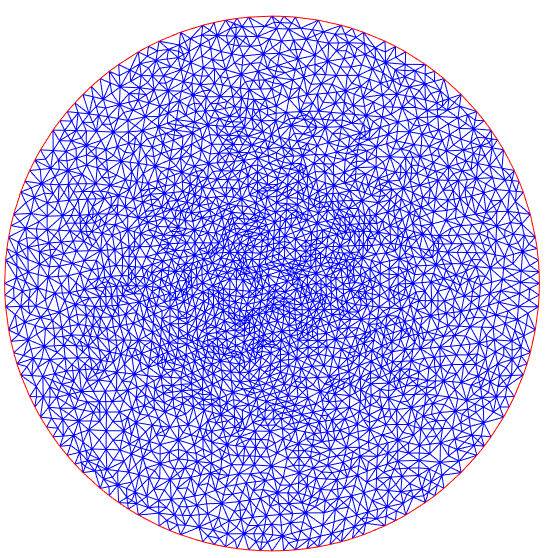}}\hspace{0.5cm}
    \subcaptionbox{$\mathcal{T}_7$ (32569)\label{fig:UnitDiskP2_5Mesh8TH}} {\includegraphics[scale=0.29]{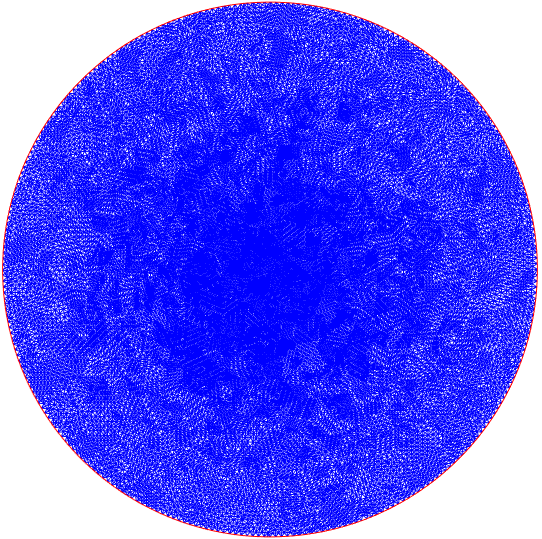}}\hspace{0.5cm}
    \subcaptionbox{$u_7$ \label{fig:UnitDiskP2_5U8TH}} {\includegraphics[scale=0.29]{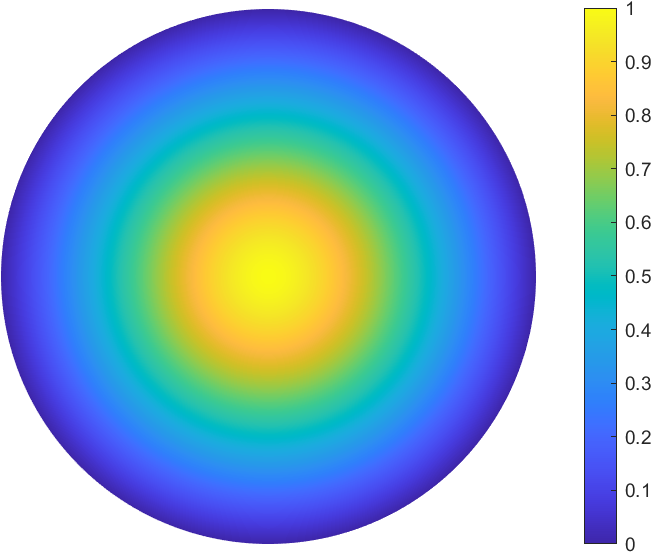}}
\\
    \subcaptionbox{$\mathcal{T}_2$ (1972)\label{fig:UnitDiskP3Mesh3TH}} {\includegraphics[scale=0.3]{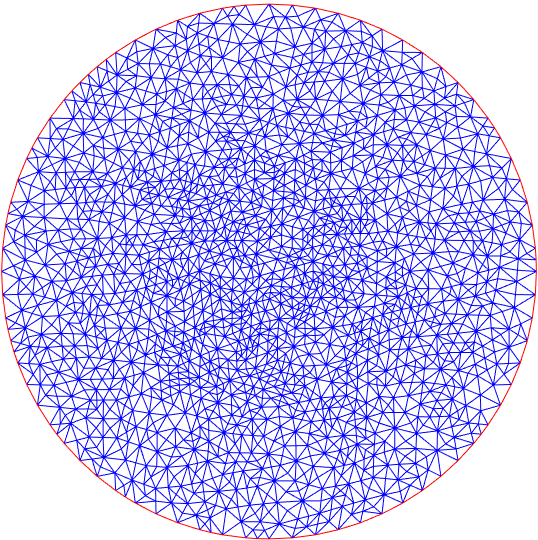}}\hspace{0.5cm}
    \subcaptionbox{$\mathcal{T}_4$ (5798) \label{fig:UnitDiskP3Mesh5TH}} {\includegraphics[scale=0.3]{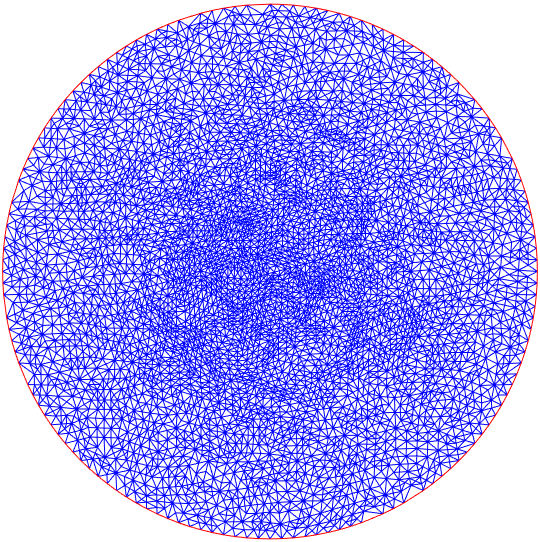}}\hspace{0.5cm}
    \subcaptionbox{$\mathcal{T}_6$ (17036) \label{fig:UnitDiskP3Mesh7TH}} {\includegraphics[scale=0.3]{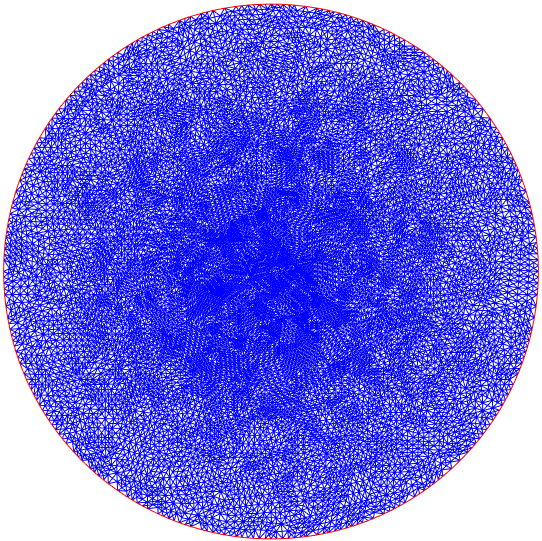}}\hspace{0.5cm}
    \subcaptionbox{$u_6$ \label{fig:UnitDiskP3U7TH}}
    {\includegraphics[scale=0.3]{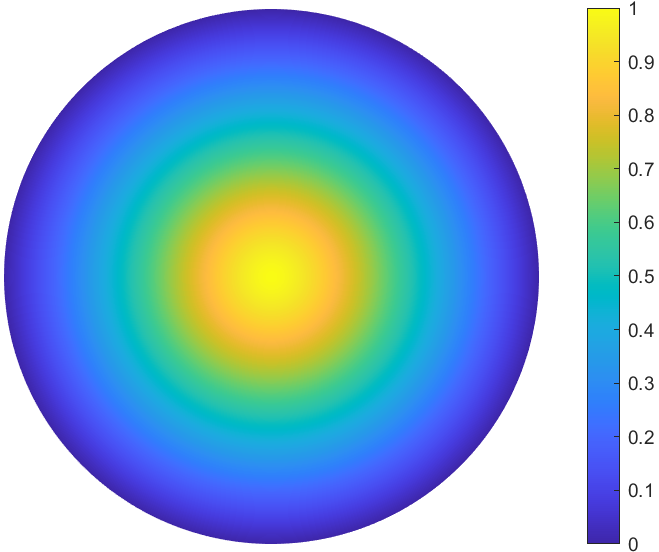}}
\\
    \subcaptionbox{$\mathcal{T}_2$ (1096) \label{fig:UnitDiskP4Mesh3TH}} {\includegraphics[scale=0.3]{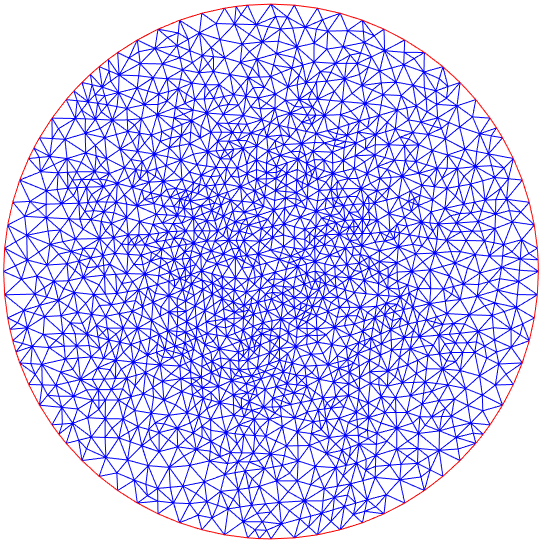}}\hspace{0.5cm}
    \subcaptionbox{$\mathcal{T}_5$ (8820) \label{fig:UnitDiskP4Mesh6TH}} {\includegraphics[scale=0.3]{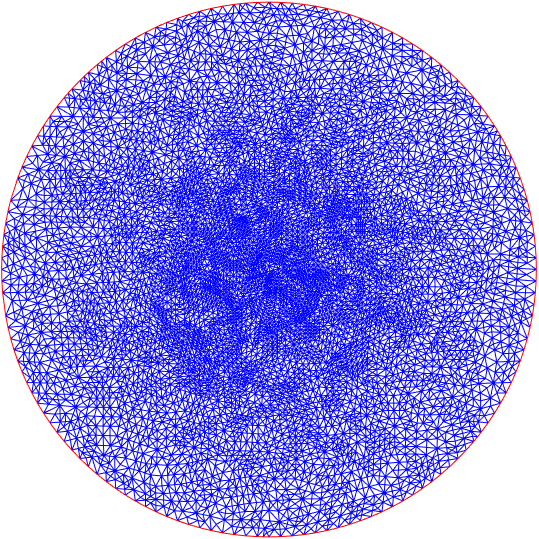}}\hspace{0.5cm}
    \subcaptionbox{$\mathcal{T}_7$ (25078) \label{fig:UnitDiskP4Mesh8TH}} {\includegraphics[scale=0.3]{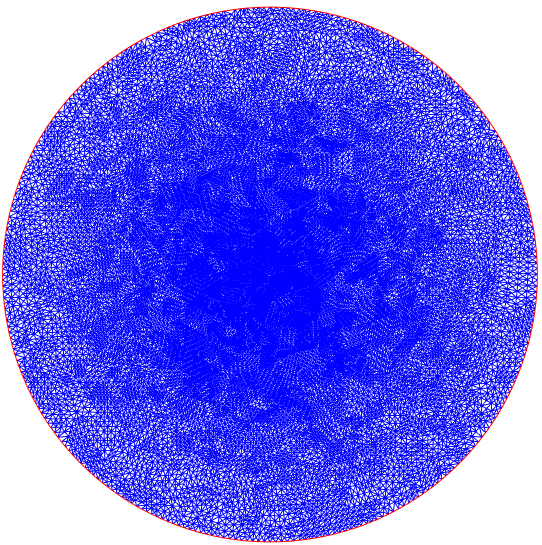}}\hspace{0.5cm}
    \subcaptionbox{$u_{9}$\label{fig:UnitDiskP4U8TH}} {\includegraphics[scale=0.3]{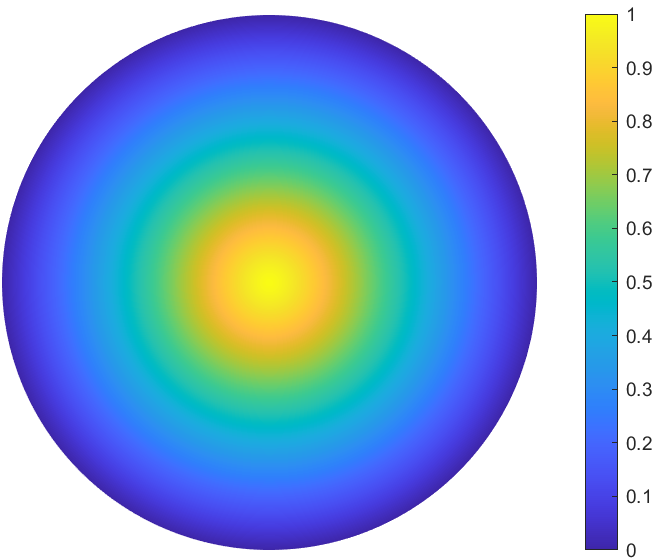}}
    \caption{Adaptive mesh refinement level with the number of vertices over each mesh as well as final computed first eigenfunctions for $p=2.5, 3, 4$ from top to bottom in Example \ref{example1}.}
    \label{UnitDiskP2to4MeshRefinement}
\end{figure}

\begin{figure}[htb!]
    \centering
    \subcaptionbox{$\mathcal{T}_5$ (5577) \label{fig:UnitDiskP10Mesh6TH}} {\includegraphics[scale=0.3]{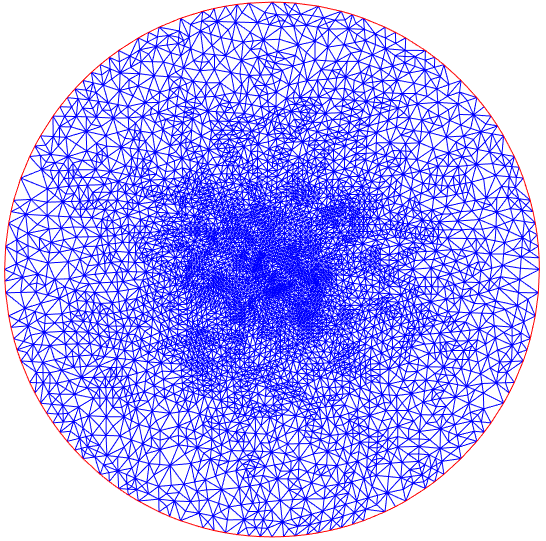}}\hspace{0.5cm}
   \subcaptionbox{$\mathcal{T}_8$ (23717) \label{fig:UnitDiskP10Mesh9TH}} {\includegraphics[scale=0.3]{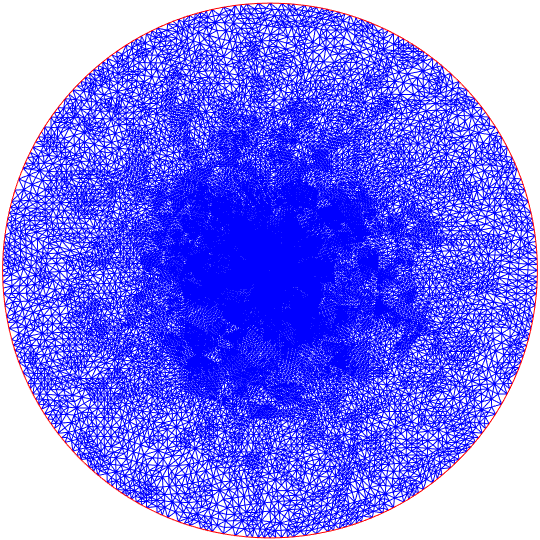}}\hspace{0.5cm}
    \subcaptionbox{$u_9$ \label{fig:UnitDiskP10U10TH}} {\includegraphics[scale=0.3]{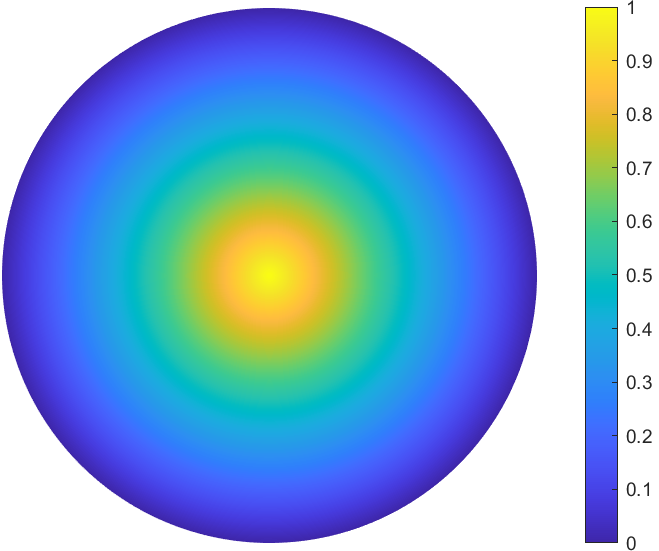}}\hspace{0.5cm}
    \subcaptionbox{reference \label{fig:UnitDiskP10U1THUniform}} {\includegraphics[scale=0.3]{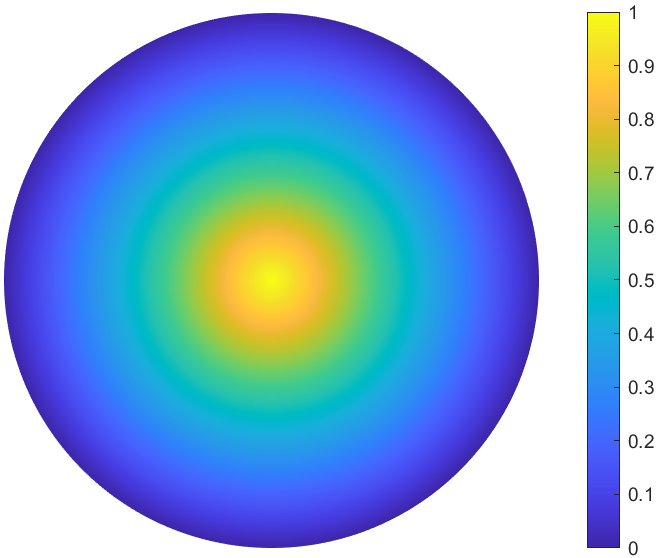}}
\\
    \subcaptionbox{$\mathcal{T}_5$ (4210) \label{fig:UnitDiskP20Mesh6TH}} {\includegraphics[scale=0.3]{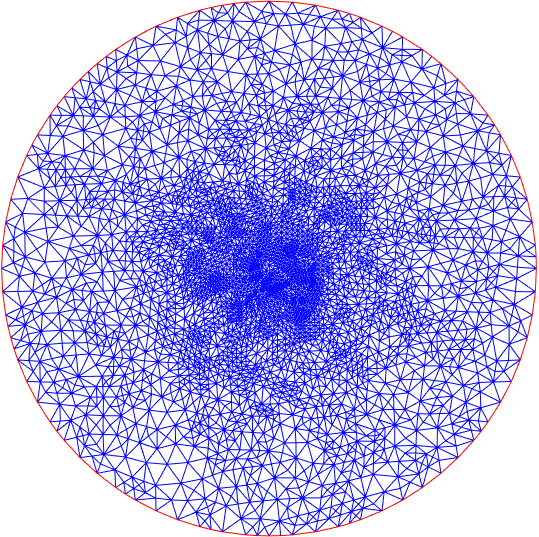}}\hspace{0.5cm}
    \subcaptionbox{$\mathcal{T}_8$ (15968) \label{fig:UnitDiskP20Mesh9TH}} {\includegraphics[scale=0.3]{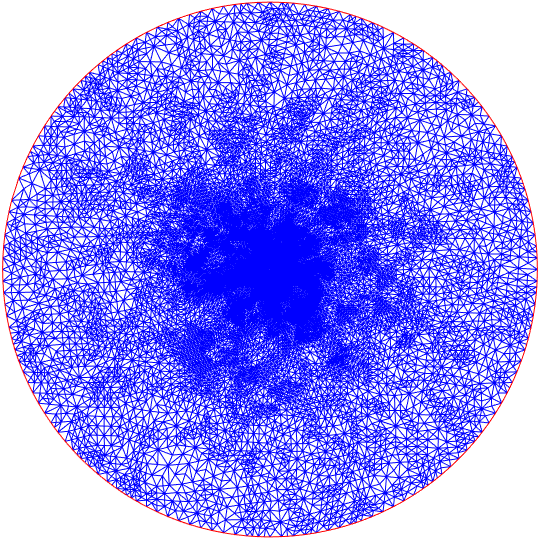}}\hspace{0.5cm}
    \subcaptionbox{$u_9$ \label{fig:UnitDiskP20U10TH}} {\includegraphics[scale=0.3]{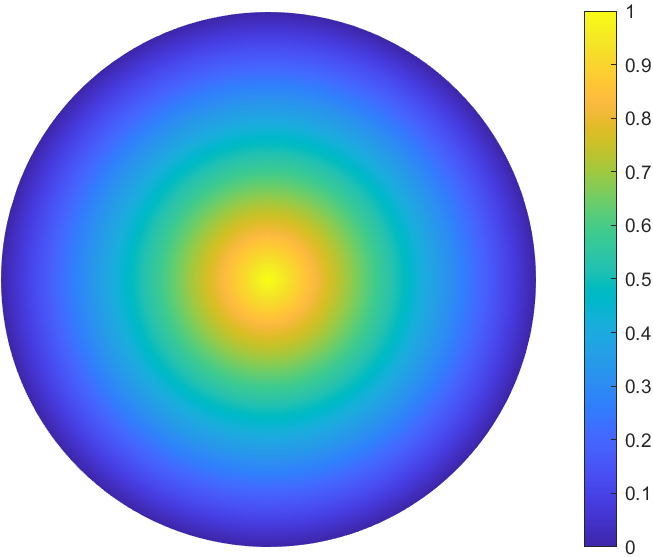}}\hspace{0.5cm}
    \subcaptionbox{reference \label{fig:UnitDiskP20U1THUniform}} {\includegraphics[scale=0.3]{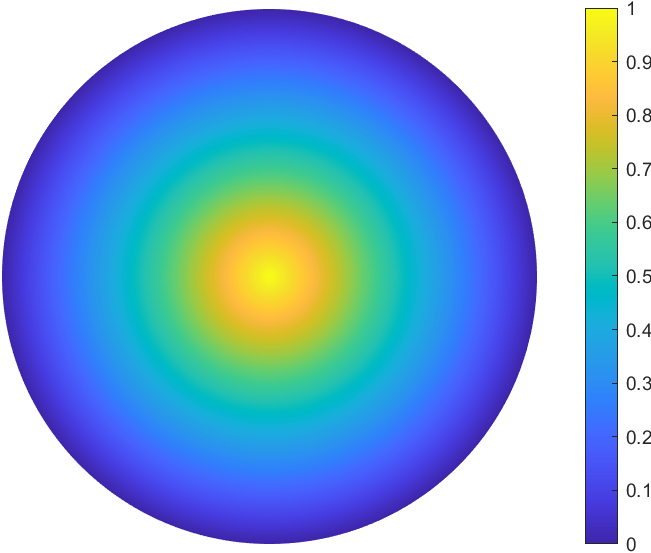}}
\\
    \subcaptionbox{$\mathcal{T}_5$ (3859) \label{fig:UnitDiskP30Mesh6TH}} {\includegraphics[scale=0.3]{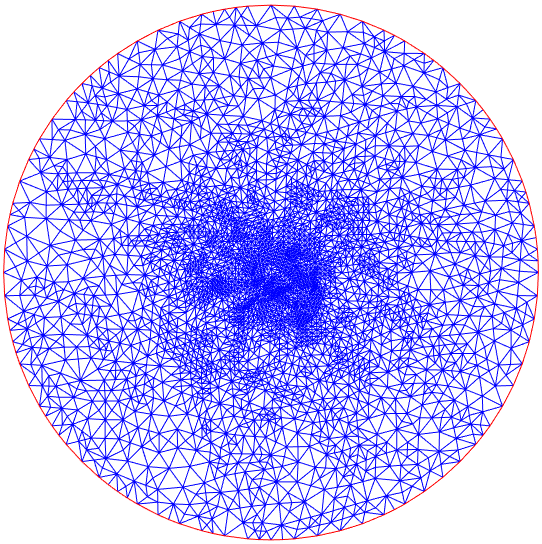}}\hspace{0.5cm}
    \subcaptionbox{$\mathcal{T}_8$ (13429) \label{fig:UnitDiskP30Mesh9TH}} {\includegraphics[scale=0.3]{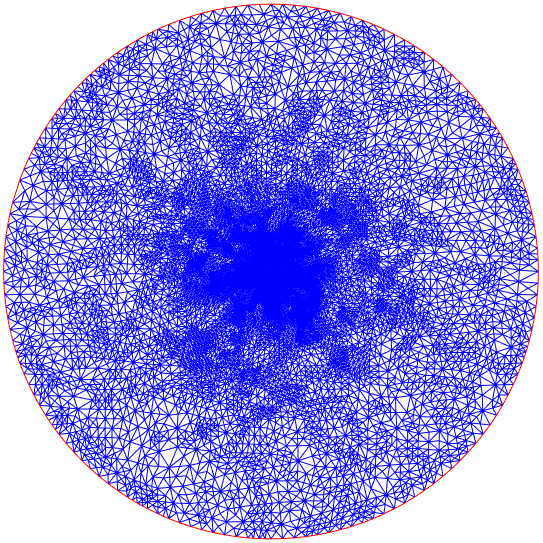}}\hspace{0.5cm}
    \subcaptionbox{$u_9$\label{fig:UnitDiskP30U10TH}} {\includegraphics[scale=0.3]{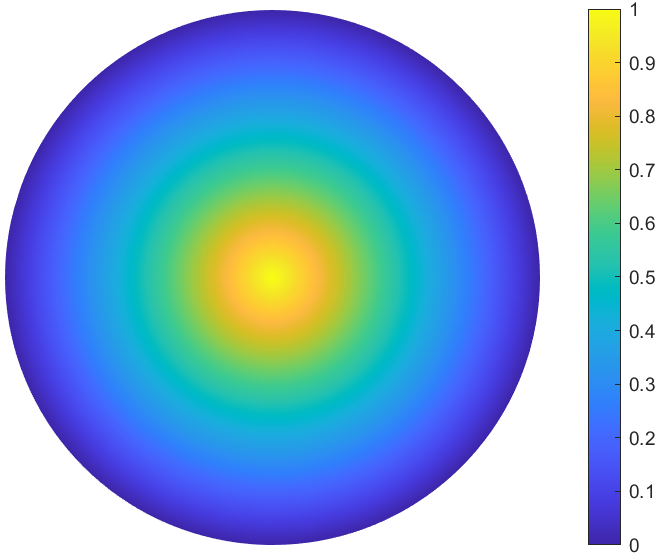}}\hspace{0.5cm}
    \subcaptionbox{reference \label{fig:UnitDiskP30U1THUniform}} {\includegraphics[scale=0.3]{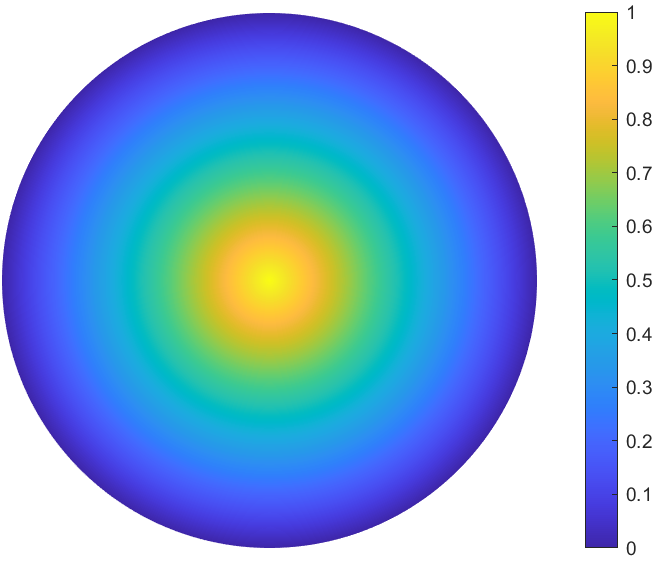}}
    \caption{Adaptive mesh refinement level with the number of vertices over each mesh as well as final computed first eigenfunctions by adaptive refinements and references for $p=10,20,30$ from top to bottom in Example \ref{example1}.}
    \label{UnitDiskP10to30MeshRefinement}
\end{figure}

\begin{example}[Unit Disk]\label{example1}
In the first example, the computational domain is a unit disk centered at the origin. Numerical experiments are implemented  for 9 different values of
\begin{align*}
p \in \{1.1, 1.5, 2, 2.5, 3, 4, 10, 20, 30\}.
\end{align*}
 Tolerance $\epsilon_K$ in Algorithm \ref{afem} is set to be $5\times 10^{-4}$ for $p=4$ and $10^{-4}$ for the rest cases, among which the adaptive algorithm terminates if the counter $k> K=9$ for $p=10, 20, 30$, while tolerance $\epsilon_M$ in Algorithm \ref{alg:alg1} is $5\times 10^{-4}$ for $p=1.1$, $10^{-5}$ for $p=1.5, 2, 2.5$ and $5\times 10^{-5}$ for the rest. 

For each $p$, approximate values, over meshes generated by the adaptive strategy, of the first eigenvalue of Problem \eqref{diff_eq} are provided in Table \ref{tab:UnitDiskEigenvalue} ($p=1.1,1.5,2,2.5,3$) and Table \ref{tab:UnitDiskEigenvalue102030} ($p=4,10,20,30$) respectively, where $k$ stands for the iteration number in Algorithm \ref{afem} while $\mu_k$ represents the approximate first eigenvalue over the $k$-th adaptive mesh. The results computed on fine meshes are listed in the last row as reference solutions. We observe that the sequence of computed eigenvalues $\{\mu_k\}_{k\geq0}$ is decreasing and approaches the reference solution for each $p$ as the adaptive mesh refinement level $k$ increases. This numerical observation confirms the convergence of Algorithm \ref{afem} as proved in Theorem \ref{thm_afem_conv}. As a benchmark, the first eigenvalue of Laplacian ($p=2$) on the unit disk is approximately 5.78321. Our adaptive eigenvalue $\mu_7=5.78402$ as listed in the fifth column ($p=2$) of Table \ref{tab:UnitDiskEigenvalue} attains a 0.014\% relative precision. One sees from Table \ref{tab:UnitDiskEigenvalue102030} that the approximate eigenvalue $\mu_7=14.6713$ for $p=4$ is smaller than the reference solution produced on a very fine mesh, which implies more accuracy with fewer degrees of freedom of Algorithm \ref{afem}.

Figures \ref{UnitDiskP1to2MeshRefinement}-\ref{UnitDiskP10to30MeshRefinement} depict a selection of adaptive meshes generated by Algorithm \ref{afem} and computed first eigenfunctions over the finest adaptive meshes. For comparison, computed first eigenfunctions over very fine meshes are displayed in the last column of Figures \ref{UnitDiskP1to2MeshRefinement} and \ref{UnitDiskP10to30MeshRefinement}. We observe that
as $p$ becomes larger, local refinements are performed in the vicinity of the origin. Moreover, as can be seen from the penultimate column of Figure \ref{UnitDiskP1to2MeshRefinement} (Figures \ref{fig:UnitDiskP1_1U4TH}-\ref{fig:UnitDiskP2U8TH}), the last column of Figure \ref{UnitDiskP2to4MeshRefinement} (Figures \ref{fig:UnitDiskP2_5U8TH}-\ref{fig:UnitDiskP4U8TH}) and the penultimate column of Figure \ref{UnitDiskP10to30MeshRefinement} (Figures \ref{fig:UnitDiskP10U10TH}-\ref{fig:UnitDiskP30U10TH}), the asymptotic behaviour of adaptively computed first eigenfunctions confirms two assertions in \cite{KawohlFridman:2003} and \cite{JuutineLindqvistManfredi:1999} that the first $L^\infty(\Omega)$-normalized eigenfunction of \eqref{diff_eq} converges to 1, the characteristic function of the unit ball, as $p\to1^+$ and to the distance function to the boundary, $1-\sqrt{x_1^2+x_2^2}$, as $p\to\infty$ respectively. This  also explains the transition of additional mesh refinements from the whole domain to the center. 

\end{example}

\begin{table}[hbt!] 
\centering
\caption{Quantitative result for $p\in\{1.1, 1.5, 2, 2.5, 3\}$ in Example \ref{example1}: the number of adaptive loops, the number of vertices and the computed first eigenvalue.}\label{tab:UnitDiskEigenvalue}
\resizebox{\textwidth}{!}{
\begin{tabular}{c cc cc cc cc cc}
    \toprule
    \multirow{2}{*}{$k$}&
    \multicolumn{2}{c}{$p=1.1$}&
    \multicolumn{2}{c}{$p=1.5$}&\multicolumn{2}{c}{$p=2$}&
    \multicolumn{2}{c}{$p=2.5$}&
    \multicolumn{2}{c}{$p=3$}\\
    \cmidrule(lr){2-11} 
     ~ & vertices & $\mu_k$ & vertices & $\mu_k$ & vertices & $\mu_k$ & vertices & $\mu_k$ & vertices & $\mu_k$ \\
    \midrule
        0 &   682 & 2.57807 & 682   & 4.02918 & 682   & 5.80632  &  682  & 7.75284  & 682  & 9.90463\\
        1 &   1055 & 2.57749& 1213  & 4.02732 & 1196  & 5.80145  & 1180  & 7.74343  & 2214 & 9.89018\\
        2 &  1637 & 2.57671& 2188   & 4.02431 & 2091  & 5.79601  & 2052  & 7.73324  & 1972 & 9.86859\\
        3 &  2498 & 2.57662& 3927   & 4.02137 & 3721  & 5.79006  & 3576  & 7.72451  & 3347 & 9.85650\\
        4 &       &        & 6935   & 4.01992 & 6532  & 5.78718  & 6236  & 7.71776  & 5798 & 9.85249\\
        5 &       &        & 12469  & 4.01888 & 11639 & 5.78542  & 10850 & 7.71546  & 9946 & 9.84022\\
        6 &      &         & 22277  & 4.01842 & 20538 & 5.78440  & 18843 & 7.71253  & 17036& 9.83989\\
        7 &      &         & 39173  & 4.01822 & 35714 & 5.78402  & 32569 & 7.71199 \\
    \midrule
    reference & 24505 & 2.56642 & 65130 & 4.01790 & 65130 & 5.78342 & 65130 & 7.71122 & 65130 & 9.83481\\
    \bottomrule
    \end{tabular}}
\end{table}

\begin{table}[hbt!]
\centering
\caption{Quantitative result for $p\in\{4, 10, 20, 30\}$ Example \ref{example1}: the number of adaptive loops, the number of vertices and the computed first eigenvalue.}\label{tab:UnitDiskEigenvalue102030}
\begin{tabular}{c cc cc cc cc}
    \toprule
    \multirow{2}{*}{$k$}&
    \multicolumn{2}{c}{$p=4$}&
    \multicolumn{2}{c}{$p=10$}&\multicolumn{2}{c}{$p=20$}&
    \multicolumn{2}{c}{$p=30$}\\
    \cmidrule(lr){2-9}
     ~ & vertices & $\mu_k$ & vertices & $\mu_k$ & vertices & $\mu_k$ & vertices & $\mu_k$ \\
    \midrule
        0 & 682   & 14.8676 &  682  & 65.5562 &   682 & 270.463 & 682   & 792.061  \\
        1 & 1096  & 14.8350  &  936  & 64.0533 &   942 & 246.373 & 909   & 646.565  \\
        2 & 1838  & 17.7682 & 1400  & 62.6093 &  1269 & 224.980 & 1260  & 538.993  \\
        3 & 3041  & 14.7447 & 2184  & 62.0110 &  1854 & 218.470 & 1801  & 506.405  \\
        4 & 5179  & 14.7279 & 3460  & 61.7149 &  2765 & 210.819 & 2625  & 466.290  \\
        5 & 8820  & 14.7044 & 7935  & 61.4366 &  4210 & 207.801 & 3859  & 449.793  \\
        6 & 14818 & 14.6769 & 5577  & 61.1423 &  6515 & 205.064 & 5768  & 440.326  \\
        7 & 25078 & 14.6713 & 8991  & 60.8983 & 10146 & 203.304  & 8763 & 431.970  \\
        8 &       &         & 23717 & 60.8331 & 15968 & 202.576 & 13429 & 428.450   \\
    	9 &       &         & 38698 & 60.6828 & 25387 & 201.803 & 20911 & 425.581   \\
    \midrule
    reference & 86001 & 14.6833 & 76492 & 60.6684 & 76492  &  201.691  & 76492  & 423.168 \\
    \bottomrule
    \end{tabular}
\end{table}


\begin{figure}[htb!]
    \centering
    \subcaptionbox{$\mathcal{T}_5$ (4098) \label{fig:UnitSquareP1_2Mesh6TH}} {\includegraphics[scale=0.3]{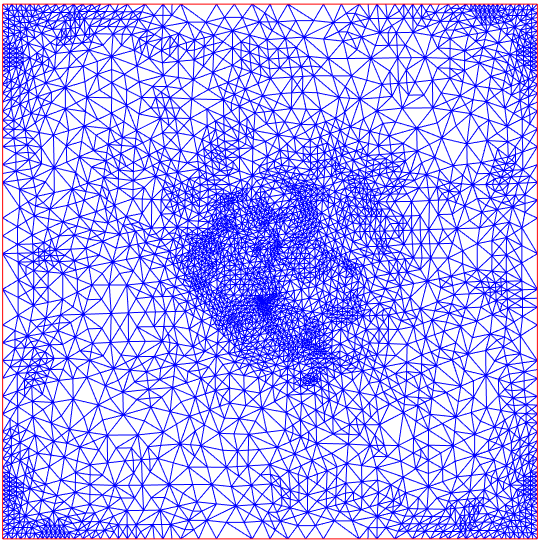}}\hspace{0.5cm}
    \subcaptionbox{$\mathcal{T}_7$ (9596) \label{fig:UnitSquareP1_2Mesh8TH}} {\includegraphics[scale=0.3]{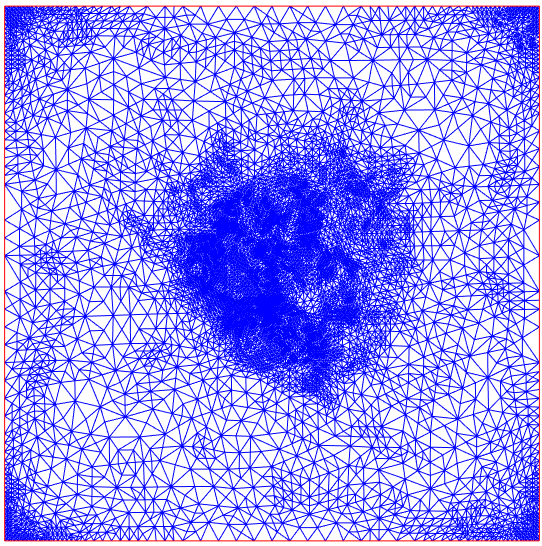}}\hspace{0.5cm}
    \subcaptionbox{$u_7$ \label{fig:UnitSquareP1_2U9TH}} {\includegraphics[scale=0.3]{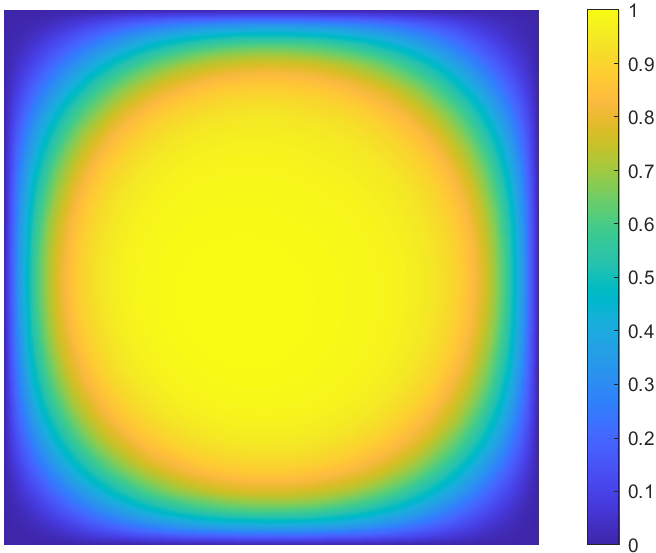}}\hspace{0.5cm}
    \subcaptionbox{reference \label{fig:UnitSquareP1_2U1THUniform}} {\includegraphics[scale=0.3]{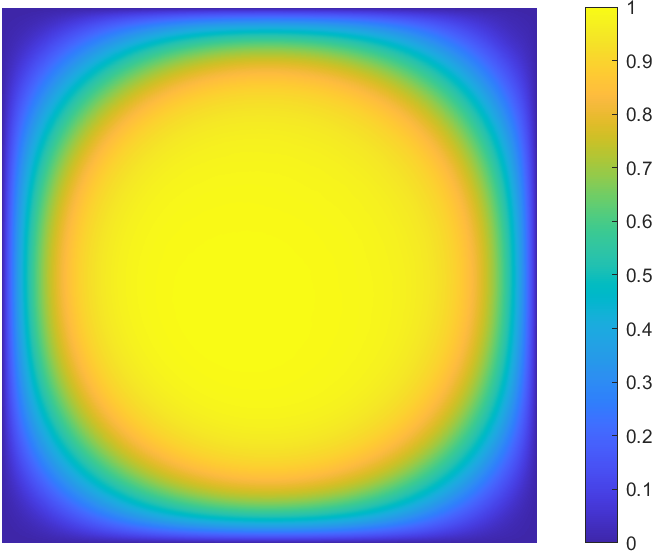}}
\\
    \subcaptionbox{$\mathcal{T}_5$ (5071) \label{fig:UnitSquareP1_5Mesh6TH}} {\includegraphics[scale=0.295]{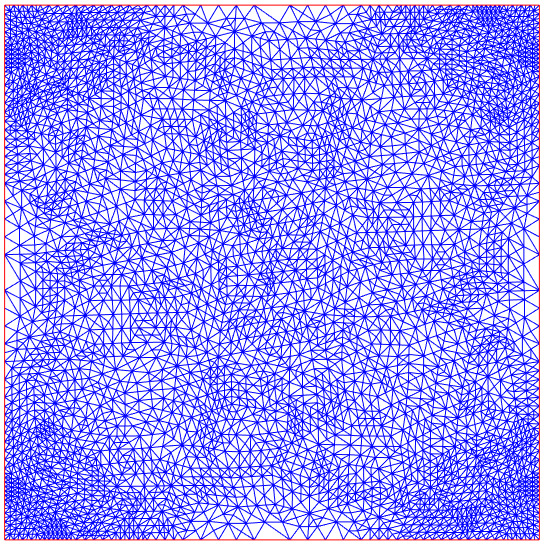}}\hspace{0.5cm}
    \subcaptionbox{$\mathcal{T}_8$ (26287) \label{fig:UnitSquareP1_5Mesh9TH}} {\includegraphics[scale=0.295]{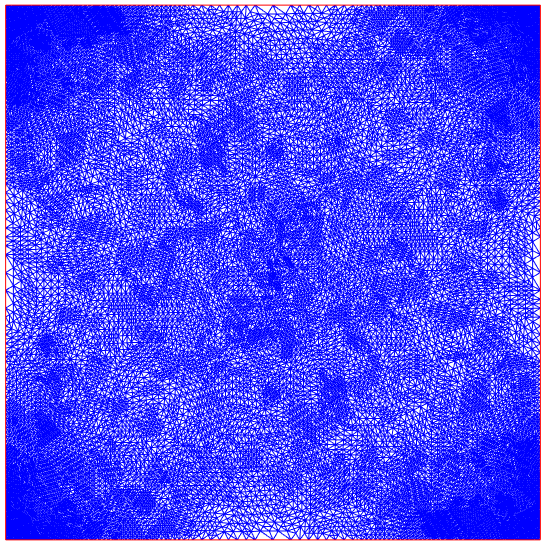}}\hspace{0.5cm}
    \subcaptionbox{$u_8$ \label{fig:UnitSquareP1_5U9TH}} {\includegraphics[scale=0.295]{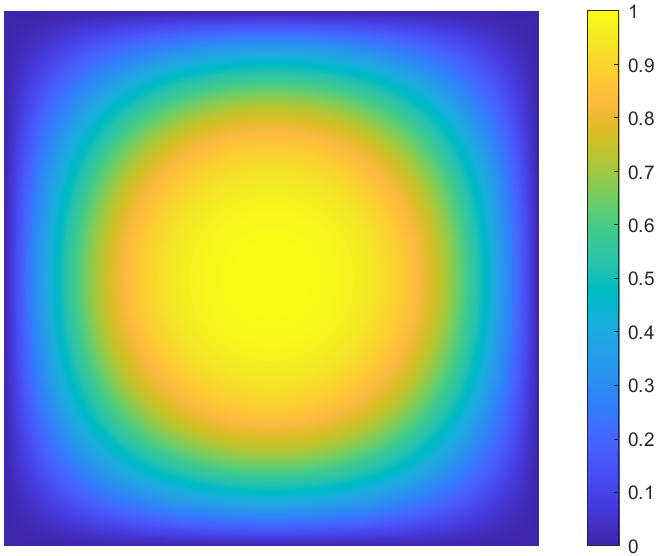}}\hspace{0.5cm}
    \subcaptionbox{reference \label{fig:UnitSquareP1_5U1THUniform}} {\includegraphics[scale=0.295]{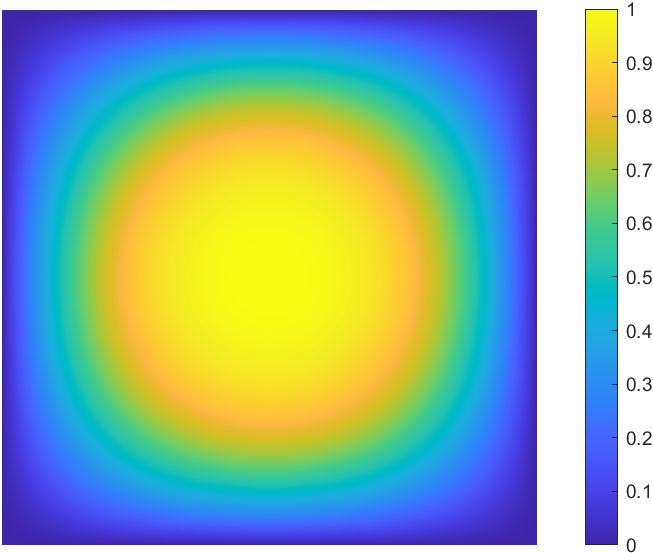}}
\\
    \subcaptionbox{$\mathcal{T}_6$ (9908) \label{fig:UnitSquareP2Mesh7TH}} {\includegraphics[scale=0.3]{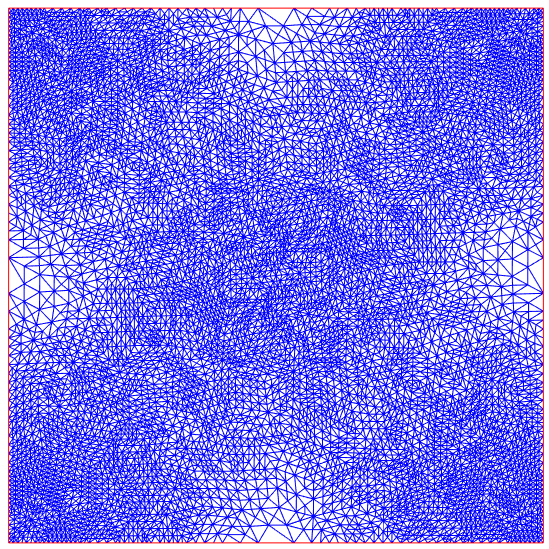}}\hspace{0.5cm}
    \subcaptionbox{$\mathcal{T}_8$ (29834) \label{fig:UnitSquareP2Mesh9TH}} {\includegraphics[scale=0.3]{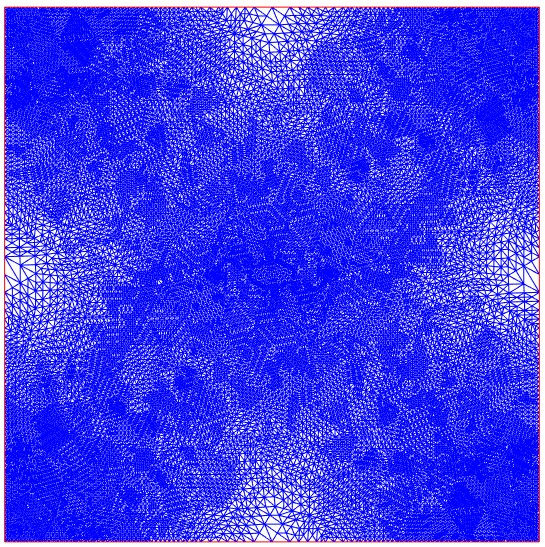}}\hspace{0.5cm}
    \subcaptionbox{$u_9$ \label{fig:UnitSquareP2U10TH}}
    {\includegraphics[scale=0.3]{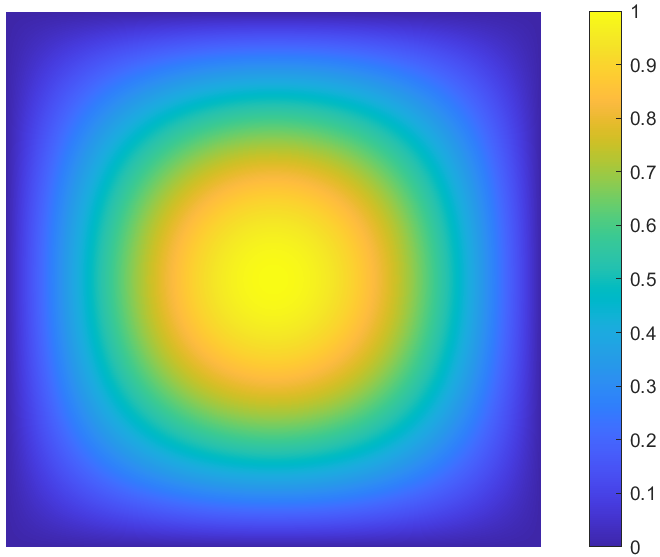}}\hspace{0.5cm}
    \subcaptionbox{reference \label{fig:UnitSquareP2U1THUniform}} {\includegraphics[scale=0.3]{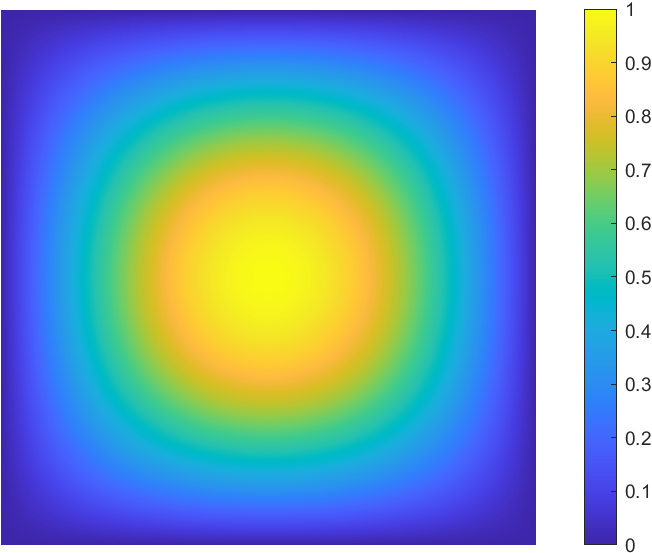}}
    \caption{Adaptive mesh refinement level with the number of vertices over each mesh as well as final computed first eigenfunctions by adaptive refinements and references for $p=1.2, 1.5, 2$ from top to bottom in Example \ref{example2}.}   
    \label{UnitSquareP1to2MeshRefinement}
\end{figure}

\begin{figure}[htb!]
    \centering
    \subcaptionbox{$\mathcal{T}_2$ (1098) \label{fig:UnitSquareP2_5Mesh3TH}} {\includegraphics[scale=0.3]{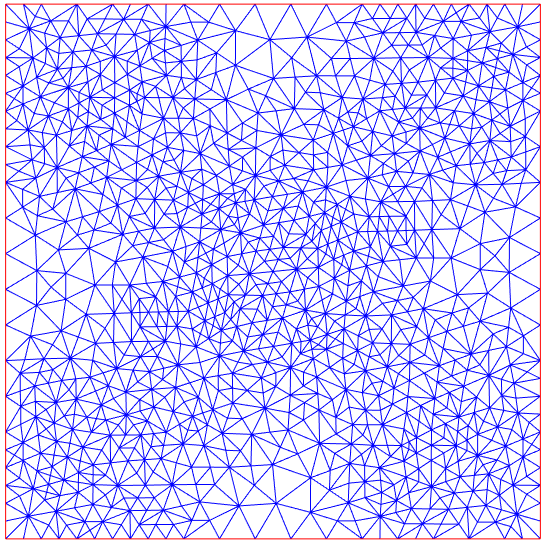}}\hspace{0.5cm}
    \subcaptionbox{$\mathcal{T}_6$ (9928) \label{fig:UnitSquareP2_5Mesh7TH}} {\includegraphics[scale=0.3]{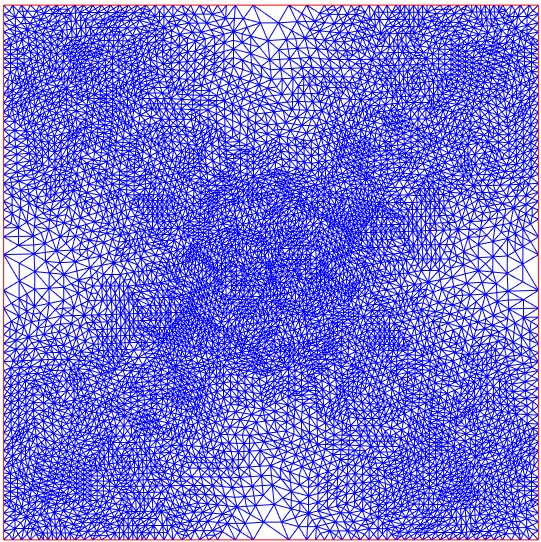}}\hspace{0.5cm}
    \subcaptionbox{$\mathcal{T}_8$ (28760) \label{fig:UnitSquareP2_5Mesh9TH}} {\includegraphics[scale=0.3]{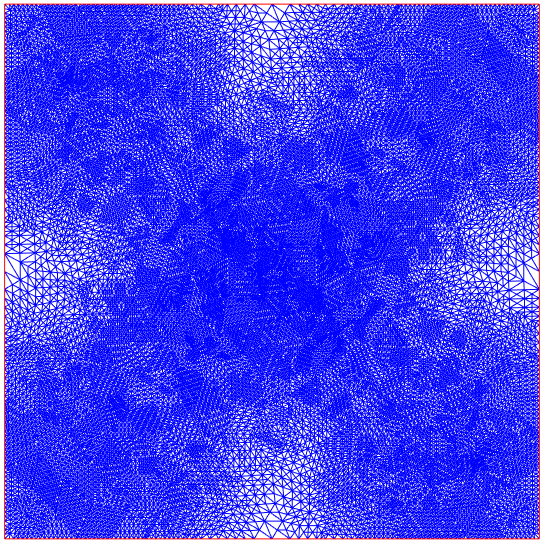}}\hspace{0.5cm}
    \subcaptionbox{$u_9$ \label{fig:UnitSquareP2_5U10TH}}
    {\includegraphics[scale=0.3]{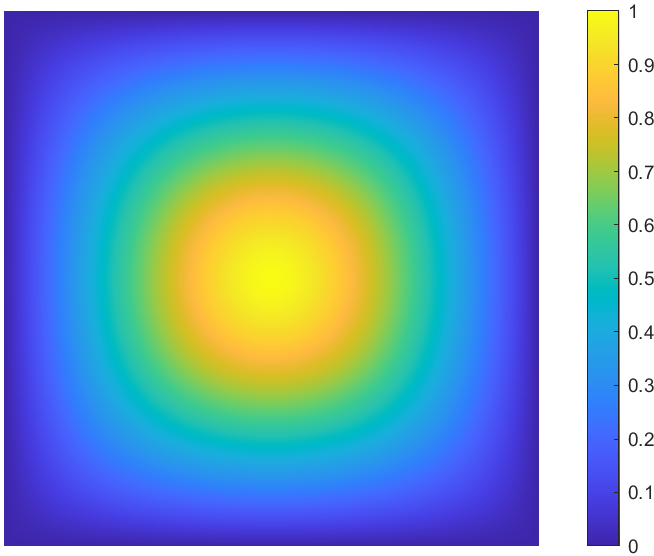}}
    \\
    \subcaptionbox{$\mathcal{T}_2$ (1056) \label{fig:UnitSquareP3Mesh3TH}} {\includegraphics[scale=0.295]{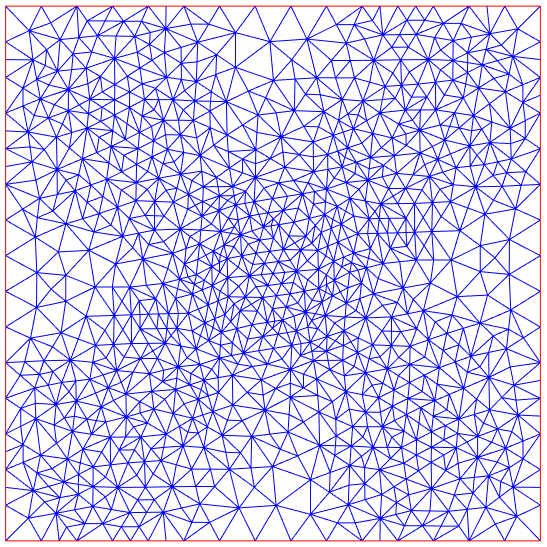}}\hspace{0.5cm}
    \subcaptionbox{$\mathcal{T}_6$ (9289)\label{fig:UnitSquareP3Mesh7TH}} {\includegraphics[scale=0.295]{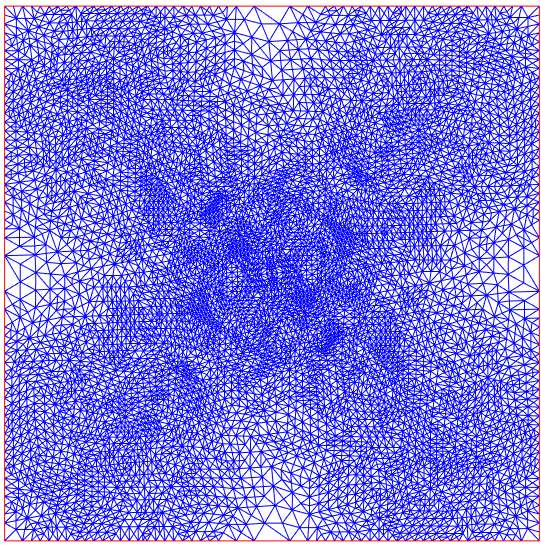}}\hspace{0.5cm}
    \subcaptionbox{$\mathcal{T}_7$ (15780) \label{fig:UnitSquareP3Mesh8TH}} {\includegraphics[scale=0.295]{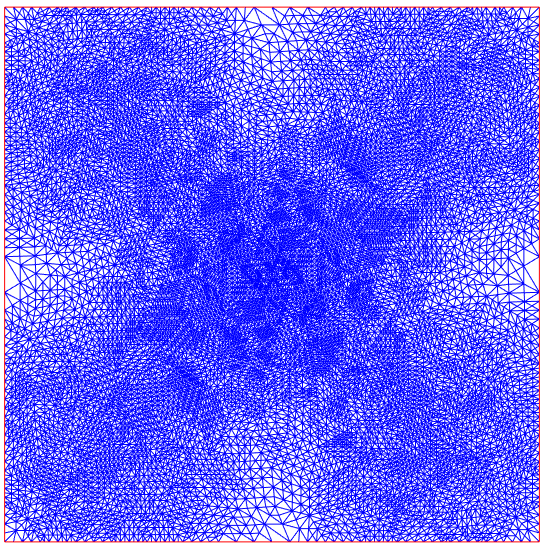}}\hspace{0.5cm}
    \subcaptionbox{$u_7$ \label{fig:UnitSquareP3U8TH}}
    {\includegraphics[scale=0.295]{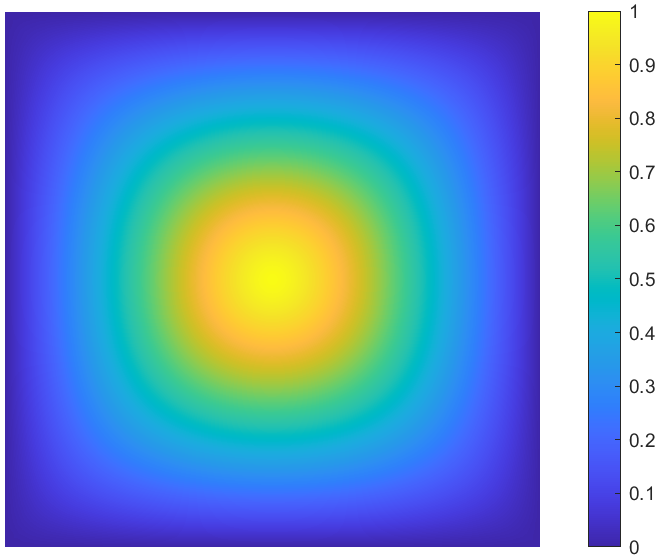}}
    \\
    \subcaptionbox{$\mathcal{T}_2$ (1002) \label{fig:UnitSquareP4Mesh3TH}} {\includegraphics[scale=0.3]{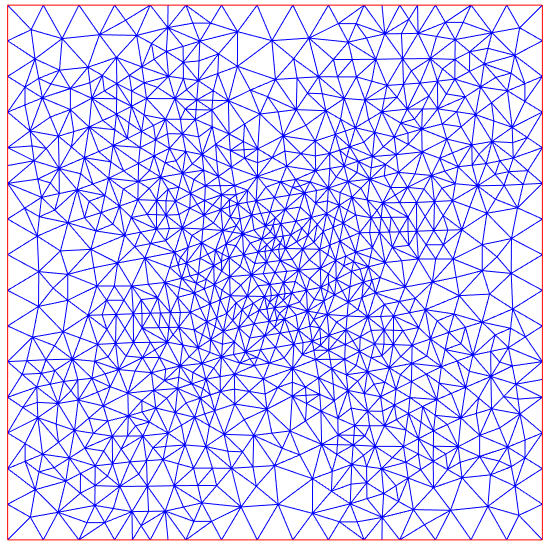}}\hspace{0.5cm}
    \subcaptionbox{$\mathcal{T}_5$ (4850) \label{fig:UnitSquareP4Mesh6TH}} {\includegraphics[scale=0.3]{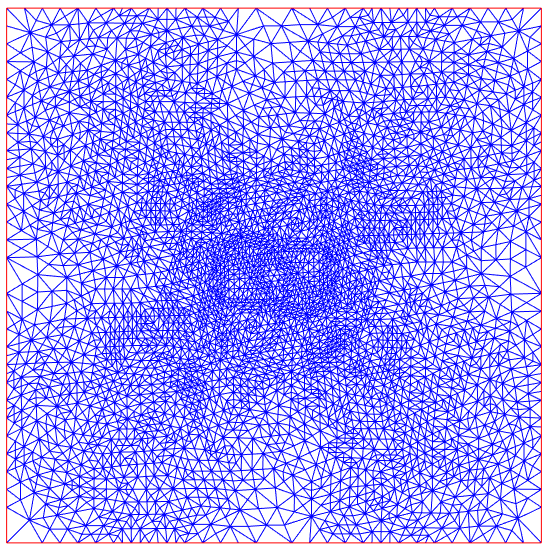}}\hspace{0.5cm}
    \subcaptionbox{$\mathcal{T}_6$ (8108) \label{fig:UnitSquareP4Mesh7TH}} {\includegraphics[scale=0.3]{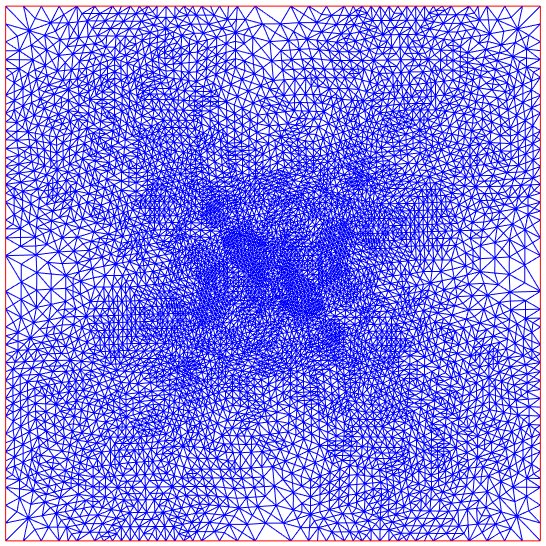}}\hspace{0.5cm}
    \subcaptionbox{$u_6$ \label{fig:UnitSquareP4U7TH}}
    {\includegraphics[scale=0.3]{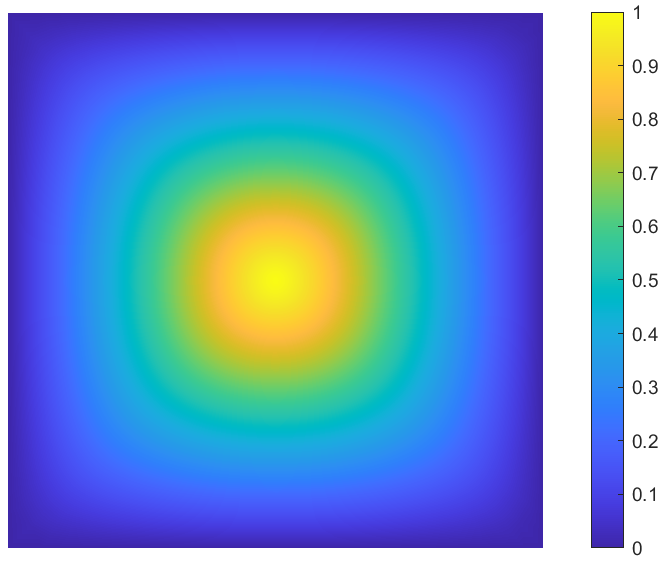}}
    \caption{Adaptive mesh refinement level with the number of vertices over each mesh as well as final computed first eigenfunctions by adaptive refinements for $p=2.5, 3, 4$ from top to bottom in Example \ref{example2}.}
    \label{UnitSquareP2to4MeshRefinement}
\end{figure}

\begin{figure}[htb!]
    \centering
    \subcaptionbox{$\mathcal{T}_6$ (4359) \label{fig:UnitSquareP10Mesh7TH}} {\includegraphics[scale=0.3]{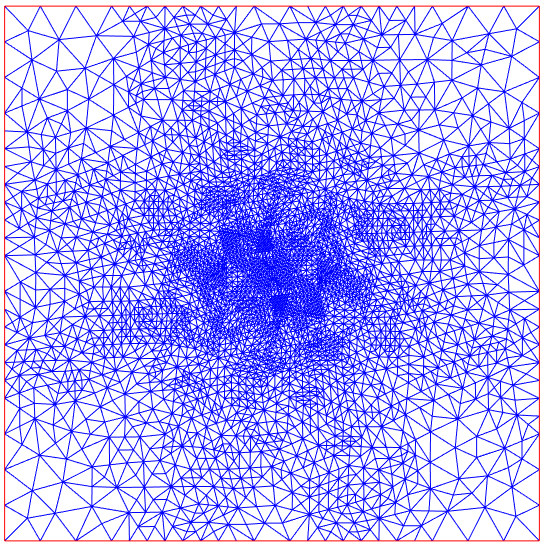}}\hspace{0.5cm}
    \subcaptionbox{$\mathcal{T}_{9}$ (19286) \label{fig:UnitSquareP10Mesh10TH}} {\includegraphics[scale=0.3]{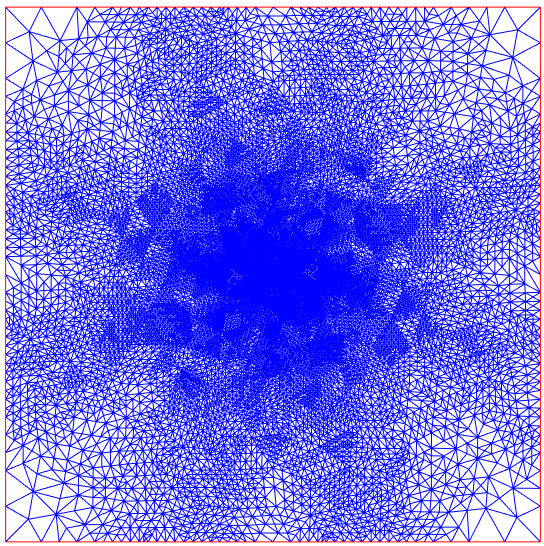}}\hspace{0.5cm}
    \subcaptionbox{$u_{9}$ \label{fig:UnitSquareP10U10TH}} {\includegraphics[scale=0.3]{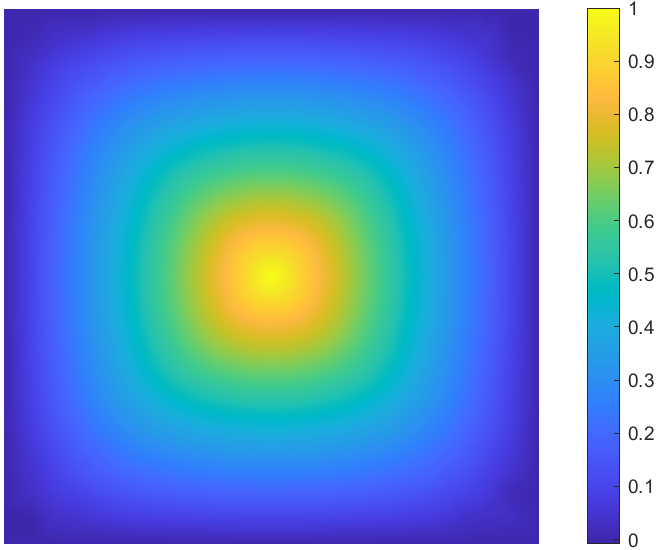}}\hspace{0.5cm}
    \subcaptionbox{reference \label{fig:UnitSquareP10U1THUniform}}{\includegraphics[scale=0.3]{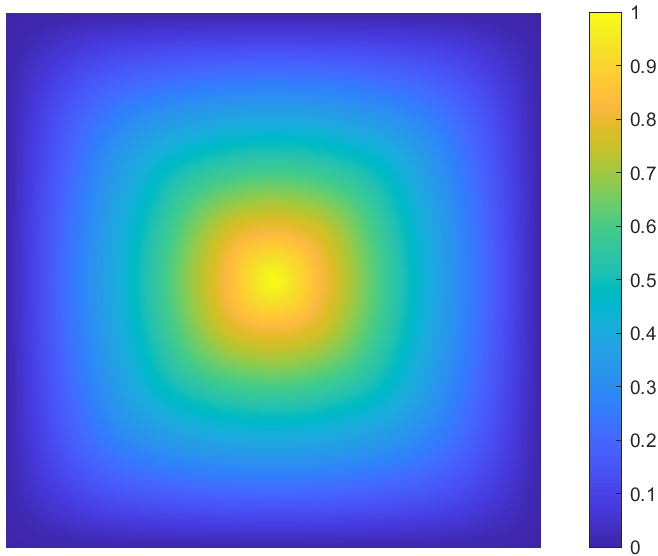}}
\\
    \subcaptionbox{$\mathcal{T}_6$ (3089) \label{fig:UnitSquareP20Mesh7TH}} {\includegraphics[scale=0.3]{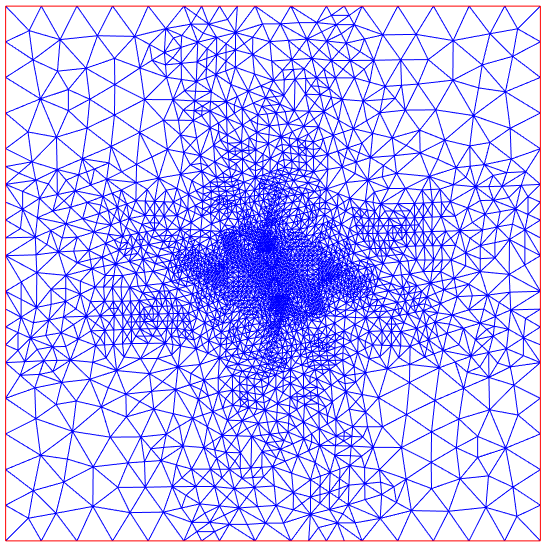}}\hspace{0.5cm}
    \subcaptionbox{$\mathcal{T}_{9}$ (12756) \label{fig:UnitSquareP20Mesh10TH}} {\includegraphics[scale=0.305]{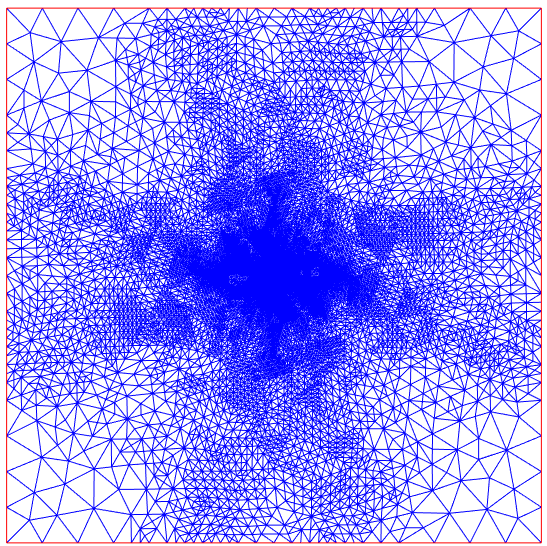}}\hspace{0.5cm}
    \subcaptionbox{$u_{9}$ \label{fig:UnitSquareP20U10TH}}{\includegraphics[scale=0.3]{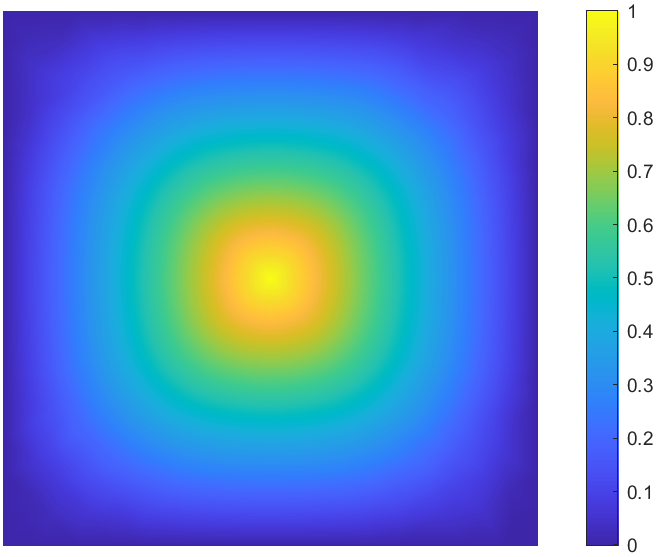}}\hspace{0.5cm}
	\subcaptionbox{reference \label{fig:UnitSquareP20U1THUniform}}{\includegraphics[scale=0.3]{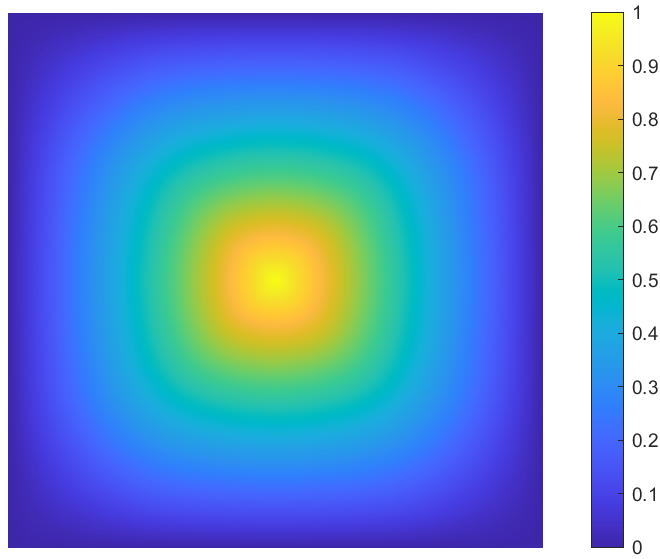}}
\\
    \subcaptionbox{$\mathcal{T}_6$ (2617) \label{fig:UnitSquareP30Mesh7TH}} {\includegraphics[scale=0.3]{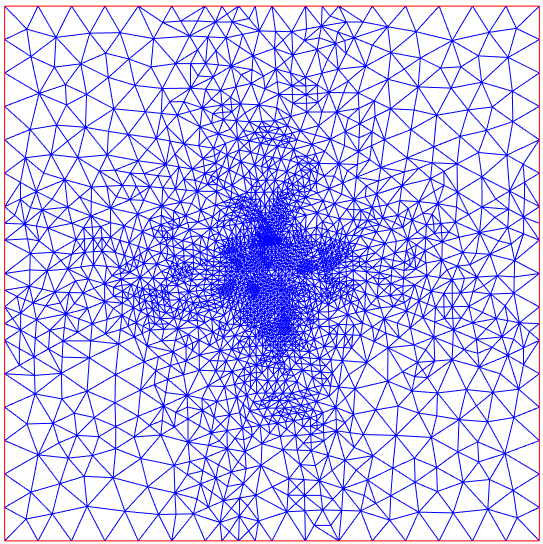}}\hspace{0.5cm}
    \subcaptionbox{$\mathcal{T}_{9}$ (9026) \label{fig:UnitSquareP30Mesh10TH}} {\includegraphics[scale=0.3]{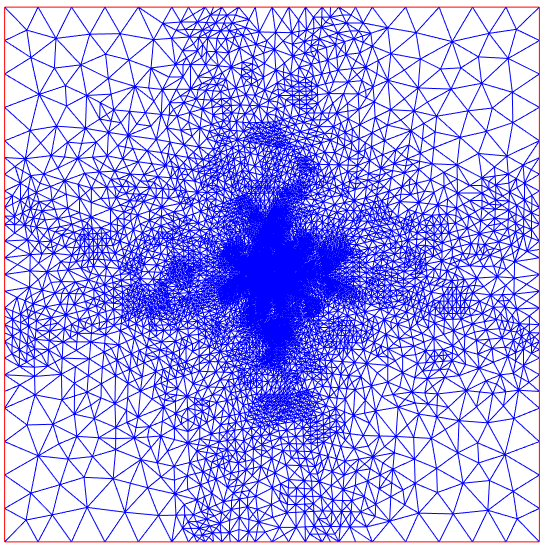}}\hspace{0.5cm}
    \subcaptionbox{$u_{10}$ \label{fig:UnitSquareP30U11TH}}{\includegraphics[scale=0.3]{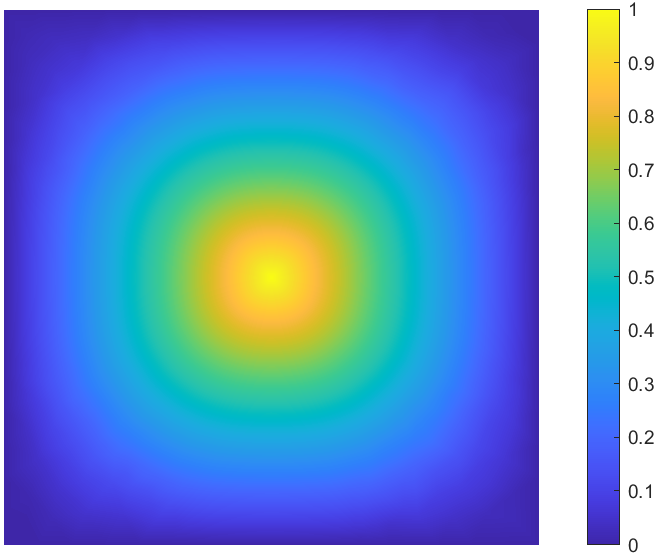}}\hspace{0.5cm}
    \subcaptionbox{reference \label{fig:UnitSquareP30U1THUniform}}{\includegraphics[scale=0.3]{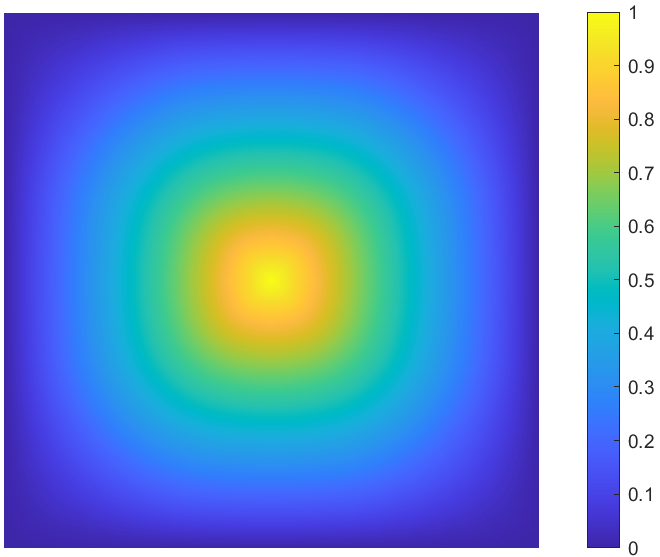}}
    \caption{Adaptive mesh refinement level with the number of vertices over each mesh as well as final computed first eigenfunctions by adaptive refinements and references for $p=10, 20, 30$ from top to bottom in Example \ref{example2}.}
    \label{UnitSquareP10to30MeshRefinement}
\end{figure}

\begin{example}[Unit Square]\label{example2}

We next consider the first eigenpair of \eqref{diff_eq} in a unit square $(0,1)^2$ with $p = 1.2, 1.5, 2, 2.5, 3, 4, 10, 20, 30$. 
In this example, tolerance $\epsilon_K$ in Algorithm \ref{afem} is chosen as $5\times 10^{-4}$ for $p=3$ and $10^{-4}$ for the rest $p$ with a maximum adaptive refinement number $K=9$ and $K=10$ imposed for $p=10, 20$ and $p = 30$ respectively while $\epsilon_M=10^{-5}$ in Algorithm \ref{alg:alg1} except for $p=1.2$ and $p=4$, in both of which $\epsilon_M=5\times 10^{-5}$. 

Table \ref{tab:UnitSquareEigenvalue} and Table \ref{tab:UnitSquareEigenvalue102030} contain all computed first eigenvalues over adaptively generated meshes and very fine meshes as reference solutions. As in the previous example, the sequence of adaptive eigenvalues for each $p$ strictly decreases to the reference solution. Noting that the exact first eigenvalue of Laplacian ($p=2$) on the unit square is $2\pi^2\approx 19.7392$, we see that our result $\mu_9 = 19.7410$ has an relative error less than $10^{-4}$.

Sequences of adaptively generated meshes and computed first eigenfunctions over the final adaptive meshes as well as very fine meshes are displayed in Figures \ref{UnitSquareP1to2MeshRefinement}-\ref{UnitSquareP10to30MeshRefinement}. The mesh is essentially refined in the vicinity of the center and around four corners as before for small $p$ and then near two crossing diagonals for moderate $p>2$. Finally, additional local refinements take place around the center for large $p$. This phenomenon is due to  concentration of large gradients of numerical solutions as displayed in the penultimate column of Figure \ref{UnitSquareP1to2MeshRefinement} (Figures \ref{fig:UnitSquareP1_2U9TH}-\ref{fig:UnitSquareP2U10TH}), the last column of Figure \ref{UnitSquareP2to4MeshRefinement} (Figures \ref{fig:UnitSquareP2_5U10TH}-\ref{fig:UnitSquareP4U7TH}) and the penultimate column of Figure \ref{UnitSquareP10to30MeshRefinement} (Figures \ref{fig:UnitSquareP10U10TH}-\ref{fig:UnitSquareP30U11TH}). As stated in \cite{KawohlFridman:2003}, the first eigenvalue of $p$-Laplacian converges to the Cheeger constant of $\Omega$ as $p\to 1^+$ and the characteristic function of the Cheeger domain is the associated eigenfunction of $1$-Laplacian. In the current case $\Omega=(0,1)^2$, the Cheeger domain is the unit square with each of its four corners rounded off by circular arcs of radius $1/(2+\sqrt{\pi})$ \cite[Theorem 3.9 and Remark 3.11]{Leonardi:2015}. The computed first eigenfunction in Figure \ref{fig:UnitSquareP1_2U9TH} is obviously an approximation of the characteristic function of the relevant Cheeger domain.
On the other hand, further inspection of Figures \ref{fig:UnitSquareP30Mesh7TH}-\ref{fig:UnitSquareP30Mesh10TH} for $p=30$ reveals that refinements are largely performed towards singularities in the vicinity of the center. The observation indirectly confirms the assertion in a recent paper \cite{BrustadLindgrenLindqvist:2023} that the $\infty$-ground state, as the limit of the first eigenfunction of $p$-Laplacian when $p \to \infty$, is $\infty$-harmonic in the viscosity sense and further continuously differentiable \cite{Savin:2005} in the unit square except on two diagonal segments lying in a symmetric neighbourhood around the center. 

\end{example}

\begin{table}[hbt!]
\centering
\caption{Quantitative result for $p\in\{1.2, 1.5, 2, 2.5, 3\}$ in Example \ref{example2}: the number of adaptive loops, the number of vertices and the computed first eigenvalue.}
\label{tab:UnitSquareEigenvalue}
\resizebox{\textwidth}{!}{
\begin{tabular}{c cc cc cc cc cc}
    \toprule
    \multirow{2}{*}{$k$}&
    \multicolumn{2}{c}{$p=1.2$}&
    \multicolumn{2}{c}{$p=1.5$}&\multicolumn{2}{c}{$p=2$}&
    \multicolumn{2}{c}{$p=2.5$}&\multicolumn{2}{c}{$p=3$}\\
    \cmidrule(lr){2-11} 
     ~ & vertices & $\mu_k$ & vertices & $\mu_k$ & vertices & $\mu_k$ & vertices & $\mu_k$ & vertices & $\mu_k$ \\
    \midrule
        0 & 365   & 6.25529 & 365   & 10.1415 & 365   & 19.8951 & 365   & 36.3181  & 365   & 63.6013  \\

        1 & 590   & 6.23466 & 622   & 10.1186 & 637   & 19.8497 & 632   & 36.2091  & 629   & 63.3302  \\

        2 & 975   & 6.22154 & 1022  & 10.1038 & 1083  & 19.8148 & 1098  & 36.1092  & 1056  & 63.1296  \\

        3 & 1605  & 6.21414 & 1744  & 10.0929 & 1885  & 19.7843 & 1908  & 36.0469  & 1842  & 62.9704  \\

        4 & 2589  & 6.20896 & 2995  & 10.0840 & 3257  & 19.7670 & 3290  & 36.0087  & 3150  & 62.8610  \\

        5 & 4098  & 6.20632 & 5071  & 10.0791 & 5640  & 19.7552 & 5719  & 35.9825  & 5416  & 62.8523  \\

        6 & 6238  & 6.20505 & 8775  & 10.0764 & 9908  & 19.7484 & 9928 & 35.9697  & 9289  & 62.7673  \\

        7 & 9569 & 6.20475 & 15097 & 10.0745 & 17197 & 19.7445 & 16952 & 35.9584  & 15780 & 62.7476  \\

        8 &      &        & 26287 & 10.0735 & 29834 & 19.7424 & 28760 & 35.9533 &        &      \\

        9 &       &        &       &         & 51244 & 19.7410 & 48867 & 35.9526 &       &         \\
    \midrule
    reference & 23972 & 6.19606 & 61431 & 10.0723 & 61431 & 19.7400 & 61431 & 35.9473 & 61431 &  62.7522  \\
    \bottomrule
    \end{tabular}}
\end{table}

\begin{table}[hbt!]
\centering
\caption{Quantitative result for $p \in \{4, 10, 20, 30\}$ in Example \ref{example2}: the number of adaptive loops, the number of vertices and the computed first eigenvalue.}
\label{tab:UnitSquareEigenvalue102030}
\resizebox{\textwidth}{!}{
\begin{tabular}{c cc cc cc cc}
    \toprule
    \multirow{2}{*}{$k$}&\multicolumn{2}{c}{$p=4$}&
    \multicolumn{2}{c}{$p=10$}&
    \multicolumn{2}{c}{$p=20$}&\multicolumn{2}{c}{$p=30$}\\
    \cmidrule(lr){2-9}
     ~ & vertices & $\mu_k$ & vertices & $\mu_k$ & vertices & $\mu_k$ & vertices & $\mu_k$ \\
    \midrule
        0 & 365  & 180.480 & 365   & $4.01445\times 10^{4}$ & 365   & $1.31523\times 10^{8}$ & 365   & $3.32338\times 10^{11}$ \\

        1 & 592  & 179.080 & 478   & $3.84212\times 10^{4}$ & 459   & $1.14423\times 10^{8}$ & 538   & $2.63117\times 10^{11}$ \\

        2 & 1002  & 178.081 & 692   & $3.72617\times 10^{4}$ & 613   & $1.04613\times 10^{8}$ & 684  & $2.29447\times 10^{11}$ \\

        3 & 1679 & 177.536 & 1055  & $3.68261\times 10^{4}$ & 870  & $1.01255\times 10^{8}$ & 905  & $2.04656\times 10^{11}$ \\

        4 & 2849 & 177.190 & 1654  & $3.64416\times 10^{4}$ &  1308  & $9.82599\times 10^{7}$ & 1273  & $1.93808\times 10^{11}$ \\

        5 & 4850 & 176.864 & 2693  & $3.62735\times 10^{4}$ & 2006  & $9.67871\times 10^{7}$ & 1810  & $1.84294\times 10^{11}$ \\

        6 & 8108 & 176.849 & 4359  & $3.61551 \times 10^{4}$ & 3089  & $9.57319\times 10^{7}$ & 2617  & $1.81342\times 10^{11}$ \\

        7 &      &          & 7134  & $3.60413\times 10^{4}$ & 4945  & $9.53923\times 10^{7}$ & 3888 & $1.78707\times 10^{11}$ \\

        8 &      &          & 11722 & $3.60351\times 10^{4}$ & 7902 & $9.49911\times 10^{7}$ & 5887 & $1.77028\times 10^{11}$ \\

        9 &      &          & 19286 & $3.59512\times 10^{4}$ & 12756 & $9.46224\times 10^{7}$ & 9026 & $1.75533\times 10^{11}$ \\
        10&      &          &       &            &       &                & 14000 & $1.75363\times 10^{11}$\\
        
    \midrule
    reference & 61431 & 176.750 & 69177  & $3.58589\times 10^{4}$ & 69177 & $9.25536\times 10^{7}$ & 69177 & $1.67213\times 10^{11}$ \\
    \bottomrule
    \end{tabular}}
\end{table}

\begin{figure}[htb!]
    \centering
    \subcaptionbox{$\mathcal{T}_4$ (5246) \label{fig:LshapedP1_1Mesh5TH}} {\includegraphics[scale=0.295]{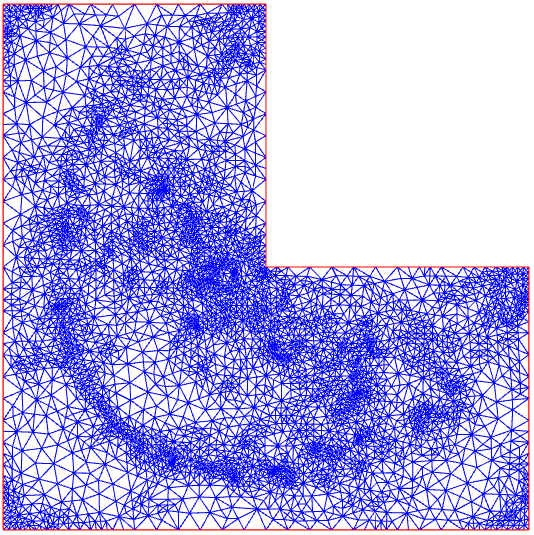}}\hspace{0.5cm}
    \subcaptionbox{$\mathcal{T}_6$ (13325) \label{fig:LshapedP_1Mesh7TH}}
    {\includegraphics[scale=0.295]{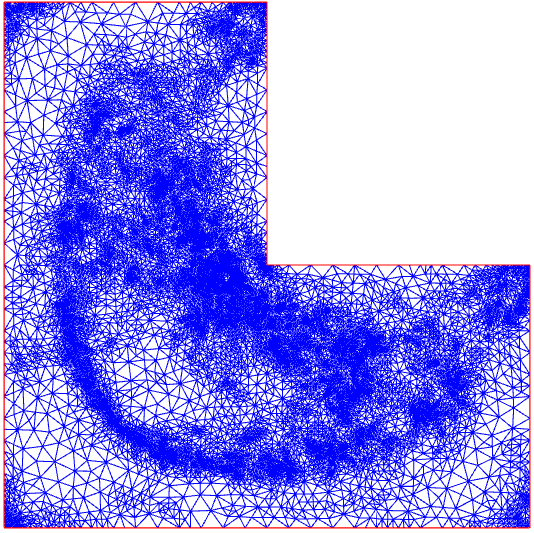}}\hspace{0.5cm}
    \subcaptionbox{$u_6$ \label{fig:LshapedP1_1U7TH}}{\includegraphics[scale=0.295]{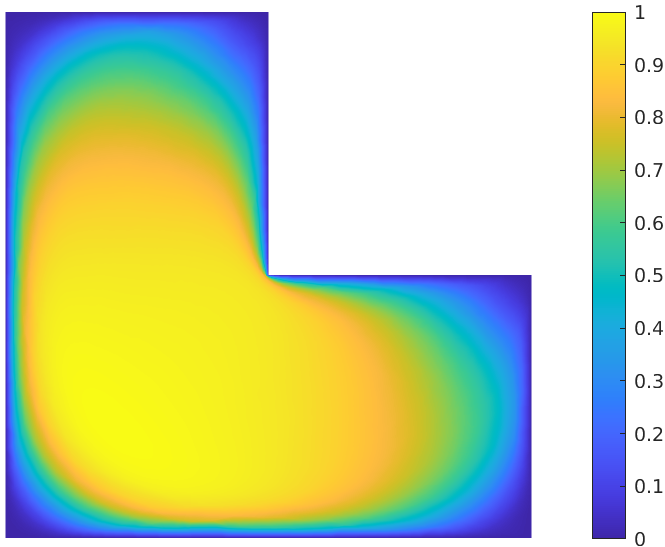}}\hspace{0.5cm}
    \subcaptionbox{reference \label{fig:LshapedP1_1U1THUniform}}{\includegraphics[scale=0.295]{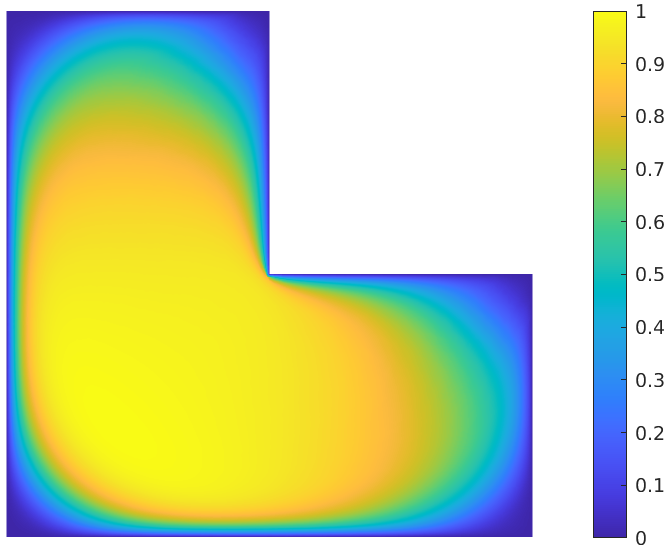}}
\\
    \subcaptionbox{$\mathcal{T}_5$ (10657) \label{fig:LshapedP1_5Mesh6TH}} {\includegraphics[scale=0.3]{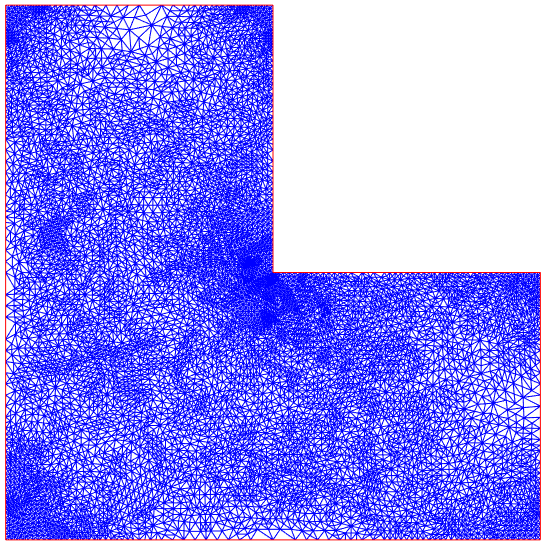}}\hspace{0.5cm}
    \subcaptionbox{$\mathcal{T}_7$ (31828) \label{fig:LshapedP1_5Mesh8TH}} {\includegraphics[scale=0.3]{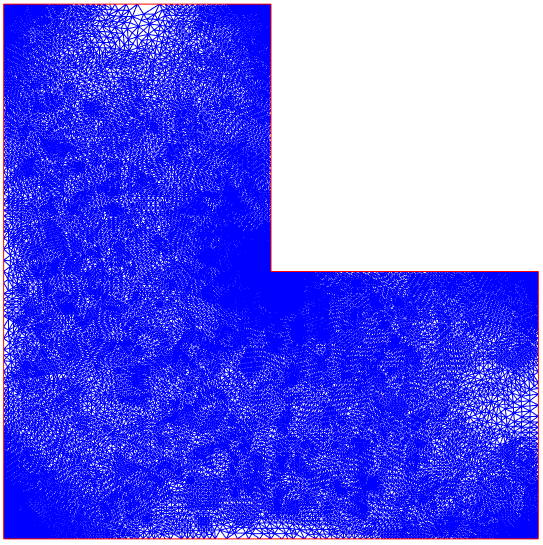}}\hspace{0.5cm}
    \subcaptionbox{$u_9$ \label{fig:LshapedP1_5U10TH}} {\includegraphics[scale=0.3]{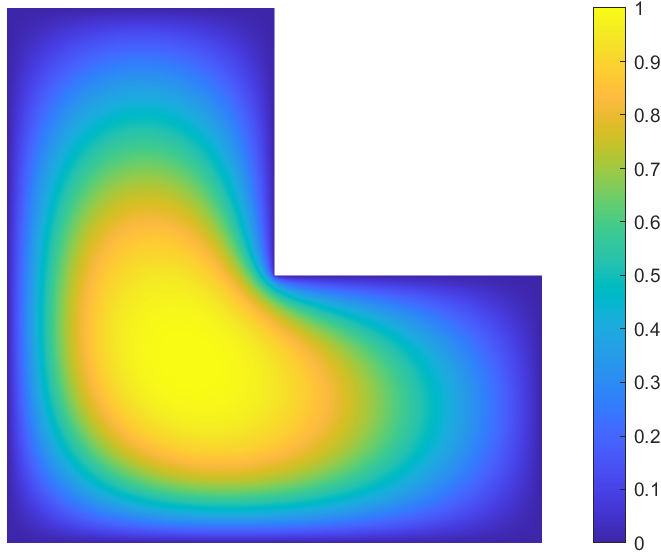}}\hspace{0.5cm}
	\subcaptionbox{reference \label{fig:LshapedP1_5U1THUniform}}{\includegraphics[scale=0.3]{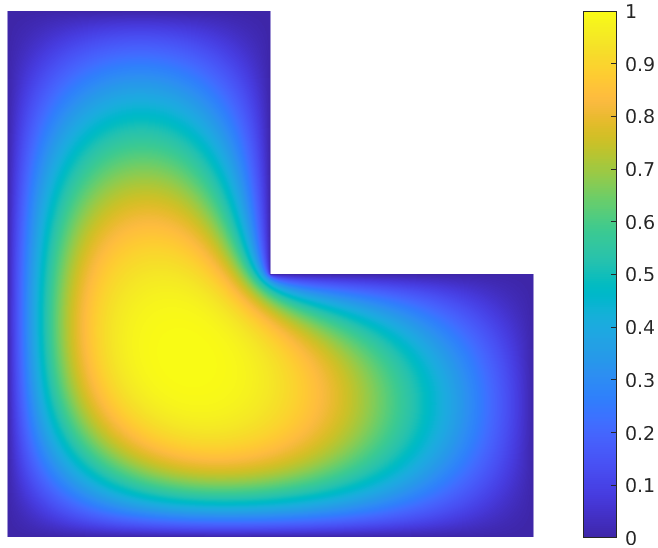}}
    \caption{Adaptive mesh refinement level with the number of vertices over each mesh as well as final computed first eigenfunctions by adaptive refinements and references for $p=1.1$ (1st row) and $p=1.5$ (2nd row) in Example \ref{example3}.}
    \label{LshapedP1to2MeshRefinement}
\end{figure}

\begin{figure}[htb!]
    \centering
    \subcaptionbox{$\mathcal{T}_3$ (2903) \label{fig:LshapedP2_5Mesh4TH}} {\includegraphics[scale=0.3]{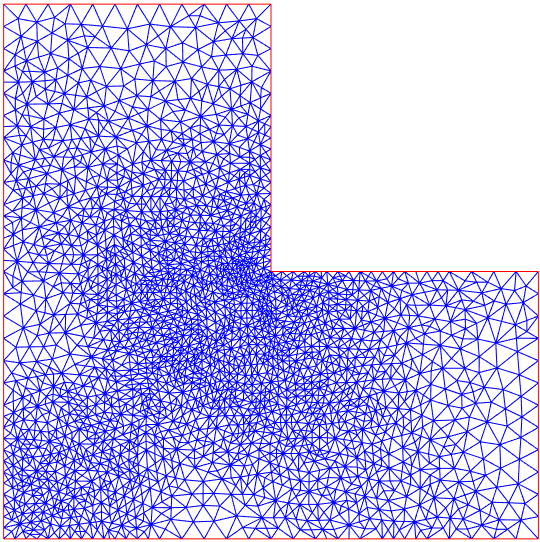}}\hspace{0.5cm}
    \subcaptionbox{$\mathcal{T}_5$ (8119) \label{fig:LshapedP2_5Mesh6TH}} {\includegraphics[scale=0.3]{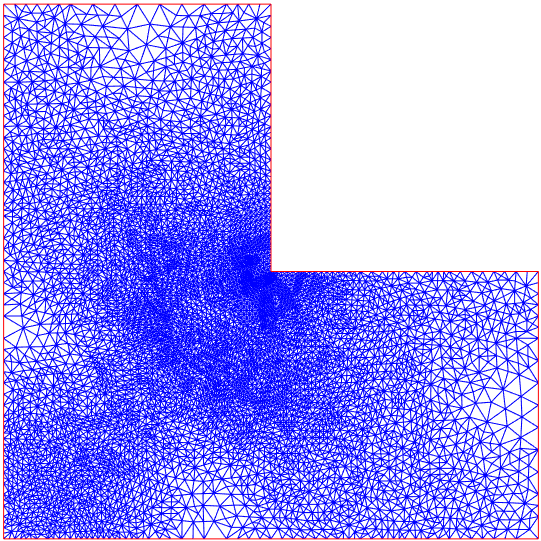}}\hspace{0.5cm}
    \subcaptionbox{$\mathcal{T}_8$ (38111) \label{fig:LshapedP2_5Mesh9TH}} {\includegraphics[scale=0.3]{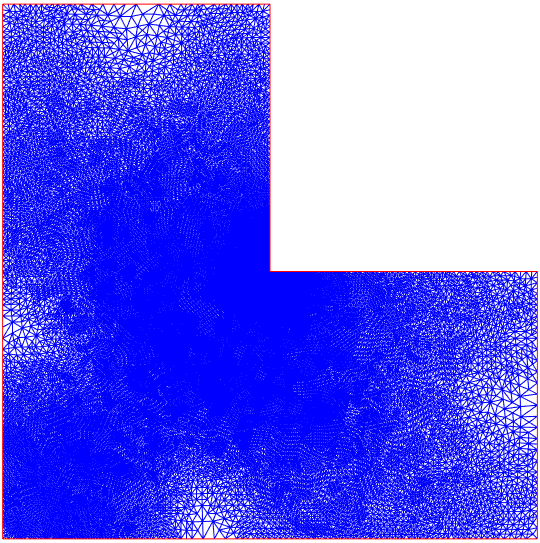}}\hspace{0.5cm}
    \subcaptionbox{$u_{11}$ \label{fig:LshapedP2_5U12TH}} {\includegraphics[scale=0.3]{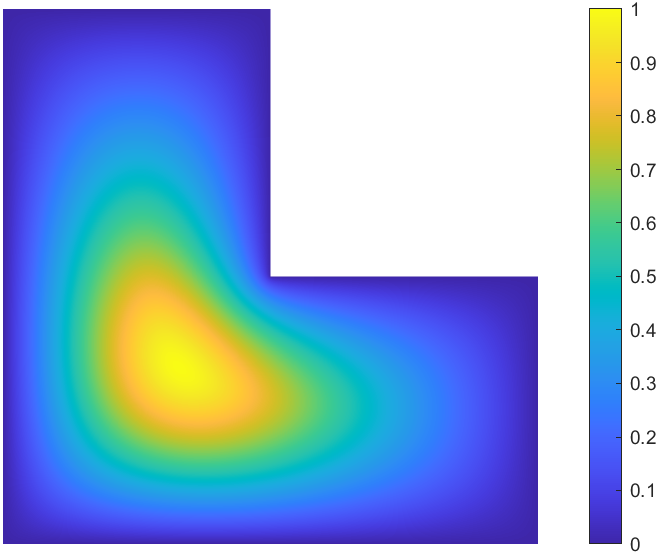}}
\\
    \subcaptionbox{$\mathcal{T}_3$ (2087) \label{fig:LshapedP4Mesh4TH}} {\includegraphics[scale=0.295]{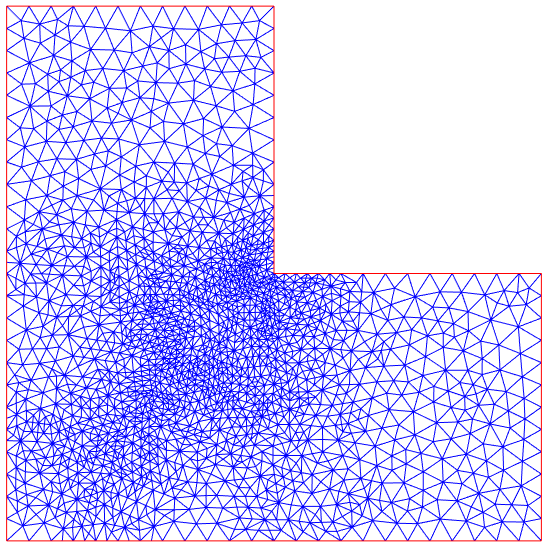}}\hspace{0.5cm}
    \subcaptionbox{$\mathcal{T}_5$ (5144) \label{fig:LshapedP4Mesh6TH}} {\includegraphics[scale=0.295]{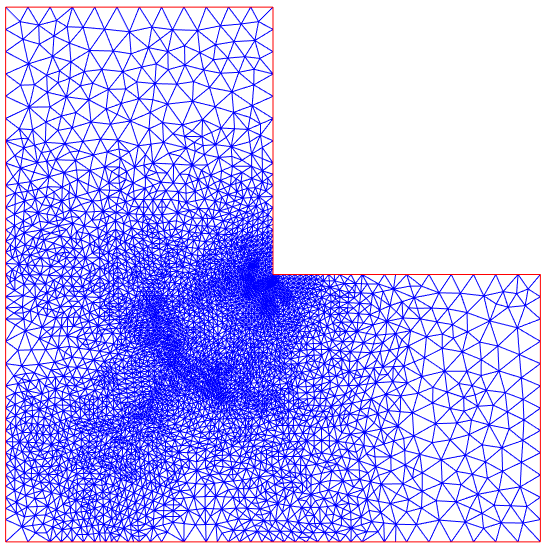}}\hspace{0.5cm}
    \subcaptionbox{$\mathcal{T}_9$ (33364) \label{fig:LshapedP4Mesh10TH}}
    {\includegraphics[scale=0.295]{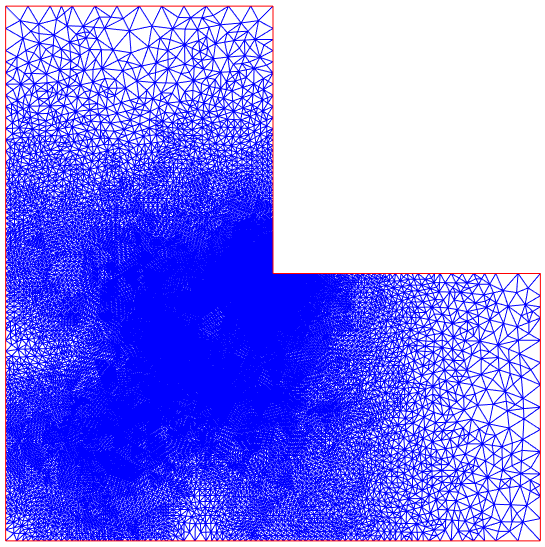}}\hspace{0.5cm}
    \subcaptionbox{$u_{12}$ \label{fig:LshapedP4U13TH}}
    {\includegraphics[scale=0.295]{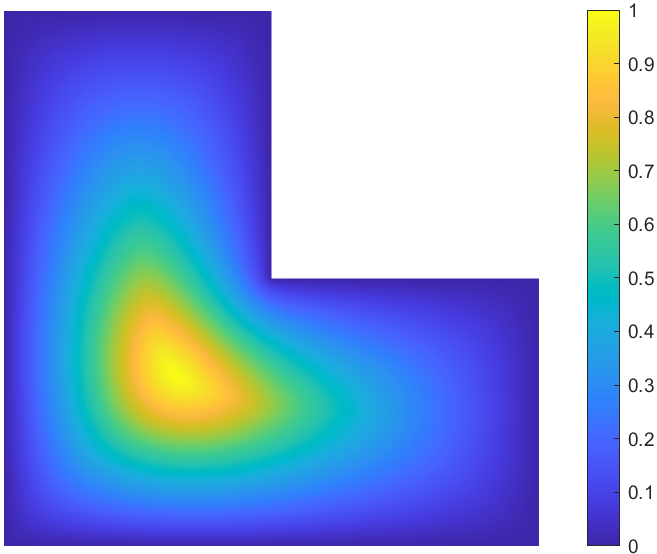}}
    \caption{Adaptive mesh refinement level with the number of vertices over each mesh as well as final computed first eigenfunctions by adaptive refinements for $p=2.5$ (1st row) and $4$ (2nd row) in Example \ref{example3}.}
    \label{LshapedP2to4MeshRefinement}
\end{figure}

\begin{figure}[htb!]
    \centering
    \subcaptionbox{$\mathcal{T}_2$ (1016) \label{fig:LshapedP10Mesh3TH}} {\includegraphics[scale=0.3]{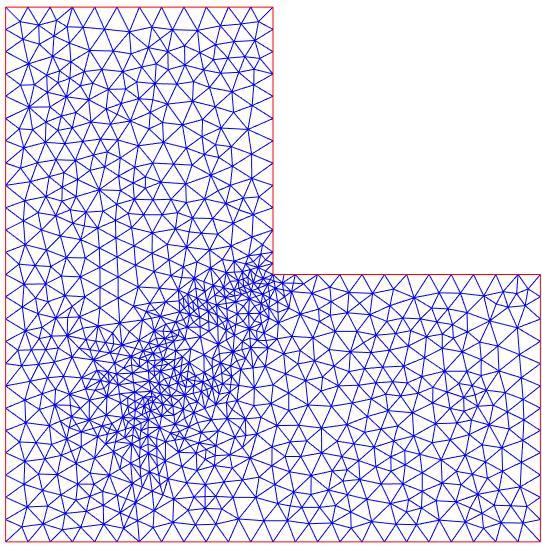}}\hspace{0.5cm}
    \subcaptionbox{$\mathcal{T}_9$ (11912) \label{fig:LshapedP10Mesh10TH}}{\includegraphics[scale=0.3]{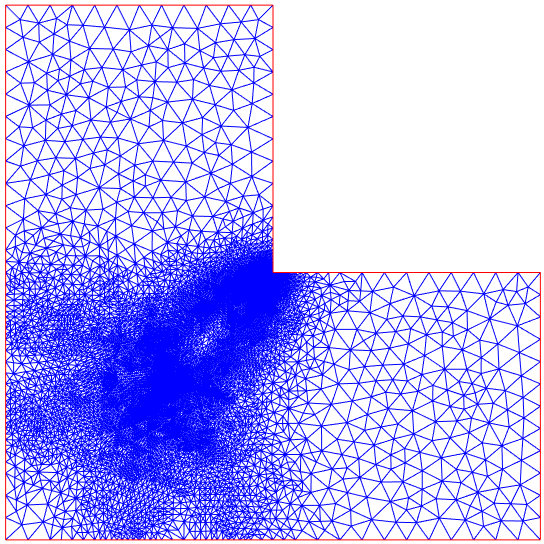}}\hspace{0.5cm}
    \subcaptionbox{$u_9$ \label{fig:LshapedP10U10TH}}{\includegraphics[scale=0.3]{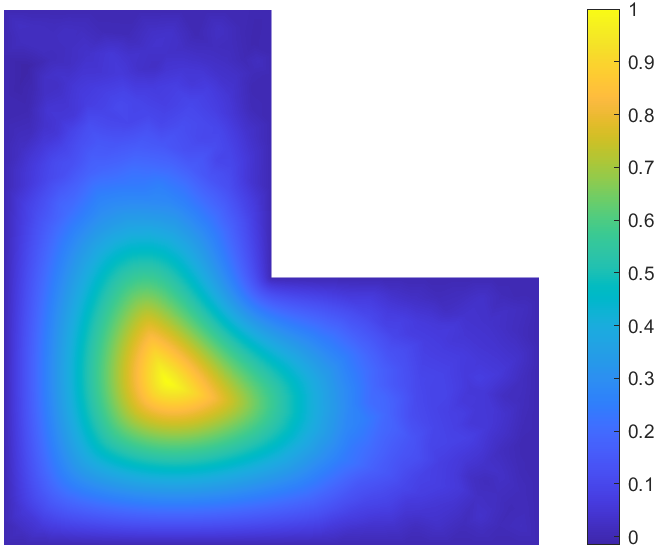}}\hspace{0.5cm}
	\subcaptionbox{reference \label{fig:LshapedP10U1THUniform}}{\includegraphics[scale=0.3]{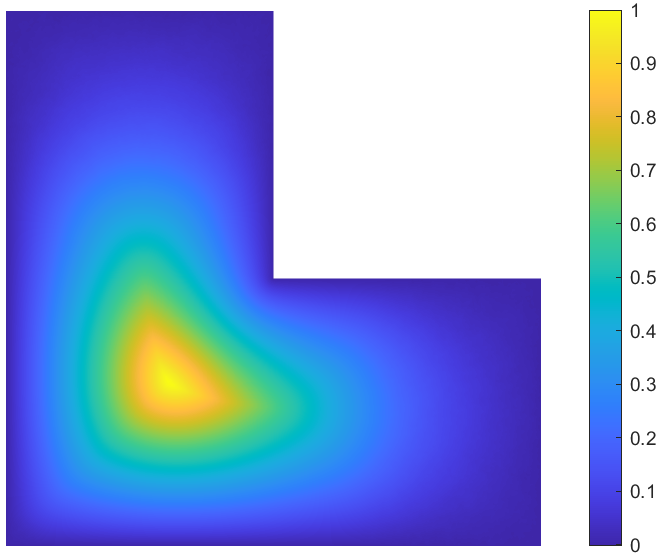}}
\\
    \subcaptionbox{$\mathcal{T}_3$ (1053) \label{fig:LshapedP30Mesh4TH}} {\includegraphics[scale=0.3]{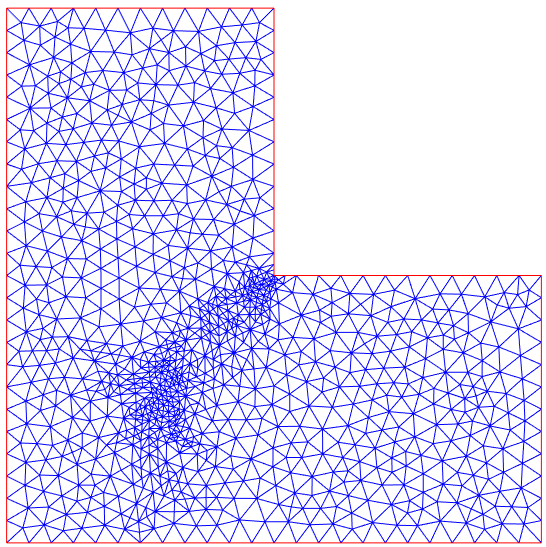}}\hspace{0.5cm}
    \subcaptionbox{$\mathcal{T}_{10}$ (9148) \label{fig:LshapedP30Mesh11TH}}{\includegraphics[scale=0.3]{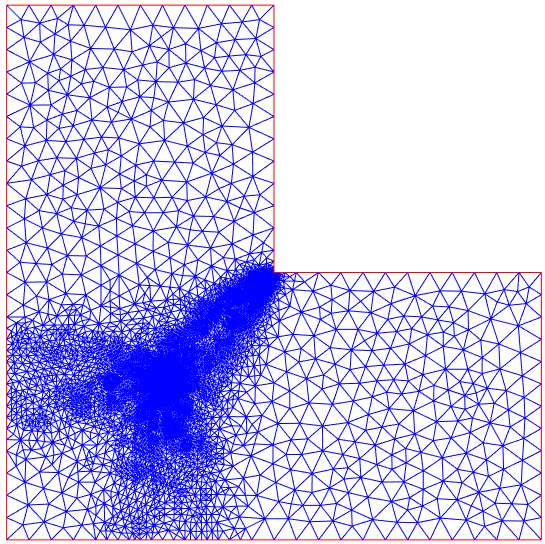}}\hspace{0.5cm}
    \subcaptionbox{$u_{11}$ \label{fig:LshapedP30U12TH}}{\includegraphics[scale=0.3]{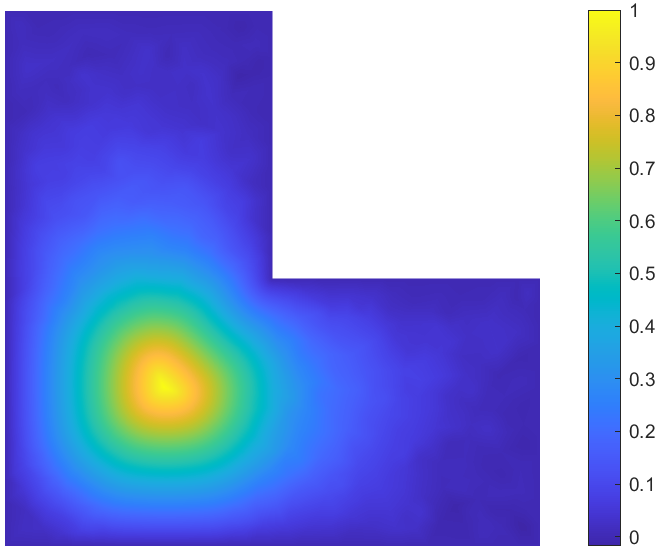}}\hspace{0.5cm}
	\subcaptionbox{reference \label{fig:LshapedP30U1THUniform}}{\includegraphics[scale=0.3]{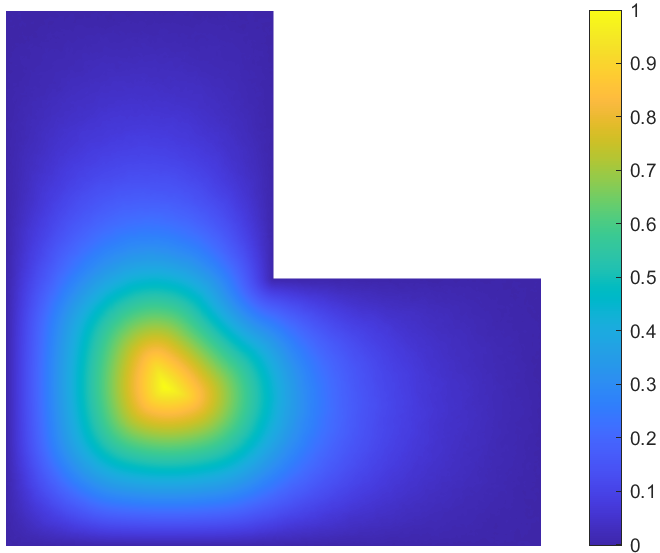}}
    \caption{Adaptive mesh refinement level with the number of vertices over each mesh as well as final computed first eigenfunctions by adaptive refinements and references for $p=10$ (1st row) and $30$ (2nd row) in Example \ref{example3}.}
    \label{LshapedP10to30MeshRefinement}
\end{figure}

\begin{example}[L-shaped Domain]\label{example3}
The third example is posed in an L-shaped domain $(0,2)^2\setminus [1,2)^2$ with the same values of $p$ as in Example \ref{example1}. Tolerances $\epsilon_K$ and $\epsilon_M$ are given by $10^{-4}$ and $10^{-5}$ for $p\in \{1.5, 2, 2.5, 3\}$, $5 \times 10^{-4}$ and $5 \times 10^{-5}$ for $p\in \{4, 10, 20, 30\}$ and $10^{-3}$ and $5 \times 10^{-4}$ for $p = 1.1$ respectively. A maximum iteration number is specified as $K=9$ for $p = 10, 20$ and $K = 11$ for $p = 30$ respectively. 

The convergence history of Algorithm \ref{afem} is presented in Tables \ref{tab:LshapedEigenvalue} and \ref{tab:LshapedEigenvalue102030}, which show similar convergence behaviour of computed first eigenvalues for each $p$. Compared with the reference solution over the fine mesh, fewer degrees of freedom by Algorithm \ref{afem} is required for more accuracy when $p$ is small. In particular, we observe from Table  \ref{tab:LshapedEigenvalue102030} with $p=4$ that the computed eigenvalues $\mu_{11}=54.6221$ and $\mu_{12} = 54.6211$ by Algorithm \ref{afem}, albeit with only $49\%$ and $78\%$ of degrees of freedom respectively, are both smaller than the reference solution over a very fine mesh.

Figures \ref{LshapedP1to2MeshRefinement}-\ref{LshapedP10to30MeshRefinement} show that the marked set from Algorithm \ref{afem} transfers from the region adjacent to the six corners to the vicinity of $(\frac{1}{2},\frac{1}{2})$ and the area around the reentrant corner as $p$ increases, and the singularities are clearly detected for each $p$. From the penultimate column of Figure \ref{LshapedP1to2MeshRefinement} (Figures \ref{fig:LshapedP1_1U7TH}-\ref{fig:LshapedP1_5U10TH}), we may deduce as in Example \ref{example2} that the computed first eigenfunctions are approximations of the characteristic function associated with some Cheeger domain \cite{KawohlFridman:2003} as the first eigenfunction of 1-Laplacian in $\Omega$. In addition, it seems that singularities for large $p$ might exist along the segment connecting $(\frac{1}{2},\frac{1}{2})$ to $(1,1)$ as displayed in Figure \ref{fig:LshapedP30Mesh11TH}. To the best of our knowledge, no reasonable explanation of this observation is available in the PDE theory. This might provide a clue about the regularity of the $\infty$-eigenvalue problem in a non-convex domain.
\end{example}

\begin{table}[hbt!]
\centering
\caption{Quantitative result for $p\in\{1.1, 1.5, 2, 2.5, 3\}$ in Example \ref{example3}: the number of adaptive loops, the number of vertices and the computed first eigenvalue.}
\label{tab:LshapedEigenvalue}
\resizebox{\textwidth}{!}{
\begin{tabular}{c cc cc cc cc cc}
    \toprule
    \multirow{2}{*}{$k$}&
    \multicolumn{2}{c}{$p=1.1$}&
    \multicolumn{2}{c}{$p=1.5$}&\multicolumn{2}{c}{$p=2$}&
    \multicolumn{2}{c}{$p=2.5$}&
    \multicolumn{2}{c}{$p=3$}\\
    \cmidrule(lr){2-11}
     ~ & vertices & $\mu_k$ & vertices & $\mu_k$ & vertices & $\mu_k$ & vertices & $\mu_k$ & vertices & $\mu_k$ \\
    \midrule
        0 & 741 & 3.29084 & 741    & 5.74153  & 741  & 9.75107 & 741   & 15.6646      & 741   & 24.4197 \\

        1 & 1240 & 3.27061 & 1245  & 5.71661  & 1188  & 9.70099 & 1111  & 15.5862  & 1064  & 24.2656 \\

        2 & 1994 & 3.26325 & 2121  & 5.70591  & 1990  & 9.68609 & 1766  & 15.5381  & 1606  & 24.1712 \\

        3 & 3286 & 3.25550 & 3643  & 5.69624  & 3385  & 9.66841 & 2903  & 15.5029  & 2556  & 24.1022\\

        4 & 5246 & 3.25223 & 6245  & 5.69123  & 5748  & 9.65753 & 4856  & 15.4807  & 4193  & 24.0693 \\

        5 & 8394 & 3.24729 & 10657  & 5.68770  & 9842  & 9.65065 & 8119  & 15.4662 & 6889  & 24.0321  \\

        6 &13325 & 3.24471 & 18444 & 5.68576 & 16767  & 9.64646 & 13653 & 15.4577  & 11378 & 24.0204 \\

        7 &      &        & 31828 & 5.68436  & 28561 & 9.64385 & 22892 & 15.4511  & 18715 & 24.0061 \\

        8 &      &        & 54608 & 5.68362  & 48554 & 9.64225 & 38111 & 15.4488  & 30495 & 24.0023 \\

        9 &       &        & 94803 & 5.68311 & 82874  & 9.64127 & 63517 & 15.4469  & 50007 & 23.9932\\

        10 &      &        &       &        & 140675 & 9.64066   & 105876 & 15.4448 & 81801 & 23.9923 \\
        11 &      &        &       &        &        &          & 175649 & 15.4444     \\
    \midrule
    reference & 20390 & 3.24002 & 122801 & 5.683402  & 209247 & 9.64097  & 256870 & 15.4449 & 256870 & 23.9903 \\
    \bottomrule
    \end{tabular}}
\end{table}

\begin{table}[hbt!]
\centering
\caption{Quantitative result for $p\in\{4, 10, 20, 30\}$ in Example \ref{example3}: the number of adaptive loops, the number of vertices and the computed first eigenvalue.}
\label{tab:LshapedEigenvalue102030}
\resizebox{\textwidth}{!}{
\begin{tabular}{c cc cc cc cc}
    \toprule
    \multirow{2}{*}{$k$}&
    \multicolumn{2}{c}{$p=4$}&\multicolumn{2}{c}{$p=10$}&
    \multicolumn{2}{c}{$p=20$}&\multicolumn{2}{c}{$p=30$}\\
    \cmidrule(lr){2-9}
     ~ & vertices & $\mu_k$ & vertices & $\mu_k$ & vertices & $\mu_k$ & vertices & $\mu_k$ \\
    \midrule
        0 & 741  & 56.1189 & 741  & $4.31754\times 10^{3}$  & 741 & $2.70810\times 10^{6}$ & 741  & $1.30341\times 10^{9}$ \\

        1 & 972 & 55.5998 & 848 & $4.15542\times 10^{3}$    & 820 & $2.46503\times 10^{6}$ & 813 & $1.11285\times 10^{9}$  \\

        2 & 1384 & 55.2551 & 1016 & $4.01369\times 10^{3}$  & 943  & $2.32909\times 10^{6}$ & 907  & $9.12838\times 10^{8}$ \\

        3 & 2087 & 55.0287 & 1296 & $3.93108\times 10^{3}$  & 1144  & $2.12411\times 10^{6}$ & 1053  & $8.74190\times 10^{8}$ \\

        4 & 3267 & 54.8835 & 1762 & $3.88738\times 10^{3}$  & 1418  & $2.04914\times 10^{6}$ & 1267  & $7.94392\times 10^{8}$ \\

        5 & 5144 & 54.7934 & 2483  & $3.85769\times 10^{3}$ & 1839  & $1.99152\times 10^{6}$ & 1616  & $7.60674\times 10^{8}$ \\

        6 & 8215 & 54.7351 & 3553 & $3.84105\times 10^{3}$  & 2477 & $1.97362\times 10^{6}$ & 2153  & $7.43803\times 10^{8}$ \\

        7 & 13142 & 54.6975 & 5245 & $3.82974\times 10^{3}$ & 3491 & $1.95886\times 10^{6}$ & 2952 & $7.33018\times 10^{8}$ \\

        8 & 20845 & 54.6695 & 7862 & $3.82284\times 10^{3}$ & 5047 & $1.94891\times 10^{6}$ & 4229 & $7.24466\times 10^{8}$ \\

        9 & 33364 & 54.6506 & 11912 & $3.81893\times 10^{3}$& 7444 & $1.94260\times 10^{6}$ & 6180 & $7.22056\times 10^{8}$ \\

        10 & 53165 & 54.6445 &     &                         &     &                        & 9148  & $7.19039\times 10^{8}$ \\

        11 & 84384   & 54.6221        &   &   &  &                                          & 13697  & $7.17379\times 10^{8}$ \\
        12 & 134431   & 54.6211       &   &   &  &   &  & \\
    \midrule
    reference & 172212 & 54.6323 & 44439 & $3.81345\times 10^{3}$ & 44439 & $1.92416\times 10^{6}$ & 44439 & $7.03153\times 10^{8}$ \\
    \bottomrule
    \end{tabular}}
\end{table}


\begin{figure}[htb!]
    \centering
    \subcaptionbox{$\mathcal{T}_2$ (3789) \label{fig:CubeP1_5Mesh3}} {\includegraphics[scale=0.3]{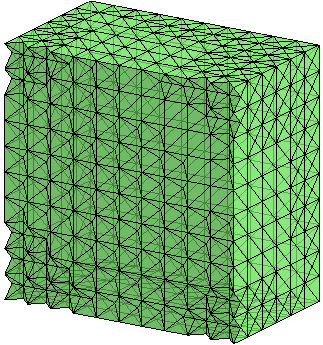}}\hspace{0.5cm}
    \subcaptionbox{$\mathcal{T}_5$ (23069) \label{fig:CubeP1_5Mesh6}} {\includegraphics[scale=0.3]{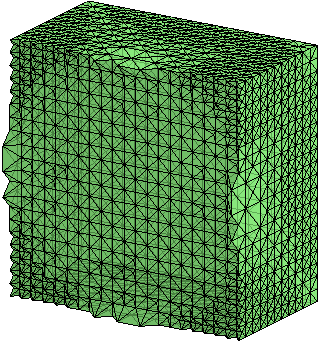}}\hspace{0.5cm}
    \subcaptionbox{$\mathcal{T}_9$ (284332) \label{fig:CubeP1_5Mesh10}} {\includegraphics[scale=0.3]{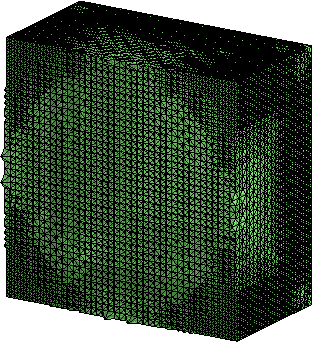}}\hspace{0.5cm}
    \subcaptionbox{$u_9$ \label{fig:CubeP1_5U10}} {\includegraphics[scale=0.3]{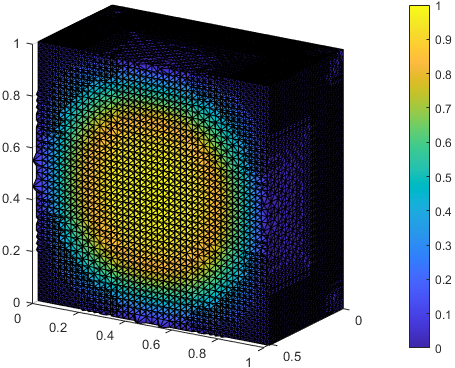}}
\\
    \subcaptionbox{$\mathcal{T}_2$ (4527) \label{fig:CubeP2Mesh3}} {\includegraphics[scale=0.3]{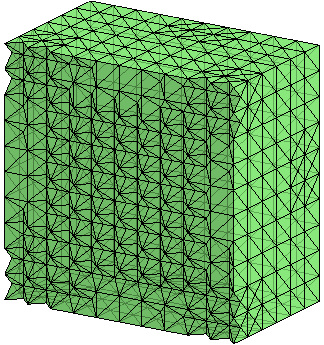}}\hspace{0.5cm}
    \subcaptionbox{$\mathcal{T}_5$ (30499) \label{fig:CubeP2Mesh6}} {\includegraphics[scale=0.3]{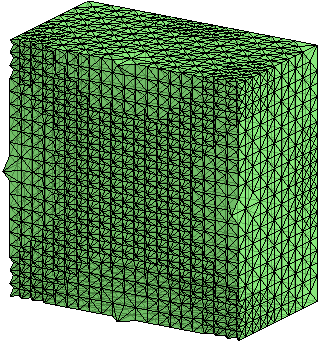}}\hspace{0.5cm}
    \subcaptionbox{$\mathcal{T}_9$ (362383) \label{fig:CubeP2Mesh10}} {\includegraphics[scale=0.3]{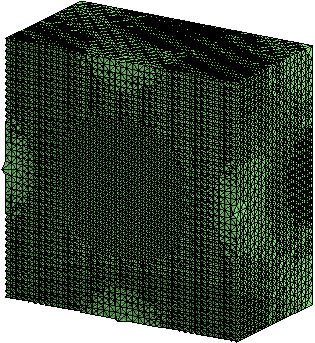}}\hspace{0.5cm}
    \subcaptionbox{$u_9$ \label{fig:CubeP2U10}} {\includegraphics[scale=0.3]{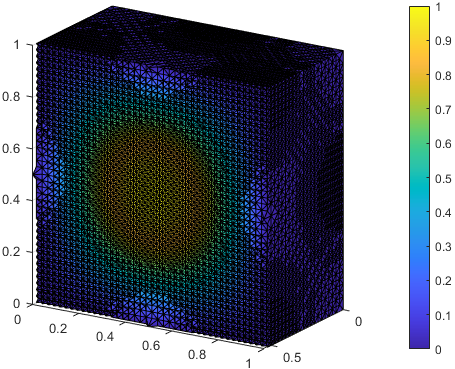}}
\\
    \subcaptionbox{$\mathcal{T}_3$ (2157) \label{fig:CubeP4Mesh4}} {\includegraphics[scale=0.3]{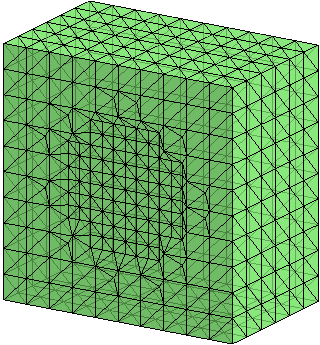}}\hspace{0.5cm}
    \subcaptionbox{$\mathcal{T}_6$ (8223) \label{fig:CubeP4Mesh7}} {\includegraphics[scale=0.3]{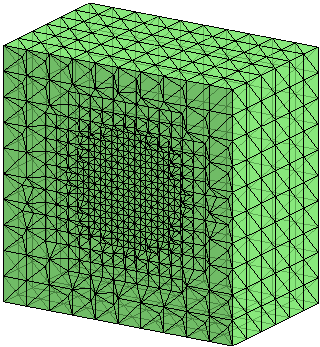}}\hspace{0.5cm}
    \subcaptionbox{$\mathcal{T}_{12}$ (221479) \label{fig:CubeP4Mesh13}} {\includegraphics[scale=0.3]{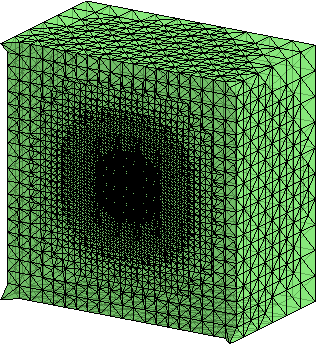}}\hspace{0.5cm}
    \subcaptionbox{$u_{12}$ \label{fig:CubeP4U13}} {\includegraphics[scale=0.3]{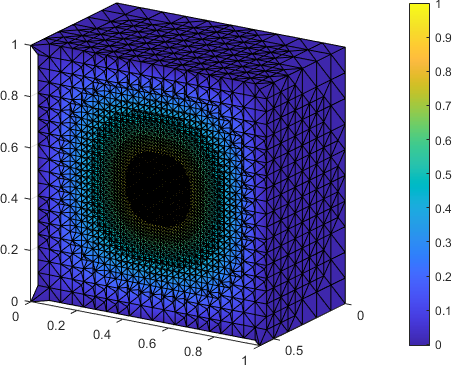}}
\\
    \subcaptionbox{$\mathcal{T}_4$ (1453) \label{fig:CubeP10Mesh5}} {\includegraphics[scale=0.3]{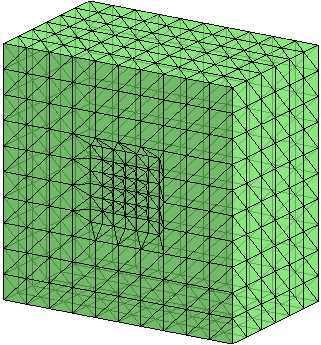}}\hspace{0.5cm}
    \subcaptionbox{$\mathcal{T}_9$ (6508) \label{fig:CubeP10Mesh10}} {\includegraphics[scale=0.3]{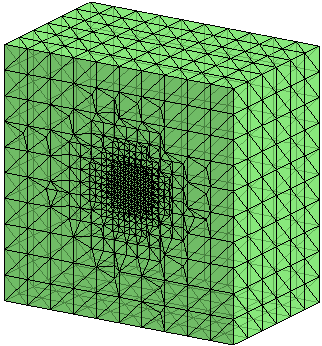}}\hspace{0.5cm}
    \subcaptionbox{$\mathcal{T}_{14}$ (92043) \label{fig:CubeP10Mesh15}} {\includegraphics[scale=0.3]{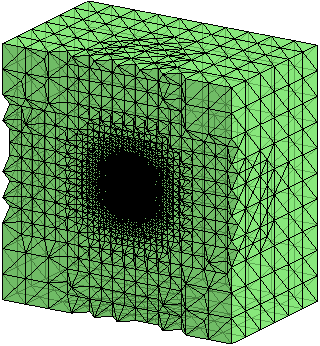}}\hspace{0.5cm}
    \subcaptionbox{$u_{14}$ \label{fig:CubeP10U15}} {\includegraphics[scale=0.3]{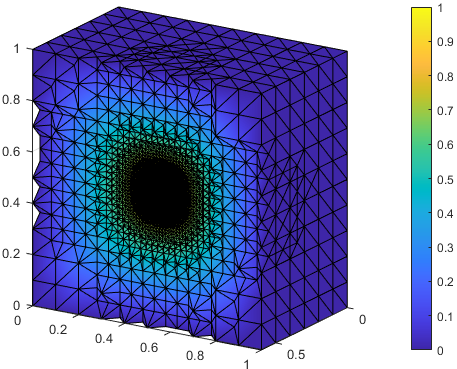}}
    \caption{Adaptive mesh refinement level with the number of vertices over each mesh as well as final computed first eigenfunctions by adaptive refinements in the cross section $x<0.5$ for $p=1.5, 2, 4, 10$ from top to bottom in Example \ref{example4}.}
    \label{CubeP1to10MeshRefinement}
\end{figure}

\begin{example}[Unit Cube]\label{example4}
The first 3-d experiment is conducted in the unit cube $(0,1)^3$ with 4 different values of $p\in\{1.5, 2, 4, 10\}$. Tolerance $\epsilon_{M}$ in Algorithm \ref{alg:alg1} is fixed at $10^{-5}$. Due to the limited computing resource, a maximum adaptive refinement number $K$ is prescribed instead of a tolerance $\epsilon_K$ for all $p$.

Table \ref{tab:CubeEigenvalue} presents the number of vertices and computed first eigenvalues over adaptively refined meshes. As shown in the previous 2-d numerical tests, the sequence of adaptive eigenvalues for each $p$ monotonically decreases. It is worth mentioning that $\mu_9 = 29.6297$ for $p=2$ attains an relative error about $7 \times 10^{-4}$ with respect to the exact first eigenvalue $3\pi^2 \approx 29.6088$ of Laplacian in the unit cube. To illustrate the efficiency of Algorithm \ref{afem}, the numerical simulation of $p=2$ yields an approximate first eigenvalue $29.9918$ over a fine mesh with $373248$ vertices whereas our adaptive algorithm can achieve higher accuracy ($\mu_8=29.6456$) at a cost of almost half vertices ($188348$). 

In Figure \ref{CubeP1to10MeshRefinement}, we depict the evolution of adaptive meshes and computed first eigenfunctions over the finest adaptive meshes in the cross section $x<0.5$ due to the symmetry. As in cases of the 2-d unit square, refinements are largely performed in the vicinity of the center and around four corners for small $p$ and then in the region where the computed first eigenfunctions are non-zero for large $p$. 

\end{example}

\begin{table}[hbt!]
\centering
\caption{Quantitative result for $p \in \{1.5, 2, 4, 10\}$ in Example \ref{example4}: the number of adaptive loops, the number of vertices and the computed first eigenvalue.}
\label{tab:CubeEigenvalue}
\begin{tabular}{c cc cc cc cc}
    \toprule
    \multirow{2}{*}{$k$}&\multicolumn{2}{c}{$p=1.5$}&
    \multicolumn{2}{c}{$p=2$}&
    \multicolumn{2}{c}{$p=4$}&\multicolumn{2}{c}{$p=10$}\\
    \cmidrule(lr){2-9}
     ~ & vertices & $\mu_k$ & vertices & $\mu_k$ & vertices & $\mu_k$ & vertices & $\mu_k$ \\
    \midrule
        0 & 1331 & 14.9719 & 1331  & 30.8578 & 1331  &  346.595 & 1331  & $1.48342\times 10^{5}$ \\

        1 & 2003 & 14.8508 & 2183  & 30.4664 & 1428  & 339.998 & 1341  & $1.32567 \times 10^{5}$ \\

        2 & 3789 & 14.6678 & 4527  & 29.9688 & 1713  & 333.102 & 1389  & $1.20605\times 10^{5}$ \\

        3 & 7193 & 14.6143 & 7798  & 29.9302 & 2157  & 329.550 & 1453  & $1.18186\times 10^{5}$ \\

        4 & 11871& 14.5777 & 14095 & 29.8301 & 3049  & 323.672 & 1605  &  $1.09896\times 10^{5}$ \\

        5 & 23069& 14.5361 & 30499 & 29.7193 & 5167  & 319.016 & 1810  & $1.06198\times 10^{5}$ \\

        6 & 42768& 14.5179 & 53979 & 29.6903 & 8223  &  315.660 & 2215  & $1.02281\times 10^{5}$ \\

        7 & 76375&  14.5053 & 92036 & 29.6710 & 13864  & 312.402 & 2899 &  $9.91547\times 10^{4}$ \\

        8 &148519&  14.4948 & 188348 & 29.6456 & 24899 & 310.919 & 4363 & $9.57893\times 10^{4}$ \\

        9 &284332&  14.4892 & 362383 & 29.6297 & 41528 & 310.086 & 6508 & $9.21262\times 10^{4}$ \\

        10 &      &          &       &         & 70878 &  309.285 & 10285 & $9.04754\times 10^{4}$ \\
	
        11 &      &          &       &         & 128276 & 308.471 & 18264 & $8.77086\times 10^{4}$ \\
        12 &      &          &       &         & 221479  & 307.792 & 30195 & $8.66154\times 10^{4}$\\
        13 &      &          &       &         &        &           & 51359 & $8.55563\times 10^{4}$ \\
        14 &      &          &       &         &        &           & 92043 & $8.51059\times 10^{4}$ \\
    \bottomrule
    \end{tabular}
\end{table}

\begin{figure}[htb!]
    \centering
    \subcaptionbox{$\mathcal{T}_2$ (3114) \label{fig:Lshape3DP1_5Mesh3}} {\includegraphics[scale=0.3]{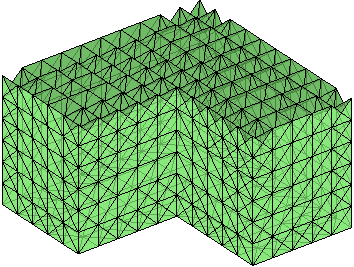}}\hspace{0.5cm}
    \subcaptionbox{$\mathcal{T}_5$ (19449) \label{fig:Lshape3DP1_5Mesh6}} {\includegraphics[scale=0.3]{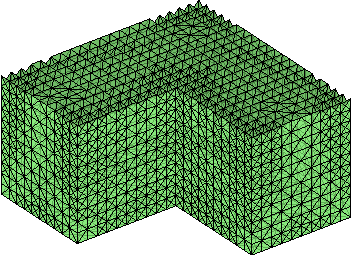}}\hspace{0.5cm}
    \subcaptionbox{$\mathcal{T}_8$ (123451) \label{fig:Lshape3DP1_5Mesh9}} {\includegraphics[scale=0.3]{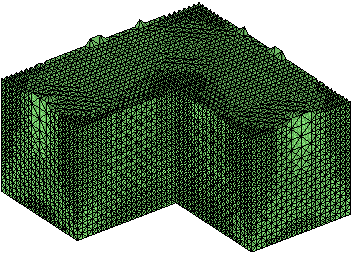}}\hspace{0.5cm}
    \subcaptionbox{$u_9$ \label{fig:Lshape3DP1_5U10}} {\includegraphics[scale=0.26]{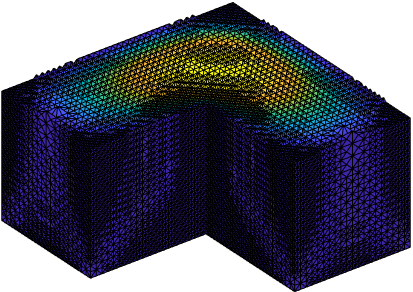}}
    \includegraphics[scale=0.25]{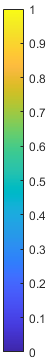}
\\
    \subcaptionbox{$\mathcal{T}_3$ (1740) \label{fig:Lshape3DP4Mesh4}} {\includegraphics[scale=0.3]{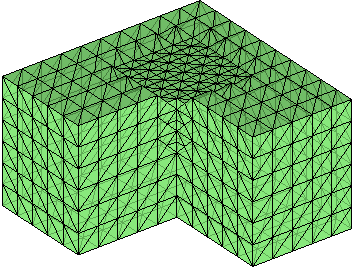}}\hspace{0.5cm}
    \subcaptionbox{$\mathcal{T}_6$ (6197) \label{fig:Lshape3DP4Mesh7}} {\includegraphics[scale=0.3]{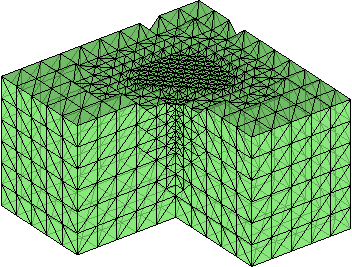}}\hspace{0.5cm}
    \subcaptionbox{$\mathcal{T}_{12}$ (162735) \label{fig:Lshape3DP4Mesh13}} {\includegraphics[scale=0.3]{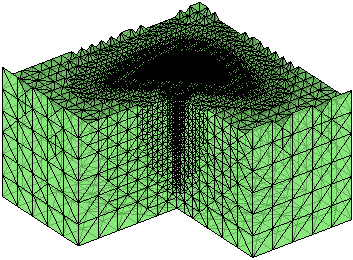}}\hspace{0.5cm}
    \subcaptionbox{$u_{12}$ \label{fig:Lshape3DP4U13}} {\includegraphics[scale=0.26]{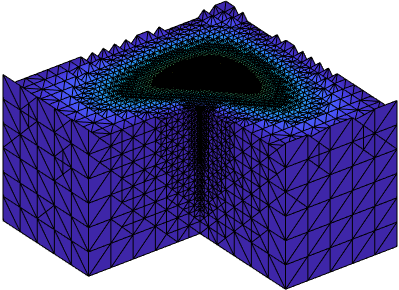}}
    \includegraphics[scale=0.25]{Lshape3DColorbar.png}
\\
    \subcaptionbox{$\mathcal{T}_4$ (1248) \label{fig:Lshape3DP10Mesh5}} {\includegraphics[scale=0.3]{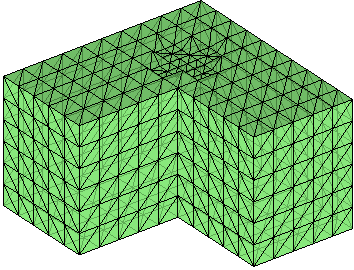}}\hspace{0.5cm}
    \subcaptionbox{$\mathcal{T}_9$ (5069) \label{fig:Lshape3DP10Mesh10}} {\includegraphics[scale=0.3]{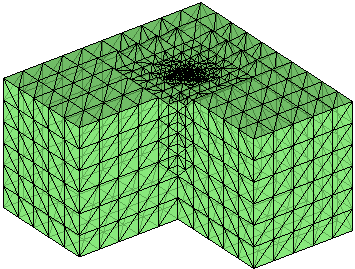}}\hspace{0.5cm}
    \subcaptionbox{$\mathcal{T}_{13}$ (38977) \label{fig:Lshape3DP10Mesh14}} {\includegraphics[scale=0.3]{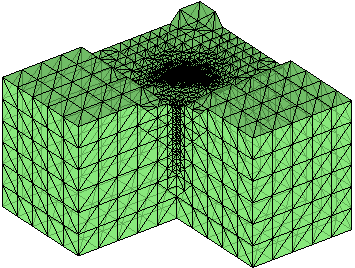}}\hspace{0.5cm}
    \subcaptionbox{$u_{13}$ \label{fig:Lshape3DP10U14}} {\includegraphics[scale=0.26]{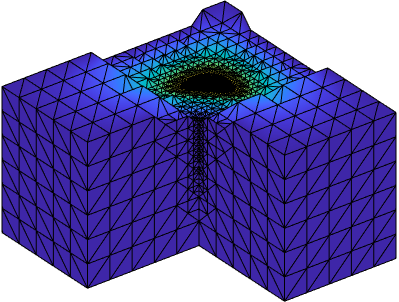}}
    \includegraphics[scale=0.25]{Lshape3DColorbar.png}
    \caption{Adaptive mesh refinement level with the number of vertices over each mesh as well as final computed first eigenfunctions by adaptive refinements in the cross section $z<1$ for Example \ref{example5}:  $p=1.5, 4, 10$ from top to bottom.} 
    \label{Lshape3DP1to10MeshRefinement}
\end{figure}

\begin{example}[3-d L-shaped Domain]\label{example5}
In this example we consider Problem \eqref{vp_eigen} in a 3-d L-shaped domain $(0,2)^3\setminus(0,1]\times (0,1] \times (0,2)$. The experiments are performed for 3 different values of $p \in \{1.5, 4, 10\}$. The parameter $\epsilon_{M}$ is set to be $10^{-5}$ and the number of loops $K$ in Algorithm \ref{afem} is specified as in Table \ref{tab:Lshape3DEigenvalue}.

Computed first eigenvalues as well as associated numbers of vertices over adaptively generated meshes are listed in Table \ref{tab:Lshape3DEigenvalue}. As before, the sequence of numerical first eigenvalues demonstrates a steadily decreasing trend with the number of vertices increasing. A selection of  adaptive meshes and the computed first eigenfunction over the finest adaptive mesh in the cross section $z<1$ for each $p$ are visualized in Figure \ref{Lshape3DP1to10MeshRefinement} due to the symmetry. In the horizontal direction, we observe refinements near the six corners for small $p$ and then evidently in the interior as well as adjacent to the reentrant corner for large $p$. This behaviour is similar to that in Example \ref{example3} for the 2-d case. In addition, local refinements occur vertically along the concave edge with strong singularities as $p$ increases. However, the refined region shrinks as the computed first eigenfunction vanishes in the most part of $\Omega$. 

\end{example}

\begin{table}[hbt!]
\centering
\caption{Quantitative result for $p \in \{1.5, 4, 10\}$ in Example \ref{example5}: the number of adaptive loops, the number of vertices and the computed first eigenvalue.}
\label{tab:Lshape3DEigenvalue}
\begin{tabular}{c cc cc cc}
    \toprule
    \multirow{2}{*}{$k$}&\multicolumn{2}{c}{$p=1.5$}&
    \multicolumn{2}{c}{$p=4$}&\multicolumn{2}{c}{$p=10$}\\
    \cmidrule(lr){2-7}
     ~ & vertices & $\mu_k$ & vertices & $\mu_k$ & vertices & $\mu_k$ \\
    \midrule
        0 & 1056 & 7.72402 & 1056  &  81.0511 & 1056  & $8.28107\times 10^{3}$ \\

        1 & 1646 & 7.59132 & 1138  & 80.0312 & 1070  & $8.40531\times 10^{3}$ \\

        2 & 3114 & 7.40312 & 1374  & 75.8082 & 1095  & $7.10282\times 10^{3}$ \\

        3 & 5905 & 7.33215 & 1740  &  73.6532 & 1118  & $6.90217\times 10^{3}$ \\

        4 & 9829 & 7.29308 & 2470  & 71.9343 & 1248    &  $6.42960\times 10^{3}$ \\

        5 & 19449& 7.24515 & 4006  & 70.3522 & 1470  & $5.96310\times 10^{3}$ \\

        6 & 36541& 7.22480 & 6197  &  69.3535 & 1775  & $5.66794\times 10^{3}$ \\

        7 & 63085& 7.21281 & 18437  & 68.6378 & 2345 &  $5.37915\times 10^{3}$ \\

        8 &123451& 7.19806 & 24899 & 68.1258 & 3423 & $5.15780\times 10^{3}$ \\

        9 &231816& 7.19169 & 30702 & 67.7104 & 5069 & $4.97964\times 10^{3}$ \\

        10 &      &        & 52381 & 67.4119 & 8046 & $4.86127\times 10^{3}$ \\
	
        11 &      &          & 95228 & 67.1837 & 14011 & $4.76670\times 10^{3}$ \\
        12 &      &          & 162735 & 67.0321 & 22848 & $4.70221\times 10^{3}$ \\
        13 &      &          &        &           & 38977 & $4.63375\times 10^{3}$ \\
    \bottomrule
    \end{tabular}
\end{table}

\begin{figure}[htb!]
    \centering
    \subcaptionbox{$\mathcal{T}_1$ (2503) \label{fig:TorusP2Mesh2}} {\includegraphics[scale=0.35]{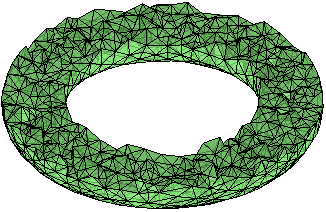}}\hspace{0.5cm}
    \subcaptionbox{$\mathcal{T}_5$ (60227) \label{fig:TorusP2Mesh6}} {\includegraphics[scale=0.35]{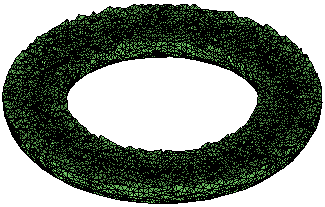}}\hspace{0.5cm}
    \subcaptionbox{$\mathcal{T}_7$ (311225) \label{fig:TorusP2Mesh8}} {\includegraphics[scale=0.35]{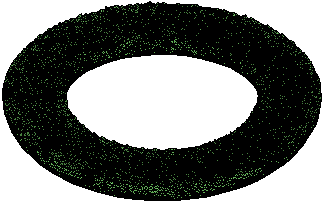}}\hspace{0.5cm}
    \subcaptionbox{$u_7$ \label{fig:TorusP2U8}} {\includegraphics[scale=0.15]{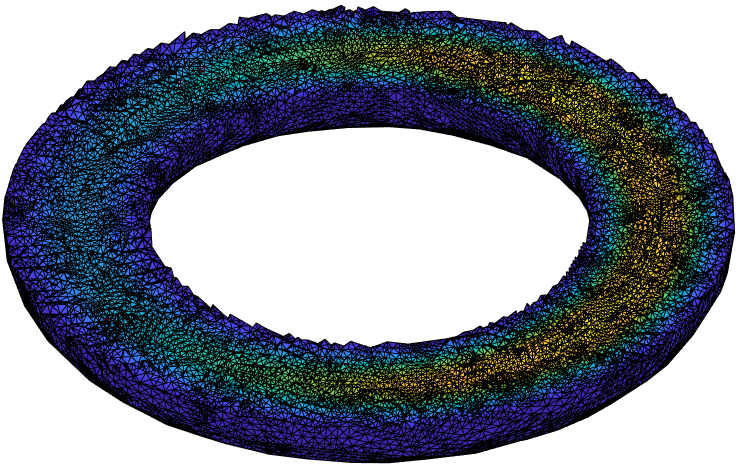}}
    \includegraphics[scale=0.25]{Lshape3DColorbar.png}
\\
    \subcaptionbox{$\mathcal{T}_1$ (1566) \label{fig:TorusP4Mesh2}} {\includegraphics[scale=0.35]{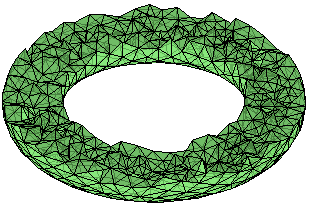}}\hspace{0.5cm}
    \subcaptionbox{$\mathcal{T}_3$ (6783) \label{fig:TorusP4Mesh4}} {\includegraphics[scale=0.35]{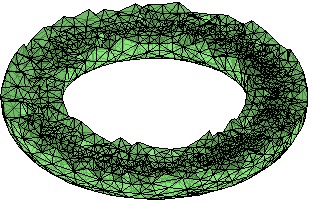}}\hspace{0.5cm}
    \subcaptionbox{$\mathcal{T}_6$ (61965) \label{fig:TorusP4Mesh7}} {\includegraphics[scale=0.35]{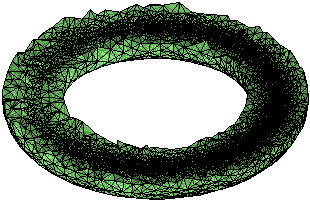}}\hspace{0.5cm}
    \subcaptionbox{$u_6$ \label{fig:TorusP4U7}} {\includegraphics[scale=0.245]{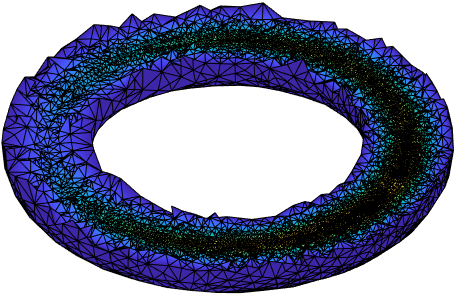}}
    \includegraphics[scale=0.25]{Lshape3DColorbar.png}
    \caption{Adaptive mesh refinement level with the number of vertices over each mesh as well as final computed first eigenfunctions by adaptive refinements in the cross section $z<0$ for $p=2$ (1st row) and $p=4$ (2nd row) in Example \ref{example6}.}
    \label{TorusP1to10MeshRefinement}
\end{figure}

\begin{example}[Torus]\label{example6}
The computational domain of the last example is a torus with major radius $0.8$ and minor radius $0.1$ centered at the origin with $p=2$ and $p=4$. Due to the limited memory of the computer, a maximum number of refinement $K$ is specified as in Table \ref{tab:Torus3DEigenvalue} and the tolerance $\epsilon_{M}$ is set to be $10^{-5}$ and $5\times 10^{-5}$ for $p=2$ and $p=4$ respectively. Table \ref{tab:Torus3DEigenvalue} documents the computed first eigenvalues featuring a gradually steady decline as the adaptive meshes evolves. Figure \ref{TorusP1to10MeshRefinement} displays the adaptive meshes generated by Algorithm \ref{afem} and visualize the computed first eigenfunctions over the finest adaptive meshes in the cross section $z<0$ due to the symmetry. As illustrated in Figure \ref{TorusP1to10MeshRefinement}, the support of the first eigenfunction is shaped like a crescent and becomes narrower with larger $p$. This results in concentrated local refinements. 
\end{example}

\begin{table}[hbt!]
\centering
\caption{Quantitative result for $p \in \{2, 4\}$ in Example \ref{example6}: the number of adaptive loops, the number of vertices and the computed first eigenvalue.}
\label{tab:Torus3DEigenvalue}
\begin{tabular}{c cc cc}
    \toprule
    \multirow{2}{*}{$k$}&\multicolumn{2}{c}{$p=2$}&
    \multicolumn{2}{c}{$p=4$}\\
    \cmidrule(lr){2-5}
     ~ & vertices & $\mu_k$ & vertices & $\mu_k$ \\
    \midrule
        0 & 959 & 158.135 & 959  & $1.14221\times 10^{4}$ \\

        1 & 2503 & 156.570 & 1566  & $1.12114\times 10^{4}$ \\

        2 & 4824 & 154.365 & 3451  & $1.08483\times 10^{4}$ \\

        3 & 11288 & 152.046 & 6783  & $1.04961\times 10^{4}$ \\

        4 & 26008& 150.206 & 13613 & $1.01963\times 10^{4}$ \\

        5 & 60227& 148.546 & 29387 & $9.94798\times 10^{3}$ \\

        6 & 142157& 146.881 & 61965 & $9.75105\times 10^{3}$ \\

        7 & 311225& 145.279 &       &           \\
    \bottomrule
    \end{tabular}
\end{table}

\section{Conclusions}
An adaptive finite element method has been designed to approximate the first eigenvalue of
the $p$-Laplacian operator. We have proved that the sequence
of discrete eigenvalues and discrete eigenfunctions converges to the exact one and the
related eigenset respectively with the help of minimization techniques in derivation of
existence result for nonlinear elliptic equations. In the process, a residual-type error estimator is available and serves in the module ESTIMATE of the adaptive algorithm. The asymptotic behavior of computed first eigenfunctions shows that our adaptive algorithm can capture the singularities as described in the PDE theory. Since the conforming finite element method only provides an upper bound for the first eigenvalue, one natural question is how to yield a lower bound. One possible choice is the nonconforming finite element method, which works for $2$-Laplacian \cite{ArmentanoDuran:2004,cg3,HuHuangLin:2014,Liu:2015,LuoLinXie:2012,YangZhangLin:2010}. In view of this, our future research topic is the study of an adaptive nonconforming method for the first eigenpair of $p$-Laplacian.

\paragraph{Funding} The research of Guanglian Li was partially supported by Hong Kong RGC through General Research Fund (project number: 17317122) and Early Career Scheme (project number: 27301921). The research of Yifeng Xu was partially supported by the National Natural Science Foundation of China (Projects 12250013, 12261160361 and 12271367), the Science and Technology Commission of Shanghai Municipality (Projects 20JC1413800 and 22ZR1445400) and  General Research Fund (Projects KF202318 and KF202468) from Shanghai Normal University. The work of Shengfeng Zhu was partially supported by the National Key Basic Research Program under grant 2022YFA1004402, the National Natural Science Foundation of China (12471377), and the Science and Technology Commission of Shanghai Municipality (Projects 22ZR1421900 and 22DZ2229014).

\end{document}